\documentclass[11pt]{article}
\usepackage{amssymb,amsmath,amsfonts,amsthm,mathtools,color}

\usepackage{mathrsfs}
\usepackage{dsfont}
 \usepackage{bm} 
 
\usepackage[numbers, square, sort&compress]{natbib}
\usepackage{tabularx}

\usepackage{comment} 
\usepackage[title,titletoc,header]{appendix}
\usepackage{graphicx,subcaption}
\usepackage{multirow}
\usepackage{float} 
\usepackage{algorithm}
\usepackage{algpseudocode}

\usepackage{empheq} 

\usepackage{paralist}

\usepackage{tikz}
\usetikzlibrary{arrows,positioning,shapes.geometric}

\graphicspath{ {./figures/} }

\usepackage{indentfirst}
\usepackage{multicol}
\usepackage{booktabs}
\usepackage{url}
\usepackage[outdir=./]{epstopdf} 

\usepackage[shortlabels]{enumitem}

\usepackage{hyperref}
\hypersetup{
    colorlinks=true, 
    linktoc=all,     
    linkcolor=blue,  
}

\numberwithin{equation}{section}

\newtheorem{Theorem}{Theorem}[section]

\newtheorem{prop}[Theorem]{Proposition}
\newtheorem{cor}[Theorem]{Corollary}

\newtheorem{remark}{Remark}

\newtheorem{assumption}{Assumption}

\theoremstyle{definition}
\newtheorem{definition}{Definition}[section]

\theoremstyle{remark}

\def\to{\rightarrow}

\def\Halmos{\mbox{\quad$\square$}}

\def\ms{\medskip}

\def\cF{\mathcal{F}}

\def\d{{\mathrm{d}}}

\def\sE{{\mathbb{E}}}

\def\sN{{\mathbb{N}}}
\def\sP{\mathbb{P}}

\def\sR{{\mathbb R}}

\DeclareMathOperator*{\argmin}{arg\,min}

\newcommand{\lc}
{\mathrel{\raise2pt\hbox{${\mathop<\limits_{\raise1pt\hbox
{\mbox{$\sim$}}}}$}}}

\newcommand{\gc}
{\mathrel{\raise2pt\hbox{${\mathop>\limits_{\raise1pt\hbox{\mbox{$\sim$}}}}$}}}

\newcommand{\ec}
{\mathrel{\raise2pt\hbox{${\mathop=\limits_{\raise1pt\hbox{\mbox{$\sim$}}}}$}}}

\def\bb{\begin{equation}} \def\ee{\end{equation}}
\def\bbn{\begin{equation*}} \def\een{\end{equation*}}

\def\beqn{\begin{eqnarray}}  \def\eqn{\end{eqnarray}}

\def\beqnx{\begin{eqnarray*}} \def\eqnx{\end{eqnarray*}}

\def\bn{\begin{enumerate}} \def\en{\end{enumerate}}

\def\bd{\begin{description}} \def\ed{\end{description}}



\definecolor{DarkGreen}{rgb}{0.2,0.6,0.2}

\def\pink#1{\textcolor{\pink}{#1}}

\definecolor{brilliantrose}{rgb}{1.0, 0.33, 0.64}

\begin{document}

\title{
Decentralized Decision-Making among Heterogeneous Autonomous Vehicles: An  \(\alpha\)-Potential Game Framework
  }

\author{
Anran Hu, Zhexin Wang\\
\scriptsize Department of Industrial Engineering and Operations Research, Columbia University \\
{\scriptsize \texttt{ah4277@columbia.edu, zw3037@columbia.edu}}
\and
Yufei Zhang\\
\scriptsize Department of Mathematics, Imperial College London,  London,  UK \\
{\scriptsize \texttt{yufei.zhang@imperial.ac.uk}} 
\and
Xuan Di*\\
\scriptsize Department of Civil Engineering and Engineering Mechanics, Data Science Institute, Columbia University \\
{\scriptsize \texttt{sharon.di@columbia.edu}}
}

\date{}

\maketitle

\begin{abstract}
We study noncooperative multi-vehicle games among heterogeneous autonomous vehicles, where  each vehicle adopts a decentralized closed-loop policy based on its own state, and optimizes an objective that depends on other vehicles through potentially asymmetric 
interaction weights. We develop an $\alpha$-potential game framework that reduces the computation of an approximate Nash equilibrium (NE) to the minimization of a single auxiliary $\alpha$-potential function. We explicitly construct this $\alpha$-potential, establish the existence of its minimizers, and characterize the equilibrium approximation error $\alpha$ in terms of interaction asymmetry. We further introduce vehicle-specific scaling to reduce the effective interaction asymmetry, thereby tightening the equilibrium approximation and, in important cases, recovering an exact NE despite asymmetric interactions. We also derive social-efficiency guarantees for the potential-selected policies, revealing how the interaction structure shapes worst-case efficiency. Numerical experiments demonstrate the flexibility of the framework in capturing heterogeneous vehicle interactions, collision and obstacle avoidance, lane changing and overtaking under different traffic configurations, and priority-based intersection crossing.
\end{abstract}


\ms

\medskip

\section{Introduction}
\label{sec:intro}

How can one design scalable control strategies for large populations of   autonomous vehicles (AVs) in complex traffic environments? Realistic AV control must account for the heterogeneous behaviors and objectives of interacting vehicles. Different vehicles may compete for limited road space, pursue different desired speeds or routes, and make conflicting decisions when merging, changing lanes, overtaking, or crossing intersections. Capturing such interactions is essential for generating realistic, human-like driving behavior in applications including platooning \citep{gong2016constrained,wei2017dynamic,zhou2017rolling,li2018nonlinear}, lane changing \citep{wang2015game,yu2018human}, merging \citep{wang2024iterative}, unsignalized intersections \citep{jing2024decentralized}, and departure-time choice on transportation networks \citep{wang2023equilibrium}.
These considerations highlight two central challenges: designing decentralized policies that scale to large AV populations and capturing the heterogeneous, self-interested behavior that arises in mixed traffic.

Scalable control of large AV systems requires decentralized decision-making. A full-state feedback policy would require each vehicle to continuously observe, communicate, and process information about an increasing number of surrounding vehicles. As the fleet size grows, this creates substantial sensing, communication, and computational burdens, making real-time implementation increasingly difficult. Decentralized policies offer a more scalable alternative by allowing each AV to make decisions primarily from its own state and locally available information \citep{di2021survey,di2025mfgreview}. By reducing reliance on global state information and centralized coordination, such policies are better suited to large-scale AV deployment. 

Meanwhile, heterogeneous and potentially conflicting objectives in mixed traffic motivate a noncooperative game-theoretic formulation of AV decision-making, in contrast to the fully cooperative control paradigm. A substantial part of the AV-control literature assumes that vehicles coordinate toward a common objective using, for example, consensus-based methods \citep{bailo2018optimal} or model predictive control \citep{daini2024traffic}. Such an assumption can be restrictive when vehicles have different preferences, risk attitudes, priorities, or incentives. Under a game-theoretic formulation, each AV is modeled as a strategic agent optimizing its own objective while responding to the behavior of others. The resulting behavior is described by a Nash equilibrium (NE), at which no vehicle can improve its objective through a unilateral deviation.

\paragraph{Challenges in  equilibrium computation.}

While the game-theoretic formulation better captures self-interested behavior in heterogeneous traffic environments, its main challenge is computational. Computing an NE  in a dynamic game is generally more demanding than solving an optimal control problem: standard procedures often rely on iterative best responses, in which each agent repeatedly solves an optimal control problem against the current strategies of the others \citep{albrecht2024multi}. The resulting computational burden can grow rapidly with the number of vehicles. Although multi-agent reinforcement learning (MARL) provides a flexible numerical alternative, it generally lacks theoretical guarantees of convergence to an NE in large-scale dynamic games. This motivates the search for structural properties of multi-agent interactions that can render equilibrium computation tractable.

Existing work on multi-vehicle games has largely pursued this goal through two theoretical frameworks that exploit particular forms of symmetry. Mean field games (MFGs) \citep{festa2018mean,huang2020game,di2025mfgreview} assume a large population of statistically similar vehicles interacting weakly and symmetrically through their empirical distribution. (Markov) potential games \citep{liu2023potential,yan2025markov}, by contrast, require exact pairwise symmetry in cross-agent interactions: the effect of vehicle \(j\)'s control on vehicle \(i\)'s objective must be matched by the reciprocal effect of vehicle \(i\)'s control on vehicle \(j\)'s objective \cite{liu2023potential}. These structural assumptions facilitate equilibrium analysis, but can be restrictive in heterogeneous traffic environments with asymmetric and strong local interactions. 
See Section~\ref{sec:strcuture_NE} 
for a more detailed discussion of these game-theoretic approaches  and their limitations.

\paragraph{Our work.} 

We develop an $\alpha$-potential game framework for decentralized control of heterogeneous AVs, in which each vehicle adopts a decentralized closed-loop policy based on its own state, and its objective  depends on the states  of other vehicles through heterogeneous interaction weights. Our model accommodates asymmetric and strong local interactions beyond the classical mean field game framework, allowing nearby vehicles to exert nonvanishing influence on one another. This feature is essential for capturing behaviors such as collision avoidance (see Section \ref{sec:application_AV}).

Our main results are as follows.
\begin{itemize}[wide]
\item 
We develop an \(\alpha\)-potential game approach for computing (approximate) NEs in the multi-vehicle game. The framework reduces equilibrium computation to a decentralized control problem with  an auxiliary objective, called an \(\alpha\)-potential function. It accommodates asymmetric interactions and   removes  the exact pairwise symmetry required by the (Markov) potential game approach. We explicitly construct the \(\alpha\)-potential function,
establish the existence of a minimizer, and show that any such minimizer yields an \(\alpha\)-NE of the original game.  We further characterize the approximation parameter \(\alpha\) analytically in terms of the game coefficients, including the length of the decision horizon, the magnitude of the interaction kernel, and the asymmetry of the vehicle interactions (Proposition \ref{prop:alpha_PG}).

\item 

We further sharpen the equilibrium approximation by applying the \(\alpha\)-potential construction to rescaled objectives with vehicle-specific scaling weights (Theorem \ref{thm:rescaled-alpha-potential}). Optimizing these weights to minimize the residual interaction asymmetry yields the tightest approximate-NE guarantee within this family of \(\alpha\)-potential functions. We prove that an optimal rescaling exists under a mild connectivity condition (Theorem \ref{thm:optimal-rescaling}). Moreover, we show that when the interaction structure satisfies a  Kolmogorov cyclic condition, this rescaling eliminates the approximation error entirely, so that any minimizer of the rescaled potential is an exact NE despite asymmetric interactions (Theorem \ref{thm:exact-rescaling}).
In the general case, we develop an efficient bisection algorithm for computing the optimal rescaling, where each iteration solves a linear feasibility problem and the algorithm converges at a geometric rate (Proposition~\ref{prop:bisection-convergence}).
\item 

We derive explicit bounds on the social-efficiency loss of the potential-selected policies, which capture decentralized, self-interested equilibrium behavior, relative to the centralized full-cooperation benchmark (Proposition \ref{prop:efficiency-rescaled-potential} and Corollary \ref{cor:pos-bound}). The bounds make explicit how the interaction structure shapes worst-case social efficiency and thus provide guidance for interaction design.
  Under pairwise symmetric interactions, the social cost of the selected self-interested equilibrium is at most twice the centralized social optimum. More generally, greater interaction asymmetry can  weaken this social-efficiency guarantee through greater dispersion in the rescaling weights, suggesting that limiting excessive asymmetry can improve social-efficiency.

\item 
On the computational side, we develop a policy gradient method with decentralized neural network policies to solve the resulting control problems. The algorithm fits naturally within the centralized training and decentralized execution paradigm widely used in MARL \citep{albrecht2024multi}. Extensive numerical experiments illustrate the modeling flexibility and the scalability of the framework. We compare weak mean field and strong local interaction regimes, capture heterogeneous vehicle interactions and collision avoidance, and reproduce lane changing, overtaking, and priority based intersection crossing across different traffic configurations. The experiments also confirm that  the proposed rescaling procedure can substantially improve the equilibrium approximation.
\end{itemize}

 \paragraph{Organization of the paper.}

The remainder of this paper is organized as follows. 
Section~\ref{sec:lit} reviews existing approaches to AV control and compares them with the proposed framework.
Section~\ref{sec:problem} formulates the multi-vehicle games with  decentralized policies. 
Section~\ref{sec:theory} presents the main theoretical results. It develops the \(\alpha\)-potential game framework and establishes equilibrium approximation guarantees, optimal vehicle-specific rescaling, conditions for recovering exact Nash equilibria, and social-efficiency bounds. It also provides an efficient algorithm for computing the optimal rescaling weights.
Section~\ref{sec:numeric} presents numerical experiments on longitudinal control, obstacle avoidance, lane changing, and intersection crossing, examining the effects of interaction strength, vehicle heterogeneity, and rescaling. 
Finally, Section~\ref{sec:conclude} summarizes the main findings and discusses future research directions.

\section{Related work}
\label{sec:lit}

We begin with a brief overview of cooperative and noncooperative approaches to AV control, and then discuss  existing methods for computing equilibria in noncooperative multi-vehicle games.

\subsection{Cooperative and noncooperative frameworks for AV decision-making}

AV decision-making has been formulated from both cooperative   and noncooperative   perspectives.

In   cooperative-control, 
AVs coordinate toward a common system-level objective, with applications including platooning \citep{han2022strategic}, lane changing \citep{fu2023cooperative,heshami2024towards,chen2025game}, merging \citep{fu2025regional}, and intersection management \citep{liu2025cooperative}.  Decentralized (or distributed) control methods have been developed to reduce the communication and computational burden  by relying more heavily on local information \citep{shen2022distributed}. While these approaches decentralize information and computation, they generally retain the premise that vehicles pursue a common or coordinated objective.

Noncooperative game-theoretic formulations instead model AVs as self-interested decision makers,  and    characterize their strategic interactions through Nash equilibria.
Such approaches have been applied to  lane changing and car-following \citep{talebpour2015modeling,wang2015game,zhang2024stackelberg,yao2025personalized}, longitudinal velocity control \citep{huang2020game,wang2024multi},  platooning and convoy control \citep{legal2023platooning,liu2024decentralized,jond2022differential}, routing \citep{shou2022markov}, merging
\citep{wei2022merging,zhang2025integrated}, and intersection interactions \citep{jing2024decentralized,huang2024robust,shu2025decision}. However, computing Nash equilibria in large-scale dynamic multi-vehicle games remains challenging unless additional structure in the game can be exploited.

For cooperative  and noncooperative   settings, MARL provides  as a flexible numerical framework for multi-vehicle  control  in dynamic and potentially partially observed environments \citep{wang2024multi,liu2024cooperative,guo2024heuristic,hua2025multi}. However, ensuring 
convergence of an MARL algorithm to an NE
in a  dynamic game    requires additional structural assumptions. 

\subsection{Tractable structures for     equilibrium computation
in multi-vehicle games}
\label{sec:strcuture_NE}

The challenges of computing NEs motivate the development of theoretical frameworks that exploit interaction structures to make equilibrium analysis and computation tractable. Two prominent approaches widely used in traffic modeling are mean field games and  (Markov) potential games. Although conceptually different, both rely on particular forms of symmetry that can be restrictive in heterogeneous AV settings.

The MFG framework considers the large-population limit of similar AVs interacting weakly and symmetrically through their empirical distribution \citep{festa2018mean,huang2020game,di2025mfgreview}.  These   assumptions can be restrictive for traffic applications:  weak interactions are less suited to capturing strong local effects such as collision avoidance, while homogeneity overlooks heterogeneous objectives and asymmetric interactions that arise in mixed traffic; see Remark 
\ref{rmk:interaction}.

(Markov) potential games, by contrast, address the finite-vehicle game directly (see e.g., \citealt{liu2023potential,jing2025decentralized,yan2025markov,chen2026deep}). Following \citet{monderer1996potential}, they posit a global potential function whose change under a player's unilateral deviation exactly matches the corresponding change in the deviating player’s objective, thereby reducing equilibrium computation to a single optimization problem of the potential function. However, such an exact potential need not exist in dynamic multi-vehicle games. Even when restricting to time-dependent controls without explicit state feedback, potential structure requires exact pairwise symmetry in cross-vehicle interactions \citep{liu2023potential}. This requirement can be restrictive for heterogeneous AV systems with asymmetric interactions.

Recent work   relaxes exact potential structure while retaining the computational advantage of optimizing a single auxiliary objective; see, e.g., \citealt{candogan2013near,varga2024upper}. In particular, \cite{guo2025markov} introduces the notion of an \(\alpha\)-potential game, in which the change in the potential under a unilateral deviation approximates the corresponding change in a player’s objective up to an error \(\alpha\). The case \(\alpha=0\) recovers an exact potential game, while minimizers of an \(\alpha\)-potential yield approximate NEs. This framework has subsequently been developed for broad classes of dynamic games (\citealt{guo2025alpha,guo2025towards,guo2025distributed}). For general stochastic differential games, \cite{guo2025alpha} constructs an \(\alpha\)-potential explicitly from the game coefficients and quantifies the approximation error in terms of asymmetry in the interaction structure.

We build on this framework for multi-vehicle games with decentralized closed-loop policies and asymmetric interactions. Exploiting the decentralized state dynamics yields a simpler potential construction, while introducing a vehicle-specific scaling further reduces effective interaction asymmetry and sharpens the equilibrium approximation.

\section{Methodology}
\label{sec:problem}

This section formulates a  dynamic Markov game with decentralized policies  for AV control.

\subsection{Problem setup}
Consider a finite-player stochastic differential game with decentralized Markov policies, defined as follows: 
$T>0$ is a fixed  time horizon, 
   $[N]=\{1,\ldots, N\}$,  $N\in \mathbb N$,
    is a finite set of
players, and for all $i\in [N]$, 
$ A_i\subset \sR^{k}$ is a nonempty compact set representing   player $i$'s action set,
and 
     $\Pi_i\subset C([0,T]\times \sR^{d}; A_i)$ is the set  of player $i$'s admissible policies that  are locally Lipschitz continuous in the spatial variable.
 We denote by 
$\Pi=\prod_{i\in [N]}\Pi_i$   the set of all admissible policy profiles of all players, and by   
$\Pi_{-i}=\prod_{j\not = i}\Pi_j$   the set of admissible policy profiles of all players except player $i$. We denote by 
 $\phi=(\phi_i)_{i\in [N]}$
 and 
 $\phi_{-i}=(\phi_j)_{j\in [N]\setminus\{i\}}$
 a generic element of $\Pi$ and $\Pi_{-i}$, respectively.  

Each player’s state evolves according to a controlled diffusion, determined by their chosen policy. The player then minimizes a cost functional over the set of admissible policies. 
More precisely, 
for any given $\phi=(\phi_i)_{i\in [N]}\in \Pi$,
player $i$'s state dynamics $X^{\phi_i}_i$ 
is governed by  the following  dynamics:
for all $t\in [0,T]$,
\begin{equation}
\label{eq:state_i}
\begin{aligned}
    \d  X_{i,t} =
& b_i(X_{i,t},\phi_i(t,X_{i,t}))
\d t 
+\sigma_i  
(X_{i,t},\phi_i(t,X_{i,t}))
\d W_t, 
\quad 
X_{i,0}=\xi_i,
\end{aligned}
\end{equation}
where $\xi_i\in \sR^d$ is a given initial state, 
 $b_i: \sR^{d}\times A_i\to \sR^{d}$
 and 
 $\sigma_i: \sR^{d}\times A_i\to \sR^{d\times m}$
 are   Lipschitz continuous functions,   and  
 $W:\Omega\times [0,T]\to \sR^m$ is  an  $m$-dimensional Brownian motion  on 
a complete  probability space 
$(\Omega,\cF,\sP)$. Given $\phi_{-i}\in \Pi_{-i}$,
 player $i$ determines their optimal policy  by minimizing the following objective function   over  $\Pi_i$:
\begin{equation}
    \begin{aligned}
        \label{eq:cost_i}
      J_i(\phi) &= \sE \Biggl[\int_0^T 
  \Biggl(
  f_i\big(
  X^{\phi_i}_{i,t},
  \phi_{i}(t,X^{\phi_i}_{i,t})\big) 
   +\sum_{j\not =i}\lambda_{ij}K(X^{\phi_i}_{i,t}-X^{\phi_j}_{j,t})
  \Biggr)\d t  + g_i(X^{\phi_i}_{i,T})\Biggl],
    \end{aligned}
\end{equation}
where 
$\lambda_{ij}\ge 0$ is a given constant, 
$f_i:\sR^d\times   A_i\to \sR$ and 
$K:\sR^d\to \sR$ 
are  given    running costs,
and 
 $g_i:\sR^d\to\sR$ is  a given  terminal cost.
 We assume that
 $f_i$ and $g_i$ are continuous and  at most of quadratic growth, 
 and 
 $K$ is a bounded continuous function satisfying  
   $K(x)=K(-x)$ for all $x\in \sR^d$  and $K\not\equiv 0$.

 We  characterize the rational behavior of players in the Markov game \eqref{eq:state_i}–\eqref{eq:cost_i} using the notion of an $\epsilon$-Nash equilibrium given below.

      \begin{definition}
     \label{def:NE}
          For any $\epsilon\ge 0$, a policy profile 
          $\bar{\phi}=(\bar{\phi}_i)_{i\in [N]}\in \Pi$  an $\epsilon$-Nash equilibrium (NE) of the
game \eqref{eq:state_i}–\eqref{eq:cost_i} if
$$
J_i(\bar{\phi})\le 
J_i(\phi_i,\bar{\phi}_{-i})+\epsilon, \quad \forall  \phi_i\in \Pi_i, i\in [N].
$$
If $\epsilon=0$,
then $\tilde \phi$ is called an NE. 
     \end{definition}

According to Definition \ref{def:NE}, a joint  policy profile $\bar \phi$ is an $\epsilon$-NE if no player can improve her objective by more than $\epsilon$ through any unilateral deviation. Note that although player $i$'s objective in \eqref{eq:cost_i} depends on the joint states of all players, player $i$ optimizes \eqref{eq:cost_i} over decentralized policies that depend only on her own private state $X_i$.  

Throughout this paper, we impose suitable compactness conditions on the class of admissible policies.

\begin{assumption}
\label{assum:compact_policy}
   For all $i\in [N]$, $\Pi_i$ is compact 
    in the compact-open topology, i.e., 
    $\Pi_i$ is  equicontinuous on every compact subset of $[0,T]\times \sR^d$.
\end{assumption}

Assumption \ref{assum:compact_policy}   reflects physical constraints on each player's control. 
It holds, for example, when $\Pi_i$ consists of neural networks with Lipschitz activation functions and bounded weights.
 Since $\Pi_i$ is generally nonconvex, standard equilibrium-existence results do not directly apply to ensure the existence of an (approximate) Nash equilibrium in the game   \eqref{eq:state_i}-\eqref{eq:cost_i}.

\subsection{Specialization to AV control} 
\label{sec:application_AV}

The game \eqref{eq:state_i}-\eqref{eq:cost_i} encompasses various dynamic games    in AV control. In this setting, each player represents a vehicle. The dynamics   \eqref{eq:state_i} can be taken as controlling velocity: 
\begin{equation}
\label{eq:AV_1}
\d  x_{i,t} =
 \phi_i(t,x_{i,t} )
\d t, \quad t\in [0,T],
\end{equation}
where $x_{i,t}$  represents 
  the position   of the $i$-th vehicle at time $t$, 
  and the control 
  $v_{i,t}\coloneqq \phi_i(t,x_{i,t})$ represents 
  the   velocity    of the $i$-th vehicle at time $t$. 
Alternatively, 
one can take the dynamics   \eqref{eq:state_i} as a double integrator with controlled acceleration:
for all $t\in [0,T]$,
\begin{align}
\label{eq:AV_2}
    \begin{split}
\d  x_{i,t} &=
 v_{i,t} \d t,
\quad 
\d  v_{i,t}  =
 \phi_i(t,x_{i,t}, v_{i,t})
\d t +\sigma_i \d W_t, 
    \end{split}
\end{align}
where $x_{i,t}$ and $v_{i,t}$  represent
  the position  and velocity   of the $i$-th vehicle at time $t$, respectively,
  and  the control 
  $u_{i,t}\coloneqq \phi_i(t,x_{i,t},v_{i,t})$ represents 
  the   acceleration    of the $i$-th vehicle at time $t$. 
  
  Vehicle $i$ determines its optimal route by minimizing the objective  \eqref{eq:cost_i}, where the terminal cost $  g_i$  specifies vehicle $ i$'s preferred target,    the running  cost $ f_i$ penalizes control effort and deviations from the target speed, while also enforcing obstacle-avoidance behavior for vehicle $ i $, and  the   kernel $ K $ together with the weights $ {(\lambda_{ij})}_{j\not = i} $  characterizes how the spatial distribution of other vehicles  influences vehicle $i$'s preferred route.
Typically, the kernel 
$K$
 decreases as the distance between two vehicles increases, modeling congestion-averse behavior.

\begin{remark}[Local and asymmetric interactions.]
\label{rmk:interaction}

    Our model accommodates both strong local interactions and asymmetric interactions among AVs, thereby going beyond   common structural assumptions in existing game-theoretic approaches. 
First, unlike the   MFG framework \citep{huang2020game}, interactions need not be  only through the empirical  distribution. Nearby vehicles may exert nonvanishing influence on one another, which is essential for modeling behaviors such as collision avoidance; see Section~\ref{sec:numerical_interaction}. 
Second, the interaction weights need not be pairwise symmetric, i.e., we allow $\lambda_{ij}\neq \lambda_{ji}$, in contrast to classical potential-game formulations that rely on symmetric cross-agent interactions.

To illustrate the strong-local-interaction feature, following \citep{oelschlager1985law}, consider
\begin{equation}
\label{eq:interaction_kernel}
\lambda_{ij}=N^{\beta-1},
\qquad
K(z)=\rho\bigl(N^{\beta/d}|z|\bigr),
\qquad z\in\sR^d,
\end{equation}
where $\rho:[0,\infty)\to[0,\infty)$ satisfies
$\lim_{|z|\to\infty}\rho(z)=0$, and $\beta\in[0,1]$ determines how the interaction strength scales with $N$.
When $\beta=0$, each player interacts with all others with strength $\mathcal O(1/N)$. This is the classical weak-interaction regime underlying mean-field approximations \citep{huang2006Largepopulation,lasry2007mean}. As $N\to\infty$, the influence of individual players vanishes and interactions are captured through the aggregate population distribution, as in the MFG formulation for traffic flow in \citep{huang2020game}. In contrast, when $\beta=1$, the decay of $\rho$ implies that each player interacts only with neighbors at distance $\mathcal O(N^{-1/d})$, but with interaction strength of order one. This strong-local-interaction regime allows nearby vehicles to have a non-negligible effect on individual decisions and is therefore well suited to modeling collision avoidance. Such interactions fall outside the classical   mean-field regime, while remaining within the scope of our framework.
\end{remark}

\section{Main   Results}
\label{sec:theory}

\subsection{$\alpha$-potential function and approximate NEs}

This section analyzes approximate NEs of the game \eqref{eq:state_i}–\eqref{eq:cost_i} by adapting the $\alpha$-potential game framework developed in \citep{guo2025alpha} to the present setting with decentralized closed-loop controls. The $\alpha$-potential game  framework reduces  the challenging problem of identifying (approximate) NEs into a single optimization problem of the associated $\alpha$-potential function.  
This connection is made precise in the following proposition, which follows directly  from  \citep[Proposition 2.1]{guo2025alpha}.

 \begin{prop}
\label{prop:alpha_pg}
  Suppose that 
  the game \eqref{eq:state_i}-\eqref{eq:cost_i} 
  is an  $\alpha$-potential game for some $\alpha\ge 0$, in the sense that there exists 
  $\Phi:\Pi\to \sR$,
  called an $\alpha$-potential function, 
  such that 
for all $i\in [N]$, $\phi_{-i}\in \Pi_{-i}$, and
$\phi_i,\phi'_i\in \Pi_i$, 
\begin{equation}
\label{eq:alpha_pg}
|(J_i(\phi'_i,\phi_{-i})
-J_i(\phi_i,\phi_{-i}))
-(\Phi(\phi'_i,\phi_{-i})
-\Phi(\phi_i,\phi_{-i}))|\le \alpha.
\end{equation}
If $\bar \phi\in \Pi$
satisfies 
$\Phi(\bar \phi)\in \argmin_{\phi\in \Pi}\Phi(\phi)$,
then $\bar \phi$ is an $\alpha$-NE of the game \eqref{eq:state_i}-\eqref{eq:cost_i}.
 
 \end{prop}

Intuitively, a function $\Phi$ is an $\alpha$-potential function for the game \eqref{eq:state_i}-\eqref{eq:cost_i}  if, whenever a player unilaterally deviates from her strategy, the resulting change in her objective differs from the corresponding change in $\Phi$ by at most $\alpha$. Once constructed, $\Phi$ serves as an auxiliary central objective for constructing approximate NEs of the game. 
In what follows, we analytically construct the $\alpha$-potential function $\Phi$. Importantly, $\Phi$ is not merely an aggregation of the players' objectives; rather, it preserves the strategic structure of the underlying non-cooperative game and captures the self-interested nature of the players.

The $\alpha$-potential game framework generalizes existing (Markov) potential game approaches to autonomous driving (see, e.g., \citet{liu2023potential,yan2025markov}) by allowing $\alpha>0$. This relaxation is particularly important for dynamic games among AVs with asymmetric interactions, which generally do not admit an exact potential function but can be characterized as $\alpha$-potential games for a suitable $\alpha>0$.  
Specifically, 
  the following proposition
  constructs analytically    an $\alpha$-potential function
for the game \eqref{eq:state_i}-\eqref{eq:cost_i} and quantifies the associated $\alpha$,
in terms of the   decision  horizon, the magnitude of the interaction kernel $K$  and   the degree of asymmetry in the interaction weights $(\lambda_{ij})_{i,j \in [N]} $.

\begin{prop}

\label{prop:alpha_PG}
 
Define $\Phi:\Pi\to \sR$ such that for all $\phi\in \Pi$,
\begin{equation}\label{eq:potential_fun_symmetric}
     \Phi(\phi) \coloneqq \sE \left[\int_0^T F(X^{\phi}_t,\phi(t, X^{\phi}_t))\d t + G(X^{\phi}_T)\right],
\end{equation}
where
$X^\phi=(X^{\phi_i}_i)_{i\in [N]}$
is the joint state process,
$\phi(t, X^{\phi}_t)
=
(\phi_i(t, X^{\phi_i}_{i,t}))_{i\in [N]}$ 
is the joint control process, and 
$F:\sR^{dN}\times A\to \sR$
and $G:\sR^{dN} \to \sR$
are given by
\begin{align*}
\begin{split}
F(x,a) &\coloneqq \sum_{i=1}^N f_i(x_i,a_i) +
\sum_{1\le i<j\le N} \frac{\lambda_{ij}+\lambda_{ji}}{2} K(x_i-x_j),
\quad 
G(x)    \coloneqq \sum_{i=1}^N g_i(x_i). 
\end{split}
\end{align*}
Then $\Phi$ is an $\alpha$-potential function of the game  \eqref{eq:state_i}-\eqref{eq:cost_i}  with 
\begin{equation}
\label{eq:alpha_bound}
    \alpha \leq T\|K\|_{L^\infty}\max_{i\in [N]}\sum_{i\neq j}|\lambda_{ij}-\lambda_{ji}|.
\end{equation}
Hence any minimizer of $\Phi$ over $\Pi$ is an $\alpha$-NE of   the game  \eqref{eq:state_i}-\eqref{eq:cost_i}. 
\end{prop}

Proposition \ref{prop:alpha_PG} highlights the importance of allowing a general $\alpha\ge 0$ in \eqref{eq:alpha_pg} to accommodate games with asymmetric interactions.
Indeed, by Proposition \ref{prop:alpha_PG},
     $\alpha=0$ (so that $\Phi$ is an exact potential function)  if and only if players'   interactions  is symmetric, i.e., $
    \lambda_{ij}=\lambda_{ji}$ for all $i,j\in [N]$.
    Consequently, the (discrete-time) potential game frameworks in \citet{liu2023potential,yan2025markov} are restricted to AVs with symmetric interactions, whereas our framework accommodates general asymmetric interactions among AVs.
  We refer to \citep{guo2025alpha,guo2025distributed} for large-population games with asymmetric interactions and small $\alpha$.

\begin{remark}[Construction of $\alpha$-potential functions]
   \label{eq:potential}
The construction of the $\alpha$-potential function in Proposition \ref{prop:alpha_PG} differs from those in \citep{guo2025alpha, guo2025distributed} in two key aspects:
(i) Proposition \ref{prop:alpha_PG} considers the Markov game \eqref{eq:state_i}–\eqref{eq:cost_i} with Lipschitz (closed-loop) policies, whereas \citep{guo2025alpha, guo2025distributed} focus on \emph{open-loop} controls without explicit state feedback;
(ii) Proposition \ref{prop:alpha_PG} exploits the independence of state dynamics   to construct the  $\alpha$-potential function \eqref{eq:potential_fun_symmetric} 
without introducing additional sensitivity equations as auxiliary states, as required in \citep{guo2025alpha, guo2025distributed}. This more tractable $\alpha$-potential function   enables   more efficient computation of approximate NEs. 

\end{remark}

Leveraging Proposition \ref{prop:alpha_PG}, the next theorem provides additional flexibility in constructing an $\alpha$-potential function  by appropriately scaling each player’s objective, thereby yielding a tighter approximate NE guarantee for games with asymmetric interactions. 

\begin{Theorem}\label{thm:rescaled-alpha-potential}
    Let $p=(p_i)_{i\in [N]}\in (0,\infty)^N$, and define 
    $\Phi^p:\Pi\to \sR$ such that for all $\phi\in \Pi$,
 \begin{equation}
 \label{eq:potential_fun_symmetric_p}
\Phi^p(\phi) \coloneqq \mathbb E\bigg[ \int_0^T F^p(X^{\phi}_t,\phi(t, X^{\phi}_t))\d t + G^p(X^{\phi}_T) \bigg],
\end{equation}
where
$X^\phi$
and 
$\phi(t, X^{\phi}_t)$ are defined as in Proposition 
\ref{prop:alpha_PG},
and 
$F^p:\sR^{dN}\times A\to \sR$
and $G^p:\sR^{dN} \to \sR$
are given by
\begin{align*}
\begin{split}
F^p(x,a) &\coloneqq \sum_{i=1}^N p_if_i(x_i,a_i) +
\sum_{1\le i<j\le N} \frac{p_i\lambda_{ij}+p_j\lambda_{ji}}{2} K(x_i-x_j),
\quad 
G^p(x)    \coloneqq \sum_{i=1}^N p_ig_i(x_i). 
\end{split}
\end{align*}
Then    $\Phi^p$ admits a minimizer over $\Pi$, each of which    is an $\alpha(p)$-NE of   the game  \eqref{eq:state_i}-\eqref{eq:cost_i}, with 
$$
\alpha(p)\le T\|K\|_{L^\infty}  \max_{i\in[N]} \frac{1}{p_i} \sum_{j\ne i} |p_i\lambda_{ij}-p_j\lambda_{ji}|.
$$

\end{Theorem}

Compared with Proposition \ref{prop:alpha_PG}, Theorem \ref{thm:rescaled-alpha-potential} provides a more general construction of the central objective for approximate NEs by introducing a scaling vector $p$ that can mitigate the asymmetry of the interaction matrix. The function $\Phi$ in \eqref{eq:potential_fun_symmetric} corresponds to the special case $p_i\equiv 1$ for all $i\in[N]$. By optimizing $p$ to minimize the residual interaction asymmetry, one can construct $\Phi^p$ whose minimizer achieves a tighter NE gap than that obtained in Proposition \ref{prop:alpha_PG}. As will be shown in Section \ref{sec:optimization_alpha}, this approach can even yield an exact NE for games with asymmetric interactions.

\subsection{Optimal  scaling weights}
\label{sec:optimization_alpha}

To obtain the tightest equilibrium approximation guarantee from
Theorem~\ref{thm:rescaled-alpha-potential}, we minimize the approximation bound
over $p\in(0,\infty)^N$. Define
\begin{equation}
\label{eq:optimal-rescaling-value}
    r^\star
     \coloneqq 
    \inf_{p\in(0,\infty)^N}
    \max_{i\in[N]}
    \frac{1}{p_i}
    \sum_{j\ne i}
    |p_i\lambda_{ij}-p_j\lambda_{ji}|,
\end{equation}
so that the smallest bound obtainable from Theorem~\ref{thm:rescaled-alpha-potential} is
\[
    \alpha^\star
    =
    T\|K\|_{L^\infty}r^\star.
\]
Since the optimization in \eqref{eq:optimal-rescaling-value} is over the open set $(0,\infty)^N$, attainment of the infimum is not immediate. We next show that, under a mild connectivity condition on the interaction structure, the problem can be equivalently reformulated over a compact set.

Associate with the interaction matrix $\Lambda=(\lambda_{ij})_{i,j\in[N]}$ the directed interaction graph \[ \mathcal G_\lambda=([N],E_\lambda), \] where there is a directed edge $j\to i$ if and only if $\lambda_{ij}>0$. Thus an edge $j\to i$ means that player $j$ influences the objective of player $i$.

\begin{assumption}[Componentwise strong connectivity] \label{ass:componentwise-connected} 
The directed   graph $\mathcal G_\lambda$ can be decomposed into a disjoint union of strongly connected subgraphs\footnote{A directed graph is strongly connected if, for any two vertices $i$ and $j$, there exists a directed path from $i$ to $j$.} with vertex sets $S_1,\ldots,S_m$. That is, \[ [N]=S_1\cup\cdots\cup S_m, \] each subgraph induced by $S_k$ is strongly connected, and there are no edges between distinct subgraphs (i.e., 
$ 
\lambda_{ij}=\lambda_{ji}=0$, for all  $i\in S_k,\ j\in S_\ell,\ k\ne\ell$). 
\end{assumption}

\begin{remark}\label{rmk:example-componentwise-connected}
Assumption~\ref{ass:componentwise-connected} is satisfied, for example, when the entire interaction graph $\mathcal G_\lambda$ is strongly connected. It is also satisfied when the interaction pattern has symmetric support, in the sense that $\lambda_{ij}>0$ if and only if $\lambda_{ji}>0$ for all $i\ne j$; in this case, the sets $S_1,\ldots,S_m$ can be taken as the connected components of the corresponding undirected interaction graph. In particular, the assumption allows interaction matrices with zero entries and disconnected groups of players.
\end{remark}

This next theorem shows that under Assumption~\ref{ass:componentwise-connected}, the optimization problem \eqref{eq:optimal-rescaling-value} admits an equivalent compact formulation, which in particular guarantees the existence of an optimal rescaling vector.

\begin{Theorem}[Existence of an optimal rescaling]
\label{thm:optimal-rescaling}
Suppose Assumption~\ref{ass:componentwise-connected} holds. Then 
\eqref{eq:optimal-rescaling-value} is equivalent to
\begin{equation}
\label{eq:optimal-rescaling-compact}
\begin{aligned}
    \min_{p\in[0,\infty)^N,\ r\in[0,r_0]}\quad & r\\
    \text{s.t.}\quad
    &
    \sum_{j\ne i}
    |p_i\lambda_{ij}-p_j\lambda_{ji}|
    \le r p_i,
    \qquad i\in[N],\\
    &
    \sum_{i\in S_k}p_i=|S_k|,
    \qquad k=1,\ldots,m,
\end{aligned}
\tag{P}
\end{equation}
where $r_0 \coloneqq 
\max_{i\in[N]}
\sum_{j\ne i}
|\lambda_{ij}-\lambda_{ji}|$. Consequently, the infimum in
\eqref{eq:optimal-rescaling-value} is obtained by some
$p^\star\in(0,\infty)^N$.
\end{Theorem}

We next identify the condition under which the optimal rescaling reduces the approximation error to zero, and hence minimizing \eqref{eq:potential_fun_symmetric_p} yields an   NE.

\begin{Theorem} 
\label{thm:exact-rescaling}
The following
statements are equivalent:
\begin{enumerate}[label=(\roman*)]
    \item $r^\star=0$;
    \item there exists $p=(p_1,\ldots,p_N)\in(0,\infty)^N$ such that
    \begin{equation}
    \label{eq:detailed-balance-lambda}
        p_i\lambda_{ij}
        =
        p_j\lambda_{ji},
        \qquad i,j\in[N];
    \end{equation}
    \item 
    $(\lambda_{ij})_{i,j\in [N]}$ satisfies 
    the following Kolmogorov cyclic condition: 
    \begin{itemize}
    \item The zero pattern is symmetric: for all $i\ne j$,
    $\lambda_{ij}=0$ if and only if $\lambda_{ji}=0$.
    \item For every finite cycle $i_1,i_2,\ldots,i_m,i_{m+1}$ with
    $i_{m+1}=i_1$,
    \begin{equation}
    \label{eq:kolmogorov-cycle}
        \prod_{k=1}^m \lambda_{i_k i_{k+1}}
        =
        \prod_{k=1}^m \lambda_{i_{k+1} i_k}.
    \end{equation}
\end{itemize}
    
\end{enumerate}
If any of the above conditions holds, 
let $p^*=(p^*_i)_{i\in [N]}$ satisfy \eqref{eq:detailed-balance-lambda}. Then 
any minimizer of $\Phi^{p^\star}$ is an NE 
  of   the game  \eqref{eq:state_i}-\eqref{eq:cost_i}.
\end{Theorem}

\begin{remark}[Rank-one interaction]
\label{rmk:rank-one-interactions}
 If 
$\lambda_{ij}=\gamma_i\tau_j$ for some $\gamma_i,\tau_i>0$,
then $(\lambda_{ij})_{i,j\in [N]}$ 
satisfies the  Kolmogorov cyclic condition. 
Indeed, the weights $p$ satisfying \eqref{eq:detailed-balance-lambda} is given by  
$p_i=\frac{\tau_i}{\gamma_i}$, $i\in[N]$,
up to a common positive scaling factor.
\end{remark}

We next turn to the computation of an optimal rescaling vector in the general case. For any fixed $r\in[0,r_0]$, the feasibility problem in \eqref{eq:optimal-rescaling-compact} is convex in $p$ and can be written as a linear feasibility problem by introducing auxiliary variables for the absolute-value terms. Moreover, feasibility is monotone in $r$. Therefore, an optimal rescaling parameter can be computed by bisection over $r\in[0,r_0]$.

\begin{algorithm}[H]
\caption{Bisection for Optimizing Rescaling Parameters}
\label{alg:bisection-rescaling}
\begin{algorithmic}[1]
\Require Interaction matrix $(\lambda_{ij})_{i,j\in[N]}$, components
$S_1,\ldots,S_m$, number of iterations $n$
\State Set $r_{\rm low}=0,\
    r_{\rm high}
    =
    r_0
     \coloneqq 
    \max_{i\in[N]}
    \sum_{j\ne i}|\lambda_{ij}-\lambda_{ji}|$.
\State Initialize $p=(1,\ldots,1)$.
\For{$k=0,\ldots,n-1$}
    \State Set
    $r_{\rm mid} \coloneqq \frac{r_{\rm low}+r_{\rm high}}{2}$.
    \State Check feasibility of
    \[
    \begin{aligned}
        \text{find}\quad &\tilde p\in[0,\infty)^N\\
        \text{s.t.}\quad&
        \sum_{j\ne i}
        |\tilde p_i\lambda_{ij}-\tilde p_j\lambda_{ji}|
        \le r_{\rm mid}\tilde p_i,
        \qquad i\in[N],\\
        &
        \sum_{i\in S_\ell}\tilde p_i=|S_\ell|,
        \qquad \ell=1,\ldots,m.
    \end{aligned}
    \]
    \If{the feasibility problem is feasible}
        \State Set $r_{\rm high}=r_{\rm mid}$ and $p=\tilde p$.
    \Else
        \State Set $r_{\rm low}=r_{\rm mid}$.
    \EndIf
\EndFor
\State \Return $p$ and $r=r_{\rm high}$.
\end{algorithmic}
\end{algorithm}

The bisection procedure is globally convergent, with an explicit geometric error bound, as stated below.

\begin{prop}[Convergence of the bisection algorithm]
\label{prop:bisection-convergence}
Suppose Assumption~\ref{ass:componentwise-connected} holds. Let
$(p^n,r^n)$ denote the output of Algorithm~\ref{alg:bisection-rescaling}
after $n$ bisection iterations. Then
\[
    0\le r^n-r^\star\le \frac{r_0}{2^n},
    \qquad
    0\le \alpha(p^n)-\alpha^\star
    \le T\|K\|_{L^\infty}\frac{r_0}{2^n}.
\]
Consequently, $r^n$ converges to $r^\star$ and $\alpha(p^n)$ converges to
$\alpha^\star$ at the   rate $O(2^{-n})$.
\end{prop}

\paragraph{Computational complexity.} 
Each iteration of Algorithm~\ref{alg:bisection-rescaling} requires solving one linear feasibility problem. After introducing auxiliary variables to linearize the absolute-value terms, the problem is a linear program with $O(N+M_\lambda)$ variables and constraints, where $M_\lambda$ denotes the number of interacting pairs. 
Hence each iteration can be solved in time polynomial in $q \coloneqq N+M_\lambda$ using standard linear-programming algorithms; see, e.g., \citealt{karmarkar1984new,vaidya1989speeding,cohen2021solving}. In particular, $q=O(N^2)$ in the dense-interaction case, while $q=O(N)$ when each player interacts with a uniformly bounded number of neighbors. Together with Proposition~\ref{prop:bisection-convergence}, an $\varepsilon$-accurate value of $r^\star$ requires only $O(\log(r_0/\varepsilon))$ such linear-programming feasibility solves.

\subsection{Social efficiency: game-theoretic versus centralized control}

After characterizing approximate Nash equilibria, we turn to their social-efficiency implications. The minimizers of the rescaled $\alpha$-potential functions  given in \eqref{eq:potential_fun_symmetric_p} capture the self-interested behavior of heterogeneous vehicles in mixed traffic and therefore need not minimize the aggregate social cost achieved under full cooperation. We quantify the resulting efficiency loss by comparing these equilibrium outcomes with the centralized social optimum, thereby measuring how much social welfare is sacrificed to accommodate decentralized, self-interested decision-making.
The resulting bounds also provide guidance for the design of vehicle interaction structures.

To quantify the efficiency loss relative to full cooperation, define 
 the social cost
 for any joint policy profile $\phi\in\Pi$ by
\begin{equation} \label{eq:social-cost} 
\operatorname{SC}(\phi)  \coloneqq \sum_{i=1}^N J_i(\phi), 
\end{equation} 
and the optimal social cost achieved in the full-cooperation setting   by
\begin{equation} \label{eq:social-optimal-cost} 
\operatorname{SC}^{\rm opt}  \coloneqq \inf_{\phi\in\Pi}\operatorname{SC}(\phi). 
\end{equation}
In the sequel, we assume without loss of generality that
$\operatorname{SC}^{\rm opt}>0$, and   measure the social efficiency of a policy profile $\phi$
through the following inefficiency ratio.

\begin{definition}[Social inefficiency] 
\label{def:social-inefficiency} 
For any $\phi\in\Pi$, define 
\begin{equation} 
\label{eq:social-inefficiency-ratio}
\operatorname{InEff}(\phi)  \coloneqq \frac{\operatorname{SC}(\phi)} {\operatorname{SC}^{\rm opt}}. 
\end{equation}
Thus $\operatorname{InEff}(\phi)\ge 1$,
with values closer to one indicating   greater social efficiency. 
\end{definition}

\begin{remark}
    The inefficiency ratio $\operatorname{InEff}(\phi)$ extends the
Price of Stability (PoS) perspective to arbitrary policy profiles, including
approximate Nash equilibria. Indeed, letting $\mathcal E$ denote the
set of NEs,
\begin{equation}
\label{eq:PoS}
\operatorname{PoS}
\coloneqq
\frac{\displaystyle\inf_{\phi\in\mathcal E}\operatorname{SC}(\phi)}
     {\operatorname{SC}^{\rm opt}}
=
\inf_{\phi\in\mathcal E}\operatorname{InEff}(\phi),
\end{equation}
so   PoS measures the social inefficiency of the most
efficient NE.

We work with $\operatorname{InEff}(\phi)$ because, under
Assumption~\ref{assum:compact_policy}, an exact NE need
not exist when the policy class is nonconvex. Even when exact equilibria
exist, identifying the socially most efficient one can be computationally
difficult, making the price of stability hard to evaluate. In contrast,
$\operatorname{InEff}(\phi)$ assesses the social performance of any
individual policy profile and therefore applies directly to the
approximate equilibria constructed in this paper. When $\alpha^*=0$ as in Theorem \ref{thm:exact-rescaling},
the potential-selected profile is an exact NE, and hence
its inefficiency ratio provides an upper bound on $\operatorname{PoS}$.
\end{remark}

  In the sequel, we derive an analytical bound on the inefficiency ratio of the selected equilibrium policy in terms of the asymmetry of the interaction weights \((\lambda_{ij})_{i,j\in[N]}\) and the rescaling vector \(p\). The key observation is that both the rescaled \(\alpha\)-potential and the social cost decompose into the same player-specific and pairwise interaction cost components, but with different coefficients. This reduces the analysis of social inefficiency in the dynamic game to a comparison of these coefficients and makes explicit how interaction asymmetry and the choice of rescaling weights affect social welfare.

To see it, define the non-interaction component of player $i$'s cost by
\[
C_i(\phi)  \coloneqq \mathbb E\left[ \int_0^T f_i\big(X^{\phi_i}_{i,t},\phi_i(t,X^{\phi_i}_{i,t})\big)\,dt + g_i(X^{\phi_i}_{i,T}) \right], 
\]
and   the pairwise interaction term by 
\[
\mathcal K_{ij}(\phi)  \coloneqq \mathbb E\left[ \int_0^T K\big(X^{\phi_i}_{i,t}-X^{\phi_j}_{j,t}\big)\,dt \right].
\]
Then 
\begin{equation} 
\label{eq:social-cost-decomposition}
\operatorname{SC}(\phi) = \sum_{i=1}^N C_i(\phi) + \sum_{1\le i<j\le N} (\lambda_{ij}+\lambda_{ji})\mathcal K_{ij}(\phi),
\end{equation}
whereas the rescaled potential associated with $p$ can be written as
\begin{equation} 
\label{eq:rescaled-potential-decomposition}
\Phi^p(\phi) = \sum_{i=1}^N p_i C_i(\phi) + \frac12 \sum_{1\le i<j\le N} (p_i\lambda_{ij}+p_j\lambda_{ji})\mathcal K_{ij}(\phi).
\end{equation}

The following proposition provides an explicit analytical bound on the inefficiency ratio of the equilibrium policy selected by minimizing the rescaled potential \eqref{eq:rescaled-potential-decomposition}, in terms of the interaction weights \((\lambda_{ij})_{i,j\in[N]}\) and the rescaling vector \(p\).

\begin{prop}[Social inefficiency guarantee] 
\label{prop:efficiency-rescaled-potential}
Suppose that $C_i(\phi)\ge 0$ and $\mathcal K_{ij}(\phi)\ge 0$ for all admissible $\phi$, $i\in[N]$, and $i<j$, and  $\operatorname{SC}^{\rm opt}>0$.
Then for all  $p\in(0,\infty)^N$
and 
$ 
\phi^p\in\arg\min_{\phi\in\Pi}\Phi^p(\phi)
$,\footnotemark
\begin{equation}
\label{eq:efficiency-guarantee}
    \operatorname{InEff}(\phi^p)
    \le
    \frac{\overline{\gamma}(p)}
         {\underline{\gamma}(p)},
\end{equation}
where
\begin{equation}
\label{eq:gamma-lower}
\underline{\gamma}(p)
 \coloneqq 
\min\left\{
    \min_{i\in[N]}p_i,\,
    \min_{\substack{1\le i<j\le N\\
                    \lambda_{ij}+\lambda_{ji}>0}}
    \frac{p_i\lambda_{ij}+p_j\lambda_{ji}}
         {2(\lambda_{ij}+\lambda_{ji})}
\right\},
\end{equation}
and
\begin{equation}
\label{eq:gamma-upper}
\overline{\gamma}(p)
 \coloneqq 
\max\left\{
    \max_{i\in[N]}p_i,\,
    \max_{\substack{1\le i<j\le N\\
                    \lambda_{ij}+\lambda_{ji}>0}}
    \frac{p_i\lambda_{ij}+p_j\lambda_{ji}}
         {2(\lambda_{ij}+\lambda_{ji})}
\right\}.
\end{equation}

\end{prop}
 
\footnotetext{
The bound \eqref{eq:efficiency-guarantee}   is nontrivial when $\lambda_{ij}+\lambda_{ji}>0$ for at least one pair $i<j$. When $\lambda_{ij}=0$ for all $i\ne j$, the game  reduces to $N$ independent control problems, and is excluded.}

When   $\alpha^*=0$ as in Theorem \ref{thm:exact-rescaling}, the bound in \eqref{eq:efficiency-guarantee} can be sharpened to yield a bound on the PoS (see \eqref{eq:PoS}).

\begin{cor} \label{cor:pos-bound}
 Suppose that there exists 
$p=(p_1,\ldots,p_N)\in(0,\infty)^N$ satisfying
\eqref{eq:detailed-balance-lambda}. Assume further that the
hypotheses of Proposition~\ref{prop:efficiency-rescaled-potential}
are satisfied.  Then 
   
   \begin{equation} 
   \label{eq:pos-kolmogorov-simple} 
   \operatorname{PoS} \le 
\operatorname{InEff}(\phi^p)\le\frac{\displaystyle\max_{i\in[N]}p_i} {\displaystyle\min_{1\le i<j\le N}\frac{p_ip_j}{p_i+p_j}} 
   \le 2\frac{\max_{i\in[N]}p_i}{\min_{i\in[N]}p_i}. 
   \end{equation}
   If, in addition, the interaction coefficients are pairwise  symmetric, i.e., $\lambda_{ij}=\lambda_{ji}$ for all $i,j\in[N]$,   then 
   $\operatorname{PoS}\le
    \operatorname{InEff}(\phi^{\mathbf 1})\le 2$. 
\end{cor}

\begin{remark}[Interaction asymmetry and social efficiency] 

The bound in \eqref{eq:pos-kolmogorov-simple} quantifies the
worst-case social-efficiency loss associated with decentralized,
self-interested behavior relative to the fully cooperative benchmark.
In the symmetric-interaction case,  the bound 
$\operatorname{InEff}(\phi^{\mathbf 1})\le 2$ implies that, the social cost of the selected
self-interested equilibrium is at most twice the centralized social
optimum.
More generally, under the Kolmogorov cyclic condition,
\eqref{eq:detailed-balance-lambda} implies that, for every interacting
pair,
\[
\frac{p_i}{p_j}
=
\frac{\lambda_{ji}}{\lambda_{ij}}.
\]
Hence asymmetry in the interaction coefficients is reflected directly
in the relative magnitudes of the rescaling weights. A larger
discrepancy between $\lambda_{ij}$ and $\lambda_{ji}$ induces greater
dispersion in the weights $p_i$, which can weaken the worst-case
social-efficiency guarantee through the ratio
$\max_i p_i/\min_i p_i$ in \eqref{eq:pos-kolmogorov-simple}.

These bounds also provide a design insight: interaction structures that avoid excessive asymmetry can lead to stronger worst-case social-efficiency guarantees.

\end{remark}

\section{Numerical Experiments} 
\label{sec:numeric}

In this section, we aim to evaluate the performance of our proposed model in a variety of traffic settings, summarized in Table~\ref{table: sim_scenarios}. 

\begin{table}[htbp]
\centering
\caption{Numerical simulation scenarios.}
\label{table: sim_scenarios}
\footnotesize
\setlength{\tabcolsep}{3pt}
\renewcommand{\arraystretch}{1.15}
\begin{tabularx}{\textwidth}{
    @{}>{\raggedright\arraybackslash}p{0.19\textwidth}
    c
    >{\raggedright\arraybackslash}p{0.12\textwidth}
    >{\raggedright\arraybackslash}X
    >{\raggedright\arraybackslash}p{0.22\textwidth}@{}
}
\toprule
Scenario (Section) & Dim. & Control & Experiment & Goal \\
\midrule
Interaction regimes (\ref{sec:numerical_interaction})
& 1D & Velocity; acceleration
& 10 vehicles; weak versus strong interactions ($\beta=0$ versus $\beta=1$) (Figures \ref{fig:first_order_1D} and \ref{fig:second_order_1D})
& Compare separation and collision avoidance \\

Scalability (\ref{sec:numerical_interaction})
& 1D & Velocity
& 30 vehicles; strong interactions  ($\beta=1$) (Figure \ref{fig:30_vehicle})
& Assess scalability to a larger number of vehicles \\

Obstacle avoidance (\ref{sec:obstacle})
& 2D & Acceleration
& 10 vehicles; no obstacle, or a circular obstacle of radius $0.1$ or $0.5$ (Figure \ref{fig:small_obstacle})
& Assess trajectory responses to obstacle size \\

Rescaling (\ref{sec:different_vehicle})
& 1D & Velocity; acceleration
& 9 vehicles with asymmetric interactions (Figure \ref{fig:different_vehicle}); exploitability before and after rescaling (Table \ref{tab:exploitability})
& Assess size-dependent equilibrium behavior; Evaluate improvement in Nash approximation \\

Lane changing (\ref{sec:lane_changing})
& 2D & Velocity
& 3 vehicles, 2 lanes; two initial layouts with identical costs and target speeds (Figure \ref{fig:comparison})
& Compare lane changes and overtaking behavior \\

Intersection crossing (\ref{sec:intersection})
& 2D & Acceleration
& 8 vehicles; same initial conditions, with versus without major--minor-road differentiation (Figure \ref{fig:crossing_comparison})
& Assess road-priority effects on yielding and crossing \\
\bottomrule
\end{tabularx}

\smallskip
\begin{minipage}{\textwidth}
\scriptsize
Dim.\ denotes spatial dimension.
\end{minipage}
\end{table}

\noindent\textbf{Algorithm.}
To minimize the $\alpha$-potential function  \eqref{eq:potential_fun_symmetric} and compute  an $\alpha$-NE, we  parametrize player $i$'s  policy  using  a sufficiently expressive neural network\footnote{The neural network architecture is feedforward fully connected with one or two hidden layers.}  $\phi_i^{\theta_i}:[0,T]\times\sR^{d}\to A_i$ with parameters $\theta_i\in\sR^{L_i}$,
and discretize the state dynamics \eqref{eq:state_i} 
on time grid $\pi_P=\{0=t_0<\dots<t_P=T\}$ for some $P\in \sN$: for any $i\in[N]$, $X^{\theta_i}_{i,0}=\xi_i$, and for any $\ell=0,\dots,P-1$, 
\begin{equation}\label{eq:param-dynamics-disc}
\begin{aligned}
        X_{i,t_{\ell+1}}^{\theta_i} &= X_{i,t_{\ell}}^{\theta_i} + b_i( X_{i,t_{\ell}}^{\theta_i}, \phi_i^{\theta_i}(t_\ell, X_{i,t_{\ell}}^{\theta_i}))\Delta_\ell      + \sigma_i( X_{i,t_{\ell}}^{\theta_i}, \phi_i^{\theta_i}(t_\ell, X_{i,t_{\ell}}^{\theta_i}))\Delta W_\ell,
\end{aligned}
\end{equation}
where 
$\Delta_\ell =t_{\ell+1}-t_\ell$, and $ 
\Delta W_\ell = W_{t_{\ell+1}}-W_{t_\ell}$. 
We minimize the following discretized $\alpha$-potential function 
\begin{equation}\label{eq:param-potential}
  \Phi(\theta) \coloneqq \sE\left[\sum_{\ell=0}^{P-1} F\!\big(X_{t_\ell}^\theta,\phi^\theta({t_\ell},X_{t_\ell}^\theta)\big)\d t+G\big(X_T^\theta\big)\right]
\end{equation}
over 
 $\theta=(\theta_i)_{i\in [N]}$.
This is achieved by using
the  Adam method, where at each iteration, 
 \eqref{eq:param-potential} is   approximated by  
\begin{equation}\label{eq:param-potential-approx}
    \begin{aligned}
            \Phi_M(\theta)&=\frac{1}{M}\sum_{m=1}^M \bigg[\sum_{\ell=0}^{P-1}F(X_{t_\ell}^{\theta,(m)},\phi^{\theta}(t_\ell,X_{t_\ell}^{\theta,(m)})) \Delta_\ell
             +G(X_{T}^{\theta,(m)})\bigg],
    \end{aligned}
\end{equation}
where 
$\left(X^{\theta,(m)}\right)_{ m\in[M]}=\left(X_{i }^{\theta_i,(m)} \right)_{ i\in[N] ,m\in[M]}$ are $M$ independent trajectories of   \eqref{eq:param-dynamics-disc}.

\subsection{Longitudinal control}

This set of experiments use a fleet of 10 AVs on a long highway (1D) or its 2D extension with a common start position ($x_{0,i}=-1$ in 1D cases and $(-1,-1)$ in 2D cases) and a target position ($z_{i}=1$ in 1D cases and $(1,1)$ in 2D cases).  

\subsubsection{Interaction regimes in cost function}
\label{sec:numerical_interaction}
Here we examine the impact of $\beta$ in \eqref{eq:interaction_kernel} for both velocity-control and acceleration-control models. Specifically, for each model, we test   two   choices  $\beta=0$ and $\beta=1$, which characterize the weak-interaction regime and strong-interaction regime, respectively. We consider a 10-player game in AV control and assume that each player has a one-dimensional state, which follows from \eqref{eq:AV_1} under velocity control and \eqref{eq:AV_2} under acceleration control with a uniform noise level $\sigma_i=0.1$. Each vehicle $i$  
minimizes its objective function defined in \eqref{eq:cost_i} with $f_i(x,a_i)=0.1a_i^2$ and  $g_{i}(x)=10|x-z_{i}|^2$, where $z_{i}=1$ is the target terminal position for vehicle $i$.  
Interaction kernel $\lambda_{ij}$ and $K$ are defined in \eqref{eq:interaction_kernel} with
$\rho:\sR^d\to [0,\infty)$  given by
$
\rho (z)={(|z|^2+1)^{-1} }.
$ 

\begin{figure}[htbp]
    \centering
    \includegraphics[width=0.473 \textwidth]{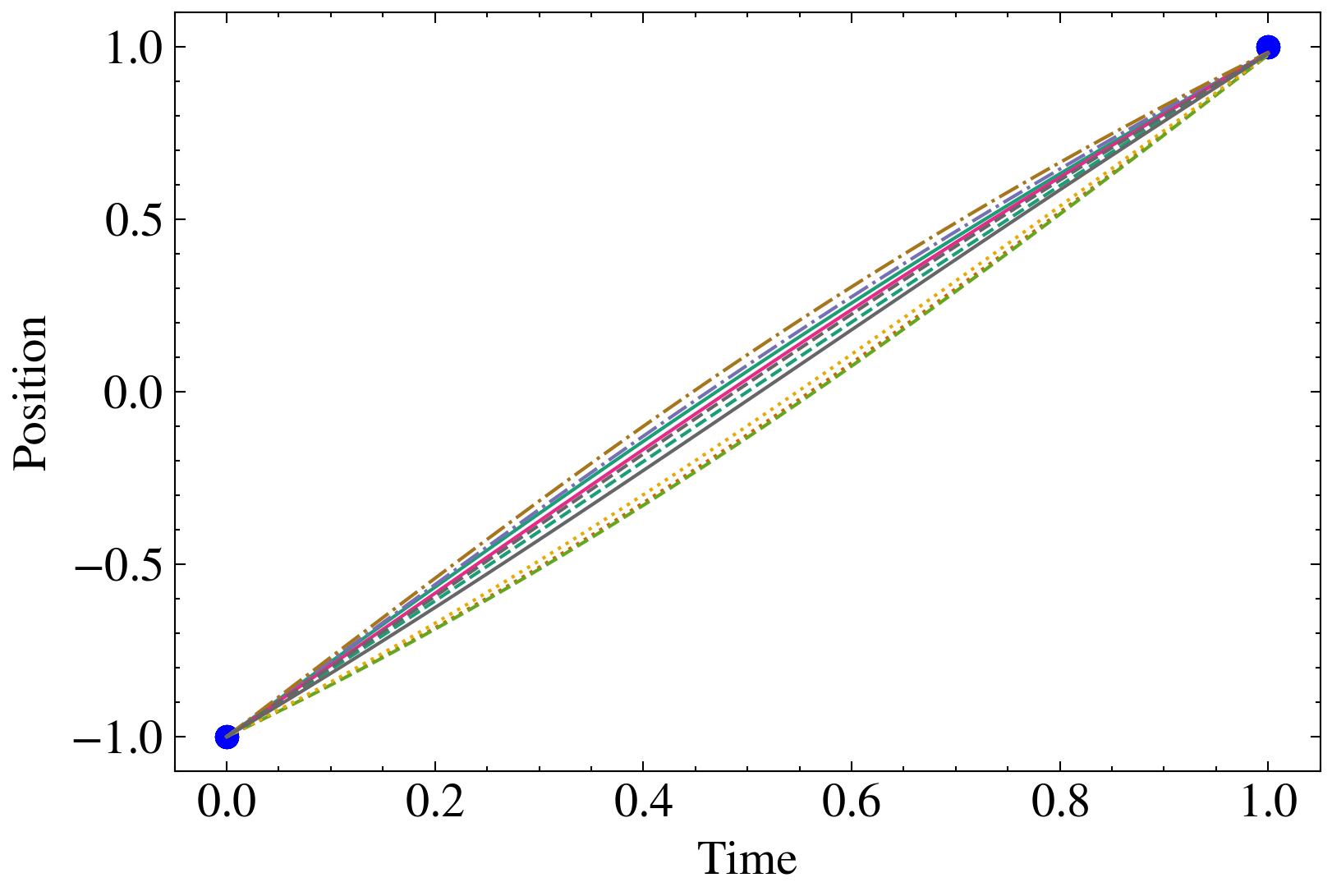}
    \quad 
    \includegraphics[width=0.47 \textwidth]{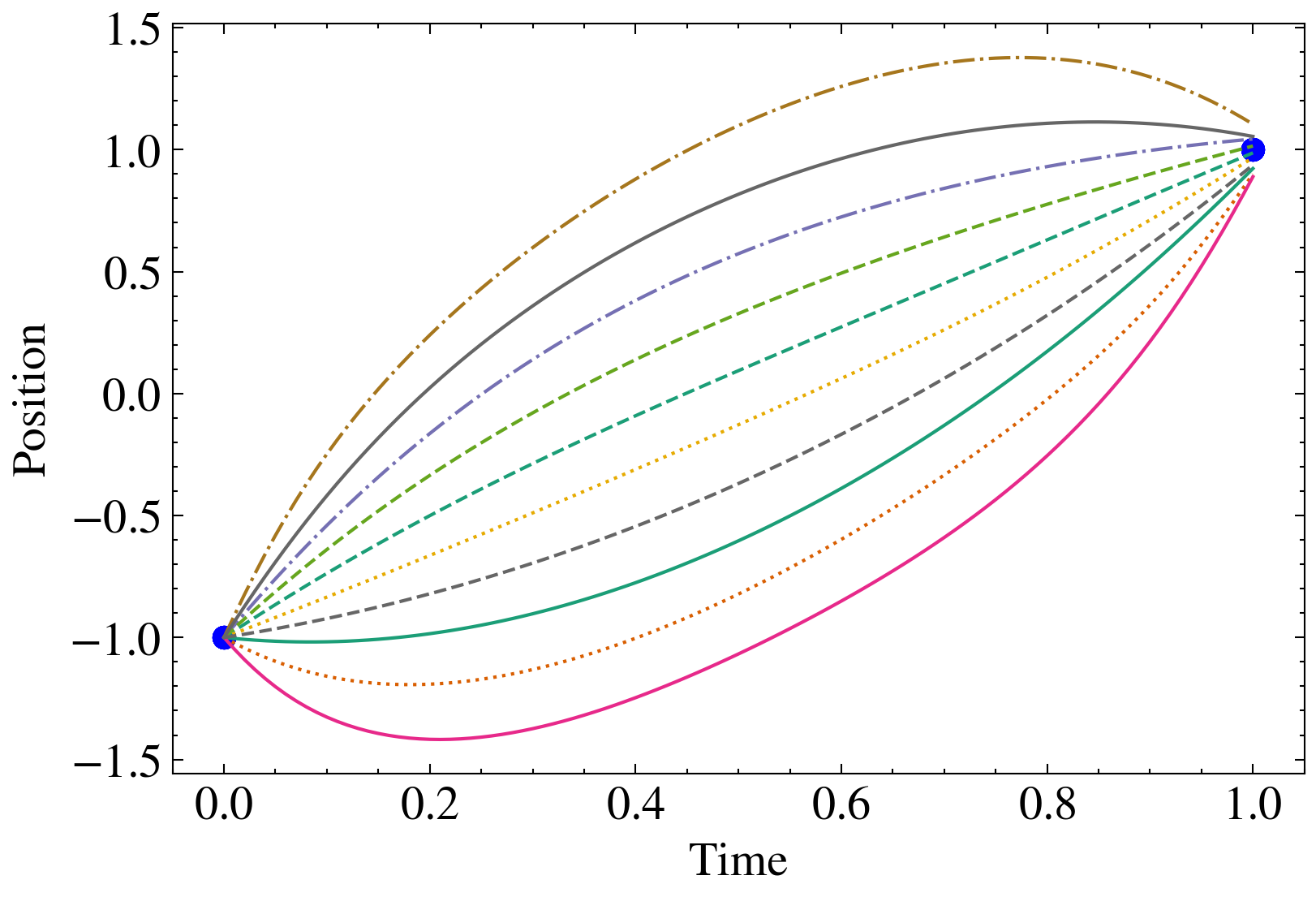}
    \caption{Vehicle trajectories of the control-velocity model for $\beta=0$ (left) and $\beta=1$ (right)}
    \label{fig:first_order_1D}
\end{figure}

\begin{figure}[htbp]
    \centering
    \includegraphics[width=0.47 \textwidth]{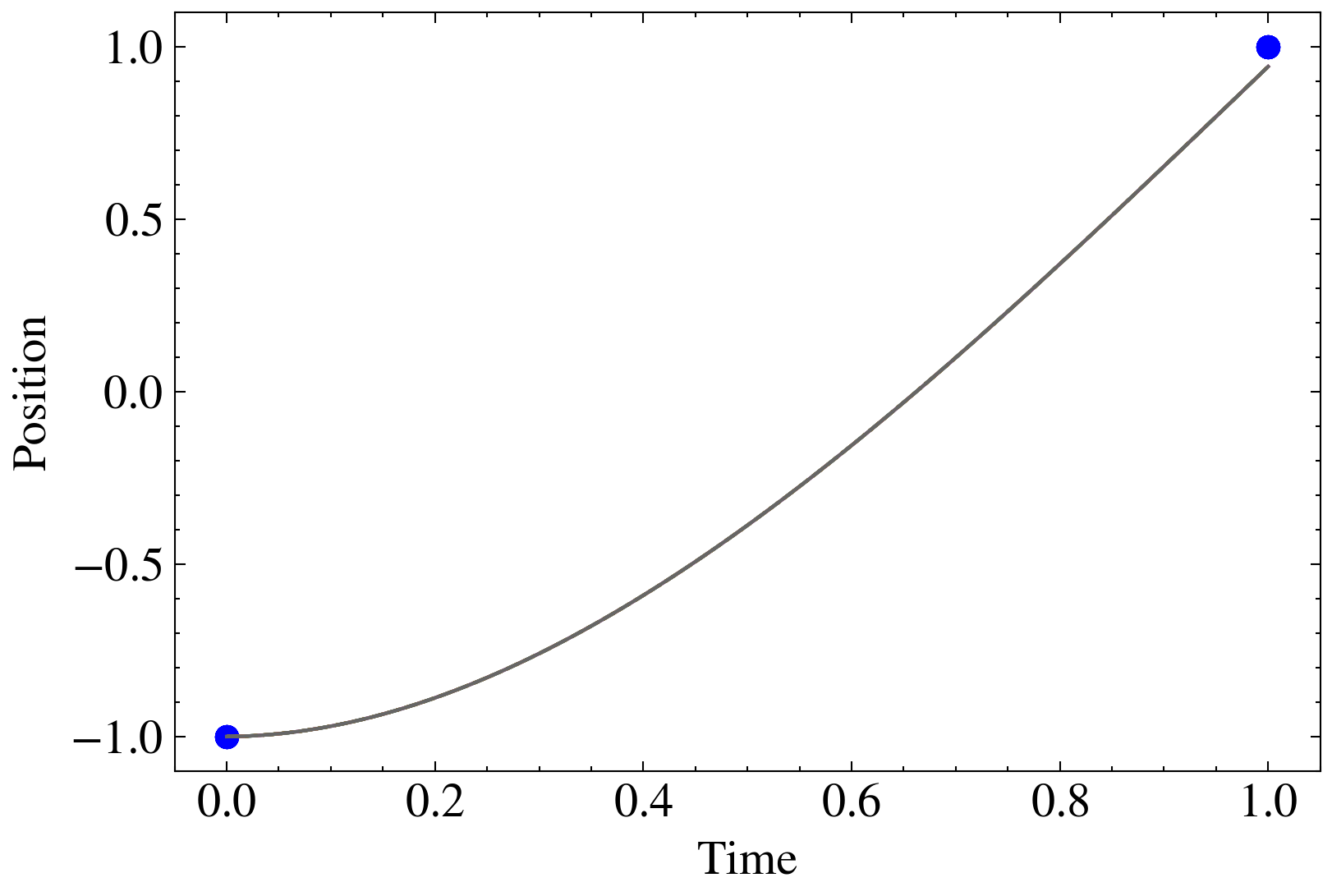}
    \quad 
    \includegraphics[width=0.47 \textwidth]{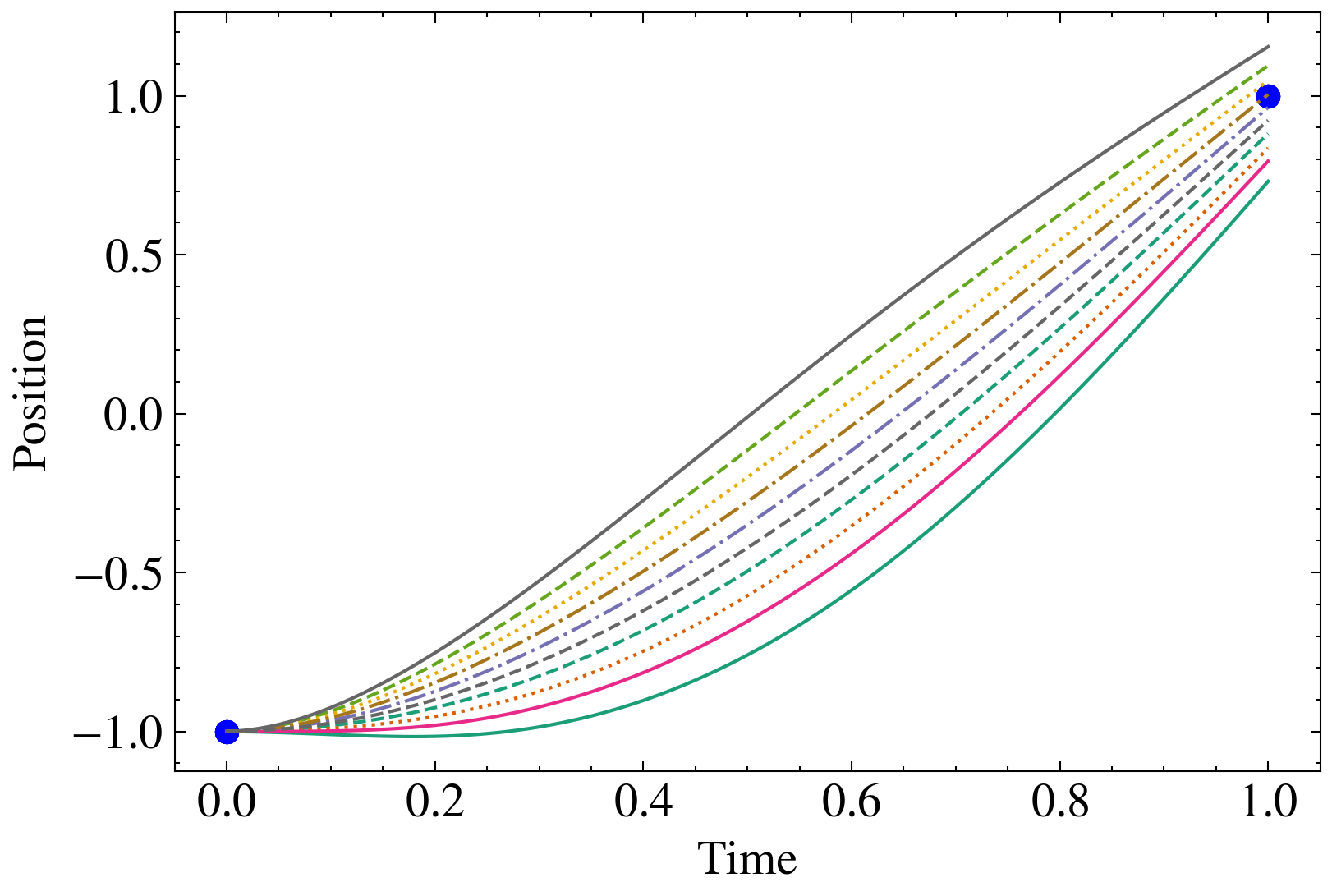}
    \caption{Vehicle trajectories of the control-acceleration model   for $\beta=0$ (left) and $\beta=1$ (right)}
    \label{fig:second_order_1D}
\end{figure}

We observe in Figures~\ref{fig:first_order_1D} and \ref{fig:second_order_1D} that, in the weak-interaction regime ($\beta=0$), both the velocity- and acceleration-control models yield near-identical controls and hence similar trajectories, reflecting the symmetric MFG structure in which agents adopt the same feedback policy. However, in the strong-interaction regime ($\beta=1$), vehicles deliberately choose different controls to give each other space, and their paths diverge. This shows that the strong-interaction regime ensures collision avoidance more effectively. 

In addition, one can see that under velocity-control, vehicles frequently detour (even briefly backtrack) to avoid others yet still reach the target accurately, whereas under acceleration-control, trajectories are smoother but exhibit dispersed terminal positions, indicating greater difficulty in hitting the target precisely. This happens because velocity control gives agents direct authority over position with cheap late corrections, so they can prioritize separation and fix position near the end; acceleration control introduces inertia, making sharp turns and last-minute corrections costly, so the optimizer smooths motion and tolerates small terminal errors.

\begin{figure}[htbp]
    \centering
    \begin{subfigure}{0.47\textwidth}
        \centering
        \includegraphics[width=\linewidth]{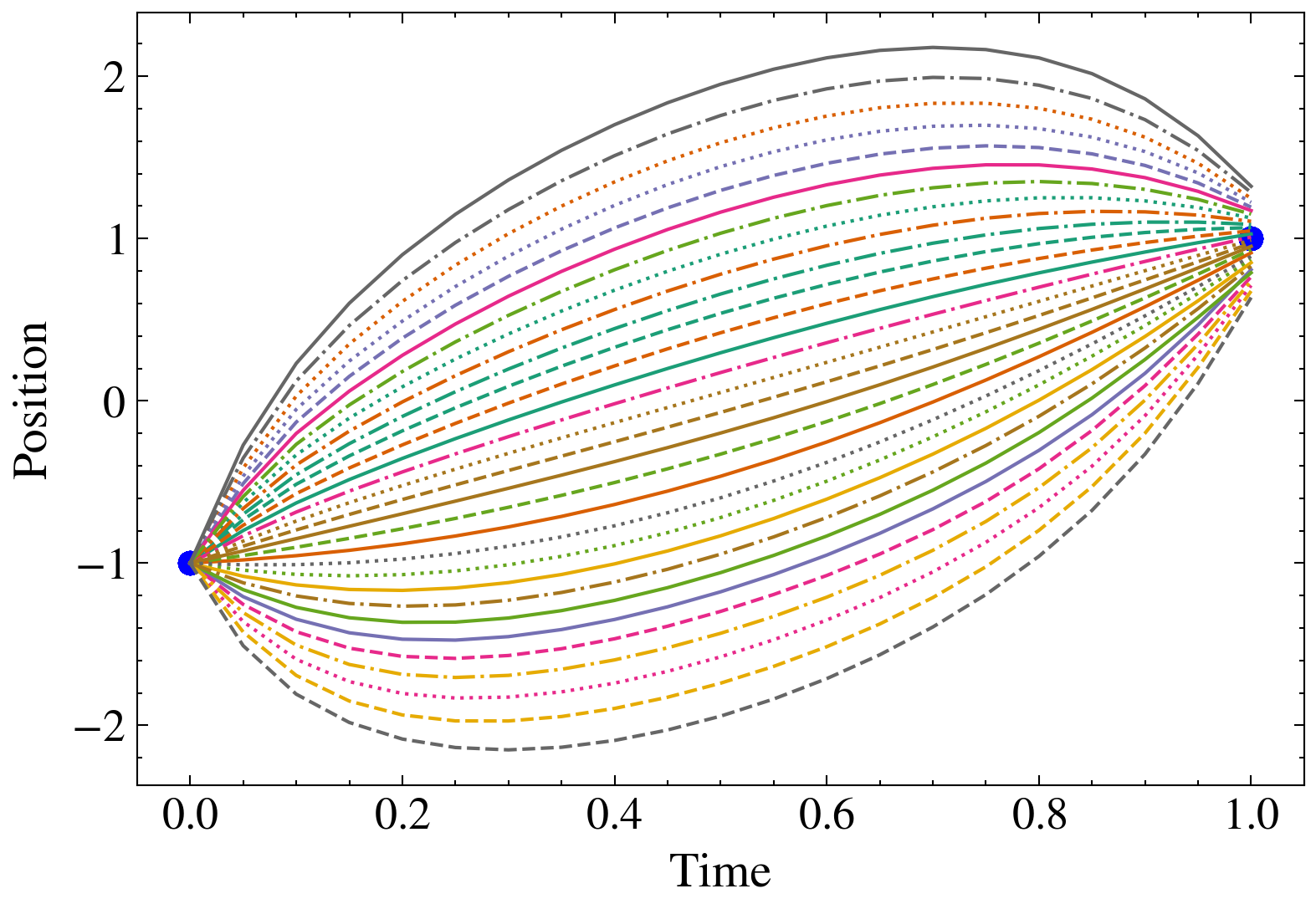}
        \caption{Trajectories of 30 vehicles}
        \label{fig:first}
    \end{subfigure}
    \hfill
    \begin{subfigure}{0.472\textwidth}
        \centering
        \includegraphics[width=\linewidth]{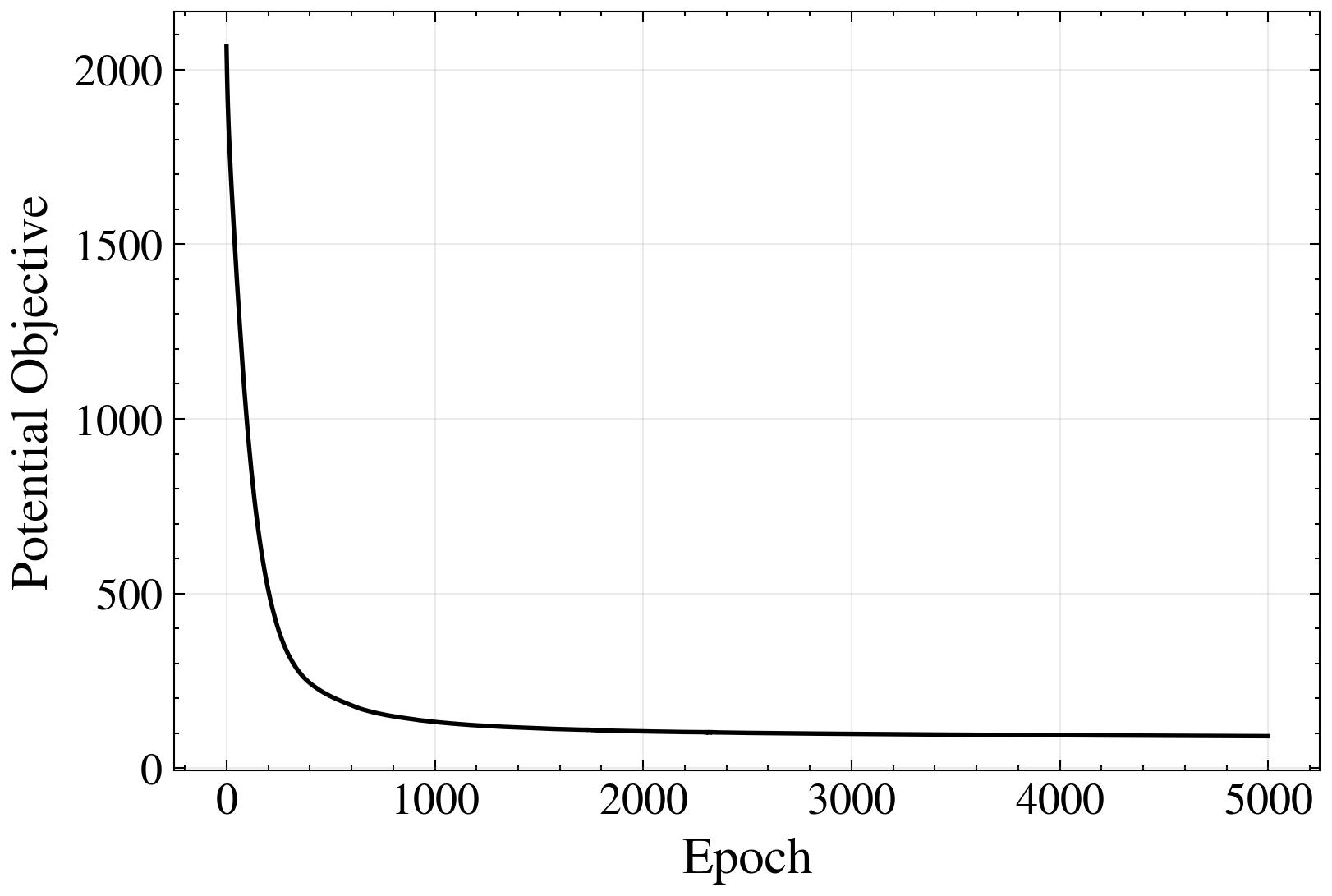}
        \caption{Potential objective during training}
        \label{fig:second}
    \end{subfigure}
    \caption{Vehicle trajectories and potential objective during training for the 30-vehicle control-velocity model.}
        \label{fig:30_vehicle}
\end{figure}

To further examine the scalability of the proposed approach, we consider a larger-scale experiment with $30$ vehicles under the control-velocity model with $\beta=1$. Figure~\ref{fig:30_vehicle} shows both the learned trajectories and the evolution of the potential objective during training. The trajectories remain well structured despite the substantially larger number of interacting vehicles, while the potential objective decreases steadily over the course of optimization. These results indicate that the proposed decentralized policy optimization framework remains effective and computationally tractable in a larger multi-vehicle setting.

\subsubsection{Two-dimensional acceleration control with obstacles}
\label{sec:obstacle}
This section  examines vehicle behavior in the presence of circular obstacles on the road.
We consider a 10-player game in two-dimensional AV control, i.e. $d=2, N=10$. The state dynamics follows \eqref{eq:AV_2} under acceleration control without noise ($\sigma_i \equiv 0$). Each vehicle $i$  minimizes its objective function   \eqref{eq:cost_i},
where
the terminal cost is 
$g_{i}(x)=2|x-z_{i}|^2$ for  $z_{i}=(1,1)$, and the interaction kernel $\lambda_{ij}$ and $ K$ are defined the same as in Sec.~\ref{sec:numerical_interaction} with $\beta=1$. 
The running cost $f_i$ depends on the position and the size of the obstacle. 
Specifically, we consider two circular obstacles centered at  $(0,0)$ with radii $0.1$ and $0.5$, respectively. To capture an obstacle avoidance behavior, we consider the running cost
$f_i(x,a_i)=0.02a_i^2+ h(x)$, with 
\begin{equation*}
    h(x)= 1000\left(1-\frac{1}{1+\exp[10(1-M|x|^2)]}\right),
\end{equation*}
where $M=100$ for the small obstacle of radius $0.1$
and $M=4$ for   the large obstacle of radius $0.5$.
The cost $h$ is a   differentiable  version of
forbidding passage through the   circular obstacles
\citep{baros2025mean}.
For comparison,  
we also consider the case without obstacles by setting
  $h\equiv 0$.

\begin{figure}[htbp]
    \centering
    \includegraphics[width=0.3\textwidth]{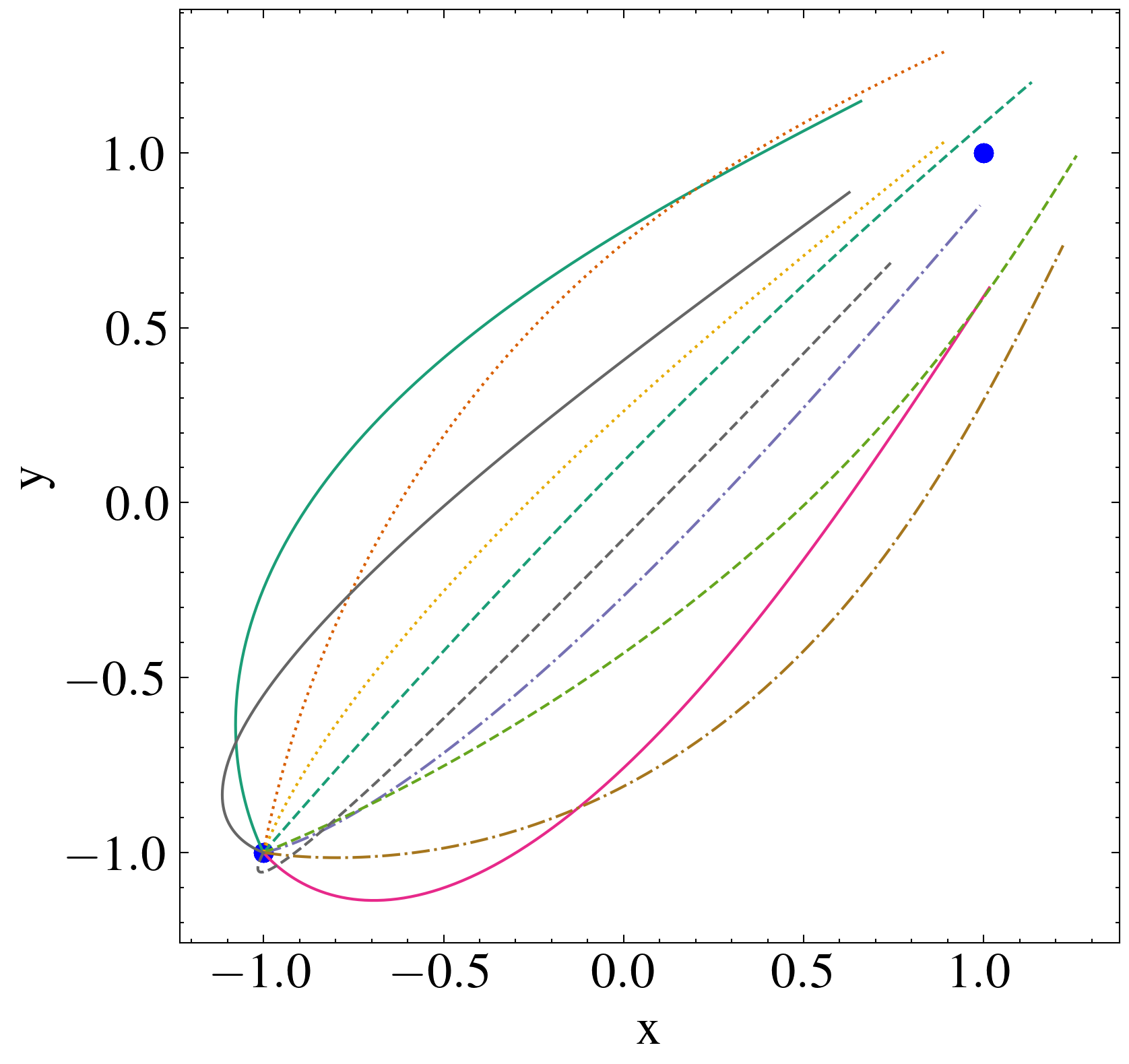}
    \includegraphics[width=0.3\textwidth]{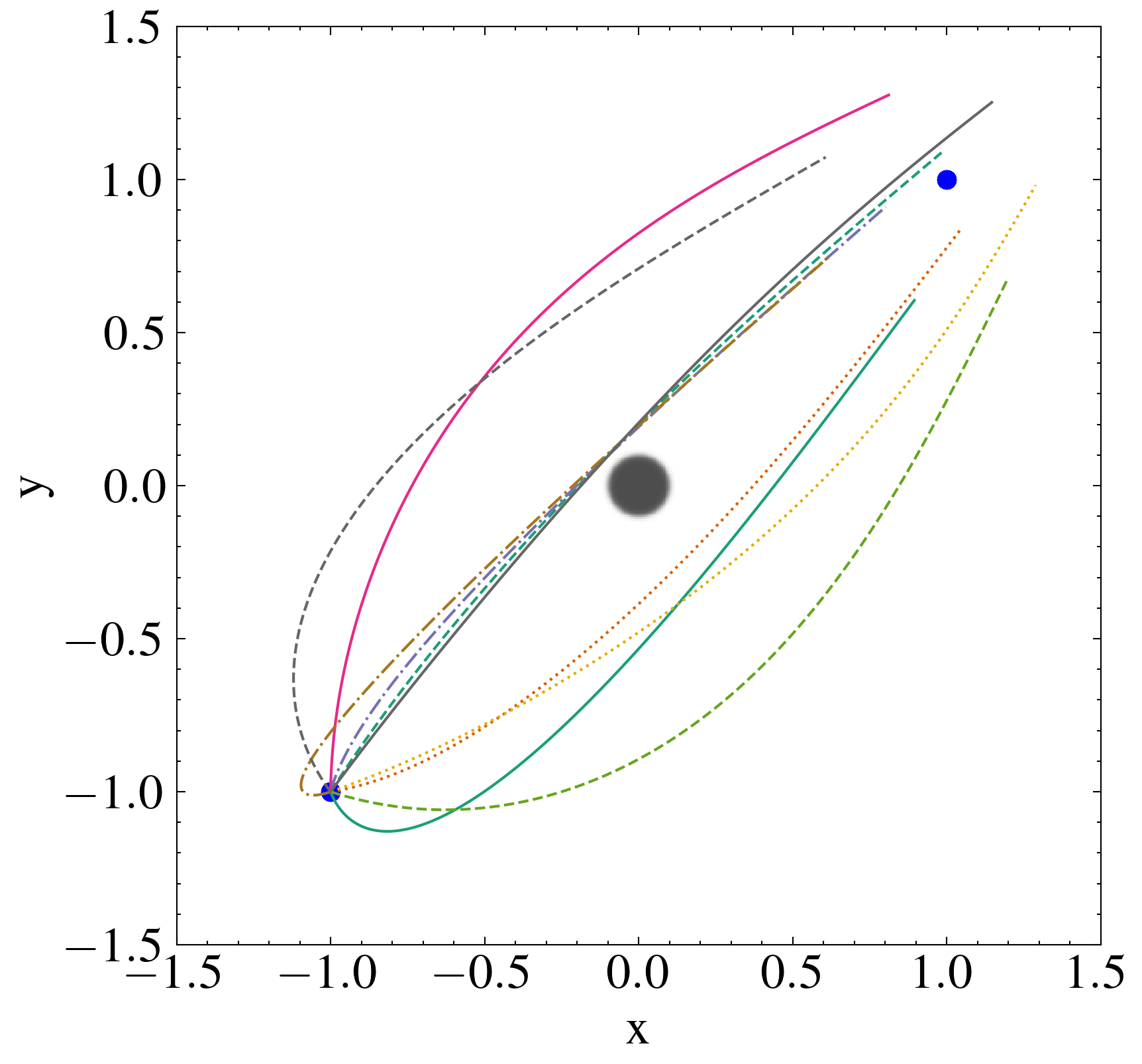}
    \includegraphics[width=0.3\textwidth]{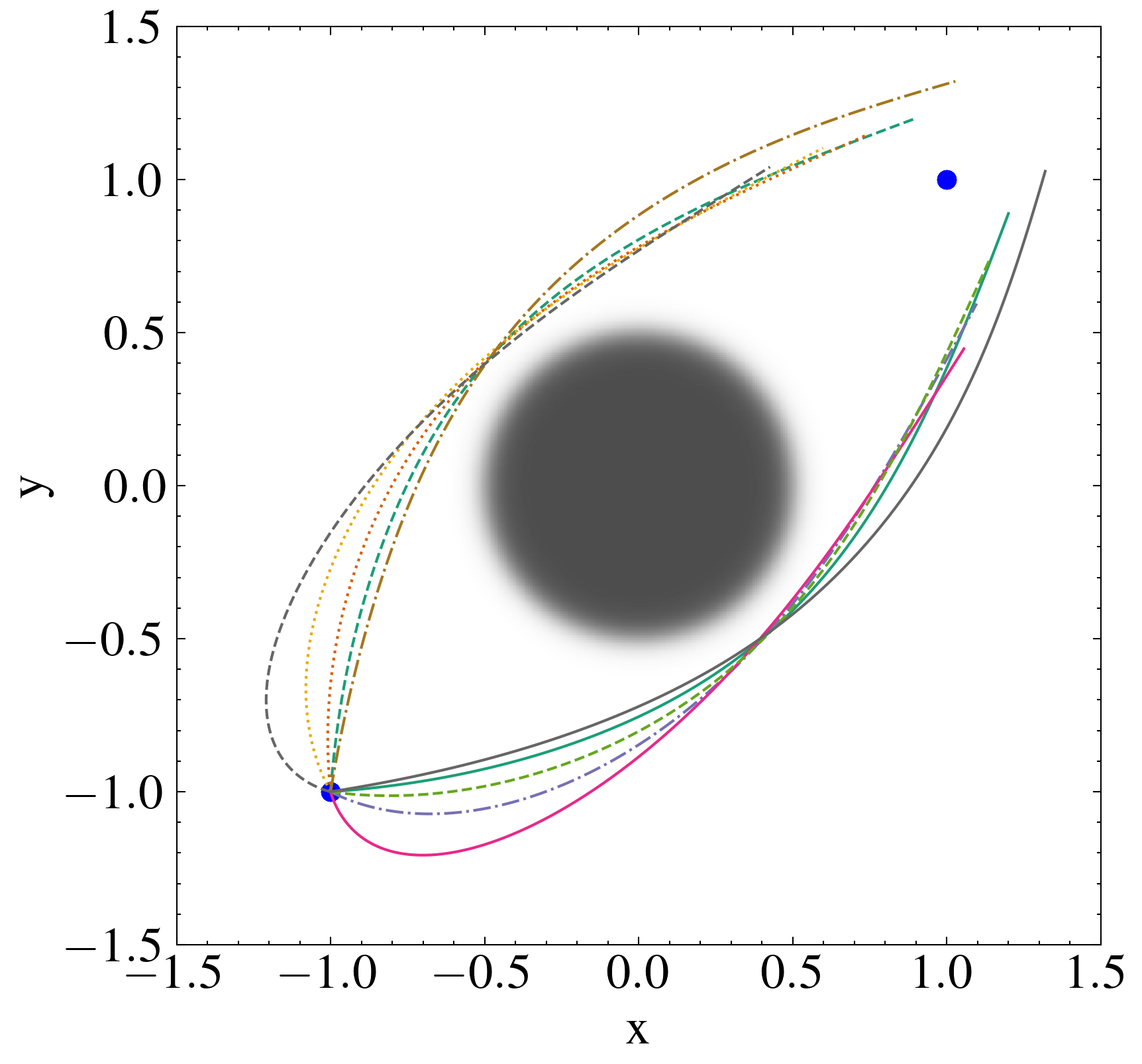}
    \caption{Vehicle trajectories under acceleration-control 
    without obstacle (left), with a small obstacle (middle), and with a large obstacle (right)}
    \label{fig:small_obstacle}
\end{figure}

Figure~\ref{fig:small_obstacle} shows that    introducing  an obstacle causes the vehicle trajectories to adjust accordingly  to satisfy obstacle-avoidance constraints. When a large obstacle is present, its wider exclusion region forces the trajectories to deviate earlier and more noticeably from their nominal paths. This results in a global reshaping of the motion patterns. In contrast, a smaller obstacle influences the trajectories only locally: only vehicles passing near its vicinity exhibit visible adjustments, while others remain close to the baseline case. These observations highlight that the trajectories respond adaptively to the obstacle, and the extent of deviation is determined by both the size and placement of the obstacle within the environment.

\subsubsection{
Impact of vehicle-specific rescaling for equilibrium computation
}
\label{sec:different_vehicle}

This section examines the effect of asymmetric  interactions and confirms that rescaling can improve the equilibrium approximation. 

\paragraph{
Rescaling recovers exact NEs.
}

We first consider a deterministic one-dimensional game with nine vehicles. We conduct this experiment under both control formulations: the vehicle dynamics follow \eqref{eq:AV_1} under velocity control and \eqref{eq:AV_2} under acceleration control, with no noise in either case.
We adopt the rank-one interaction structure introduced in Remark~\ref{rmk:rank-one-interactions}, with  $\lambda_{ij}=\gamma_i\tau_j$ for some $\gamma_i,\tau_i>0$. As shown in Remark~\ref{rmk:rank-one-interactions}, this interaction structure satisfies the Kolmogorov cyclic condition, and the rescaling $p_i=\tau_i/\gamma_i$ transforms the game into an exact potential game.
Consequently, an exact minimizer of the rescaled potential corresponds to a Nash equilibrium.

To model heterogeneous vehicle sizes, we divide the nine vehicles into three groups: large, medium, and small, with three vehicles in each group. We set $\tau_i=10$, $3$, and $0.5$ for the large, medium, and small vehicles, respectively. For each vehicle, we sample $\gamma_i\tau_i\sim\operatorname{Unif}(0.8,1.2)$ and determine $\gamma_i$ accordingly. Thus, larger vehicles have larger $\tau_i$ and smaller $\gamma_i$, so they exert stronger influence on other vehicles while being less sensitive to the influence of others. The resulting values of $\gamma_i$ and $\tau_i$ are reported in Table~\ref{table: vehicle types parameters}.  Each vehicle starts from the common initial position $x_0=-1$ and minimizes the objective function in \eqref{eq:cost_i}, with $f_i(x,a_i)=0.02a_i^2$, $g_i(x)=2|x-1|^2$, and interaction kernel $K(z)=(N^2|z|^2+1)^{-1}$. The running cost penalizes control effort, while the terminal cost encourages the vehicle's terminal position to remain close to the common target $x=1$.

\begin{table}[h]
\centering
\caption{Values of $\{\gamma_i\}_{i\in[9]}$ and $\{\tau_i\}_{i\in[9]}$.}\label{table: vehicle types parameters}
\scalebox{0.75}{\begin{tabular}{l c c c c c c c c c}
\toprule
Vehicle & 1 & 2 & 3 & 4 & 5 & 6 & 7 & 8 & 9 \\
\midrule 
$\gamma_i$ & 0.115 & 0.117 & 0.095 & 0.395 & 0.319 & 0.347 & 1.805 & 2.235 & 2.353 \\
$\tau_i$ & 10 & 10 & 10 & 3 & 3 & 3 & 0.5 & 0.5 & 0.5 \\
\bottomrule
\end{tabular}}
\end{table}

\begin{figure}[htbp]
    \centering
    \includegraphics[width=0.48 \textwidth]{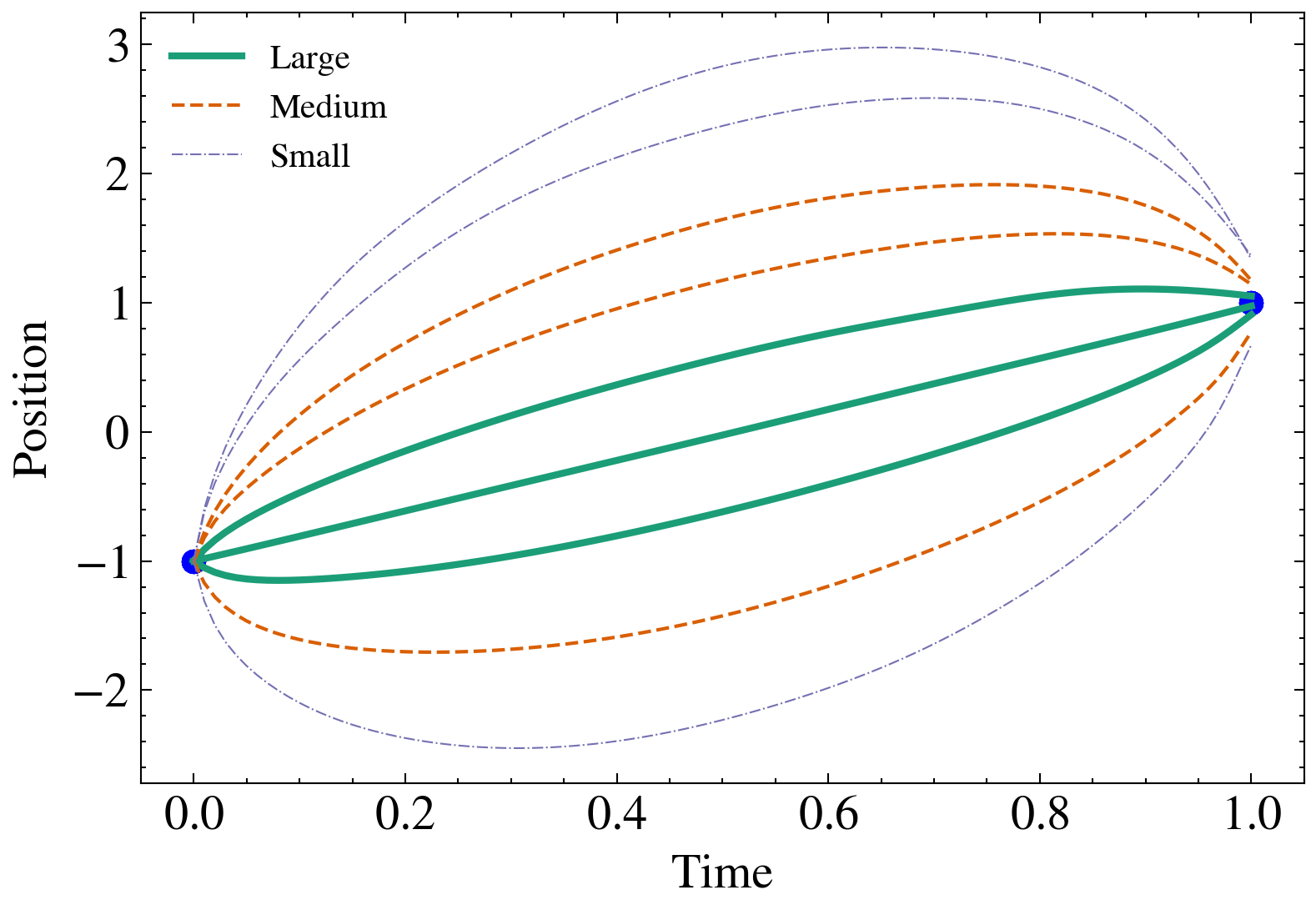}
    \quad 
    \includegraphics[width=0.48 \textwidth]{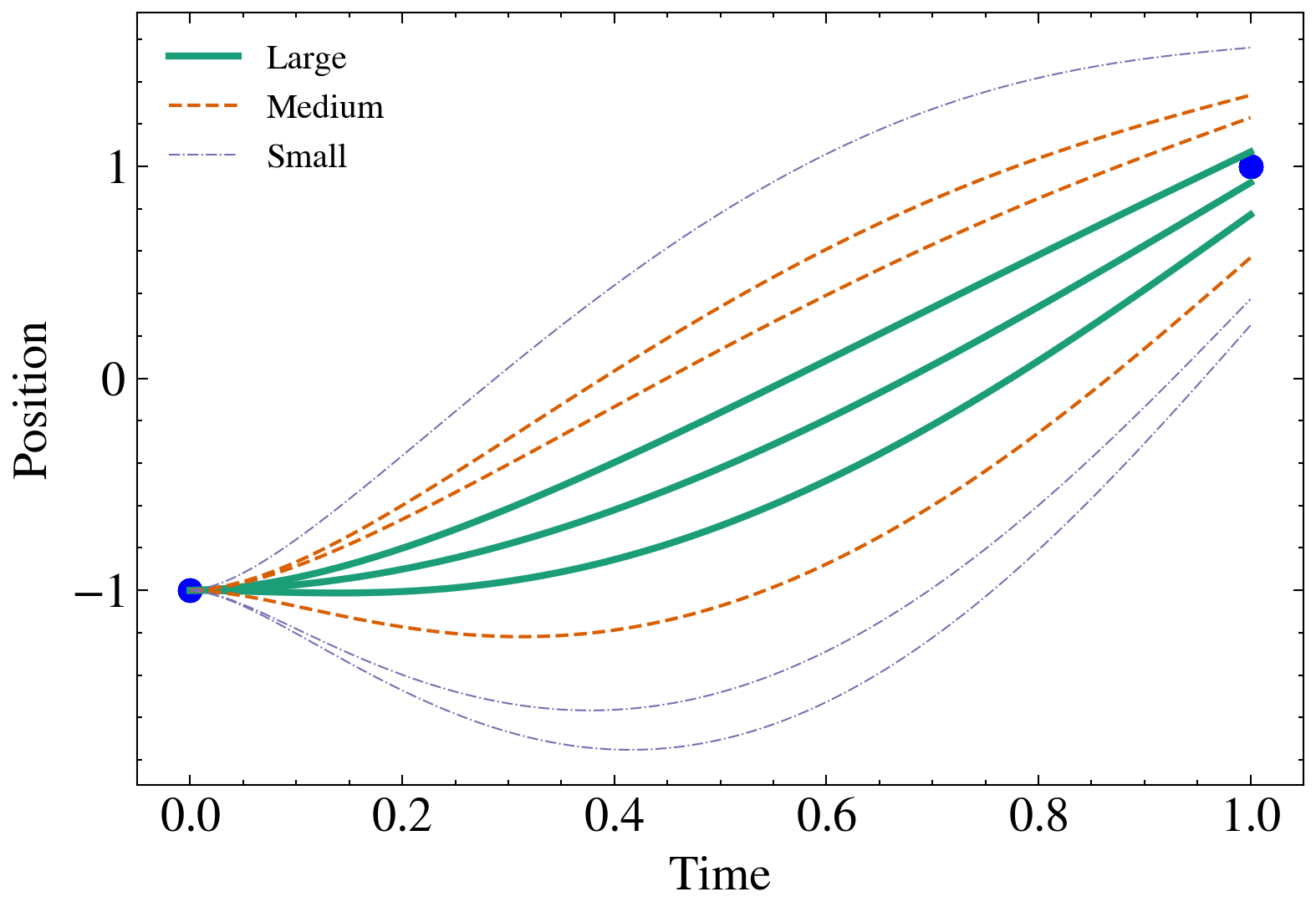}
    \caption{Vehicle trajectories of different types under control-velocity (left) and control-acceleration (right)}
    \label{fig:different_vehicle}
\end{figure}

Figure~\ref{fig:different_vehicle} shows a consistent size-dependent ordering of the vehicle trajectories under both velocity and acceleration control. Large vehicles remain near the center of the trajectory bundle, medium vehicles are displaced farther outward, and small vehicles move toward the periphery.
This pattern reflects the asymmetric interaction structure: larger vehicles exert a stronger influence on others while being less sensitive to their presence, whereas smaller vehicles are more strongly affected by the larger ones. The same qualitative ordering appears under both control formulations.

\paragraph{
Rescaling reduces the Nash gap.
}

We further consider settings in which rescaling does not recover an exact potential game and examine how it improves the equilibrium approximation. In this experiment, the vehicle dynamics follow \eqref{eq:AV_1} under velocity control with no noise. To this end, we introduce a small multiplicative Gaussian perturbation to the rank-one interaction matrix $(\lambda_{ij})_{i,j}$ by replacing each entry with $\widetilde{\lambda}_{ij}=\lambda_{ij}(1+\varepsilon\xi_{ij})$, where $\xi_{ij}\overset{\mathrm{i.i.d.}}{\sim}\mathcal{N}(0,1)$ and $\varepsilon>0$ controls the perturbation magnitude. We choose $\varepsilon$ sufficiently small so that all perturbed interaction coefficients remain positive in the realized instance. The perturbation breaks the rank-one structure, and the Kolmogorov cyclic condition no longer holds exactly.    We then apply Algorithm~\ref{alg:bisection-rescaling} to optimize the rescaling parameters $p=(p_1,\ldots,p_N)$ and minimize the corresponding upper bound on the equilibrium approximation error. 

To evaluate the resulting policy profile, for each vehicle $i$, we hold the strategies of all other vehicles fixed and compute its best-response cost
$J_i^{\mathrm{BR}}=\inf_{\phi_i'}J_i(\phi_i',\phi_{-i})$.
We define the exploitability as
$\mathcal E_i=J_i(\phi)-J_i^{\mathrm{BR}}$, which measures the cost reduction that vehicle $i$ could achieve through a unilateral deviation. At a Nash equilibrium, $\mathcal E_i=0$ for every vehicle, while smaller exploitability indicates a policy profile closer to a Nash equilibrium in terms of unilateral deviations. 
We further report the relative exploitability
$\mathcal E_i^{\mathrm{rel}}=\mathcal E_i/J_i(\phi)\times100\%$, which normalizes the exploitability by the vehicle's cost under the learned policy profile.

\begin{table}[h]
\centering
\caption{Exploitability before and after scaling.}
\label{tab:exploitability}
\scalebox{0.75}{\begin{tabular}{llccccccccc}
\toprule
& & \multicolumn{9}{c}{Vehicle} \\
\cmidrule(lr){3-11}
Scaling & Metric & 1 & 2 & 3 & 4 & 5 & 6 & 7 & 8 & 9 \\
\midrule
\multirow{2}{*}{Before}
& Cost $J_i(\phi)$
& 5.336 & 3.774 & 1.458 & 4.974 & 3.830 & 4.736 & 11.168 & 16.183 & 11.832 \\
& Relative exploitability (\%)
& 64.679 & 49.945 & 12.646 & 18.740 & 22.354 & 23.210 & 3.640 & 0.008 & 0.009 \\
\midrule
\multirow{2}{*}{After}
& Cost $J_i(\phi)$
& 1.637 & 1.712 & 1.642 & 5.258 & 4.497 & 7.148 & 8.930 & 13.070 & 13.290 \\
& Relative exploitability (\%)
& 0.180 & 0.795 & 0.928 & 1.598 & 1.145 & 0.927 & 3.717 & 0.962 & 2.275 \\
\bottomrule
\end{tabular}}
\end{table}

Table \ref{tab:exploitability} compares the  costs and relative exploitabilities before and after rescaling. Before rescaling, several vehicles have substantial incentives to deviate from the learned policy, with relative exploitabilities reaching \(64.7\%\) and \(49.9\%\) for Vehicles 1 and 2, respectively. After rescaling, the relative exploitabilities are substantially reduced for most vehicles and remain below \(4\%\) across all vehicles, with most  below \(2\%\). These results show that rescaling substantially reduces exploitability and produces a policy profile closer to a Nash equilibrium.

\subsection{Lane changing}
\label{sec:lane_changing}
In this section, we present  numerical experiments on lane-changing behavior using the previously established model framework. We consider a three-vehicle lane-changing problem on a two-lane road. The position of vehicle $i$ at time $t$ is denoted by $(x_{i,t},y_{i,t})$, where $x_{i,t}$ and $y_{i,t}$ are its longitudinal and lateral positions, respectively. The two lanes have lateral center positions $c_1$ and $c_2$, and the road boundaries are denoted by $y_{\min}$ and $y_{\max}$. Each vehicle controls its longitudinal and lateral velocities, denoted by $\phi_{i,x,t}$ and $\phi_{i,y,t}$, and seeks to travel close to its prescribed target speed while remaining within the roadway and avoiding collisions with the other vehicles.

The objective of vehicle $i$ is
\begin{equation}
\label{eq:lane_changing}
          J_i(\phi) = \mathbb{E} \Biggl[\int_0^T 
  \Biggl(f_i\big(x_{i,t},y_{i,t},\phi_{i,x,t},\phi_{i,y,t}) 
  +\sum_{j\not =i}K(x_{i,t}-x_{j,t})T(y_{i,t}-y_{j,t})
  \Biggr)\mathrm{d} t\Biggl].
\end{equation}
The individual running cost is decomposed as
$f_i=f_{i,1}+f_{i,2}$.
The first component is
\[
f_{i,1}(\cdot,\cdot,\phi_x,\phi_y)=q_1\phi_x^2+q_2(v_i^{\mathrm{des}})\phi_y^2+
q_3(\phi_x-v_i^{\mathrm{des}})^2,
\]
where $v_i^{\mathrm{des}}$ denotes the target speed of vehicle $i$. The first and third terms penalize control effort and deviations from the target longitudinal speed, respectively.  The second term penalizes lateral motion, with $q_2(z)=k_1/(k_2+k_3z)$. Since $q_2$ decreases with the target speed, a vehicle with a higher $v_i^{\mathrm{des}}$ incurs a smaller penalty for lateral motion and is therefore more willing to change lanes when doing so helps it maintain its preferred longitudinal speed.

The second component regulates the lateral position of the vehicle and is
defined by
\[
\begin{aligned}
f_{i,2}(x,y,\cdot,\cdot)
={}&
-k_4\log\left(
\sum_{j=1}^{n}
\exp\left(-k_5(y-c_j)^2\right)
\right) \\
&+
k_6\left[
\left(
\log\left(1+\exp\left(-k_7(y-y_{\min})\right)\right)
\right)^2
+
\left(
\log\left(1+\exp\left(-k_7(y-y_{\max})\right)\right)
\right)^2
\right].
\end{aligned}
\]
Here $n=2$ is the number of lanes, and $c_j$ denotes the lateral center of lane $j$. The first term smoothly penalizes deviations from the nearest lane center, thereby encouraging vehicles to travel near a lane center during normal driving. The second term penalizes departures from the roadway through soft penalties near the lower and upper road boundaries.

The pairwise interaction term in \eqref{eq:lane_changing} penalizes vehicle configurations that may lead to collisions. 
Specifically, $K(z) = \left( N^2 |z|^2 + 1 \right)^{-1}$ captures the effect of longitudinal separation, with the interaction becoming stronger as two vehicles move closer along the direction of travel, while
$T(z)=
k_8\left[
\tanh\left(k_{9}(d_{\mathrm{lat}}^2-z^2)\right)+1
\right]$ 
acts as a lateral weighting function that largely suppresses the interaction when the vehicles are separated across different lanes. 
Here $d_{\mathrm{lat}}$ is chosen slightly below the lane width, so that vehicles in the same lane interact strongly, whereas the interaction rapidly decays as their lateral separation approaches the spacing between adjacent lane centers.
Consequently,
$K(x_{i,t}-x_{j,t})T(y_{i,t}-y_{j,t})$ 
induces a strong penalty when two vehicles are longitudinally close and also occupy nearby lateral positions, thereby promoting collision avoidance. At the same time, vehicles traveling in different lanes have negligible interaction even when their longitudinal positions are close.

Table~\ref{tab:hyperparameters_lane_changing} summarizes the model parameters, initial vehicle positions, and target velocities used in the  experiments.

\begin{table}[H]
    \centering
    \caption{Parameter values and initial conditions for the lane-changing experiments in Section \ref{sec:lane_changing}.}
    \label{tab:hyperparameters_lane_changing}
\scalebox{0.75}{\begin{tabular}{c|cccccccccccccc}
    \hline
    Parameter 
    & $k_1$ & $k_2$ & $k_3$ & $k_4$ & $k_5$ 
    & $k_6$ & $k_7$ & $k_8$ & $k_9$ & $q_1$ & $q_3$ & $c_1$ & $c_2$ & $d_{\mathrm{lat}}$  \\
    \hline
    Value
    & 10.0 & 0.1 & 1.0 & 8.0 & 100.0 & 20.0 & 1.0 & 5.0 & 10.0 & 0.1 & 50 & 0.5 & 1.5 & 0.9 \\
    \hline
\end{tabular}}

\vspace{0.5em}

\scalebox{0.75}{\begin{tabular}{c|ccc}
    \hline
    Vehicle & $1$ & $2$ & $3$ \\
    \hline
    Initial position in Figure \ref{fig:comparison} (a)
    & $(-4.5, 0.5)$
    & $(-6.5, 0.5)$
    & $(0.5, 1.5)$ \\
    Initial position in Figure \ref{fig:comparison} (b)
    & $(-3.0, 0.5)$
    & $(-10.5, 0.5)$
    & $(-3.5, 1.5)$ \\
    Target velocities in Figure \ref{fig:comparison} (a) and (b)
    & $2.0$
    & $6.0$
    & $3.0$ \\
    \hline
\end{tabular}}
\end{table}

\begin{figure}[h]
    \centering
    \begin{subfigure}{0.48\textwidth}
        \centering
        \includegraphics[width=\linewidth]{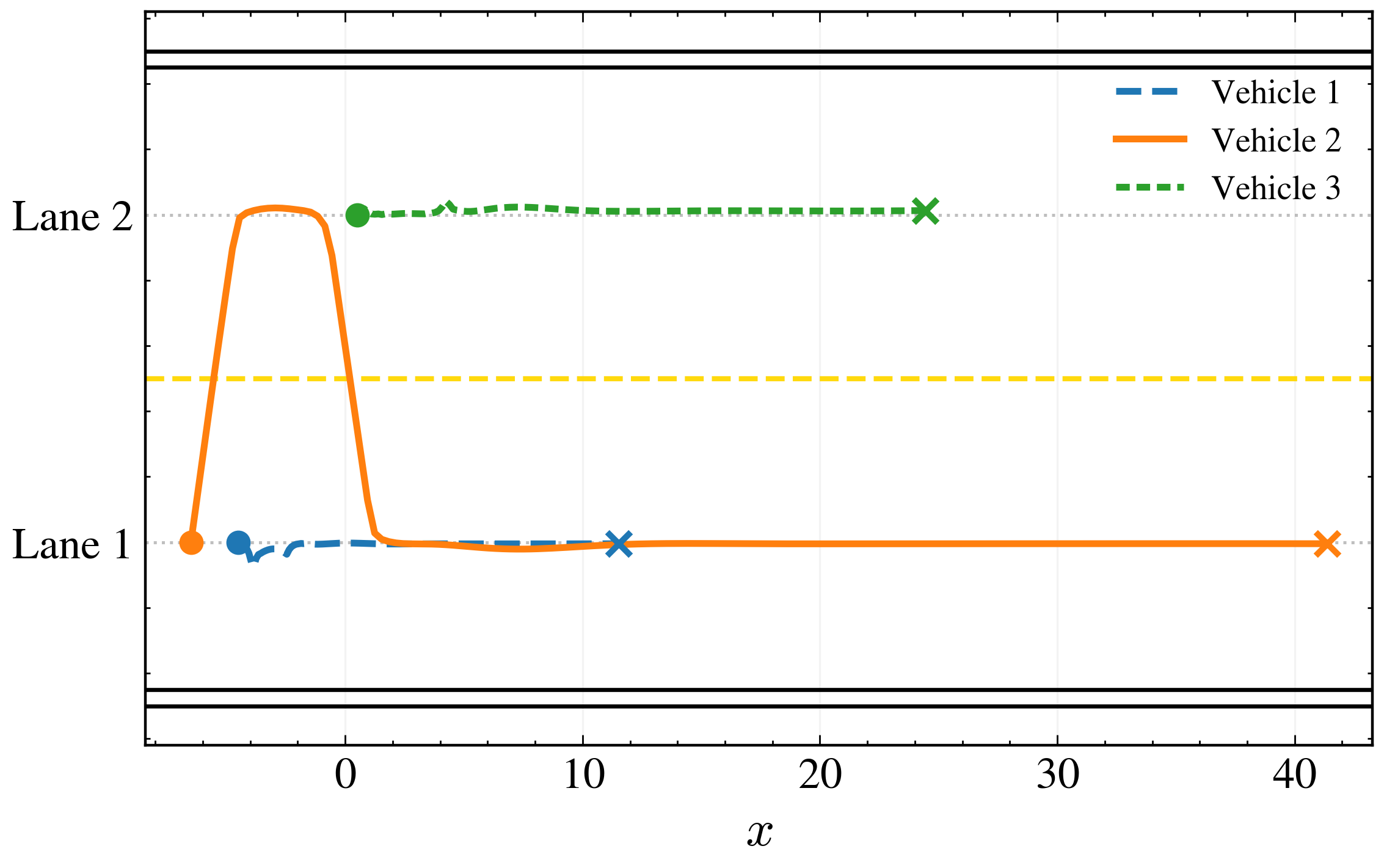}
        \caption{Standard overtaking behavior}
        \label{fig:first_lane_changing}
    \end{subfigure}
    \hfill
    \begin{subfigure}{0.48\textwidth}
        \centering
        \includegraphics[width=\linewidth]{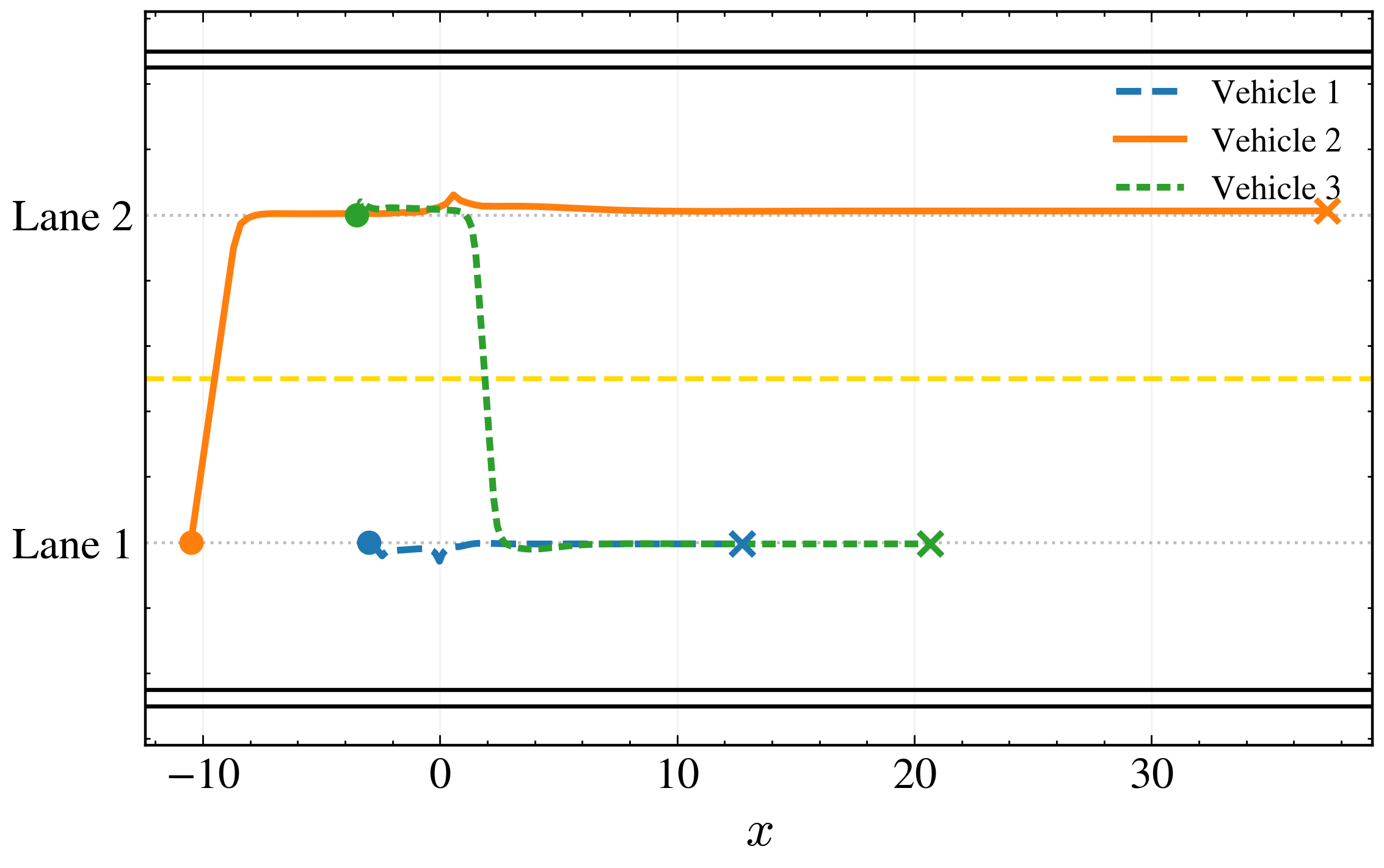}
        \caption{Alternative equilibrium lane-changing behavior}
        \label{fig:second_lane_changing}
    \end{subfigure}
    \caption{Lane changing in velocity control}
        \label{fig:comparison}
\end{figure}

We  use this model to examine how different initial traffic configurations can lead to distinct equilibrium lane-changing behaviors. 
The two scenarios use the same cost functions, cost parameters, and vehicle-specific target speeds, and differ only in the initial vehicle positions. In both scenarios, the orange vehicle has the highest target speed, followed by the green and blue vehicles. Initially, the orange vehicle travels in the right lane behind the slower blue vehicle, while the green vehicle travels in the left lane. The experiments therefore examine how the faster orange vehicle overtakes the blue vehicle while strategically interacting with the green vehicle in the adjacent lane. The behaviors are shown in
Figures~\ref{fig:comparison}(\subref{fig:first_lane_changing})
and~\ref{fig:comparison}(\subref{fig:second_lane_changing}), respectively.

Under the first initial configuration, the relative vehicle positions allow the orange vehicle to perform a standard overtaking maneuver. It first changes from the right lane to the left lane, passes the slower blue vehicle, and then returns to the right lane. The initial spacing leaves sufficient room for the orange vehicle to complete the return maneuver without incurring a large
interaction cost.

The second initial configuration leads to a qualitatively different equilibrium. The orange vehicle again changes to the left lane to overtake the blue vehicle. However, after passing the blue vehicle, returning to the right lane would place the orange vehicle too close to the blue vehicle, creating a high collision risk and therefore a large interaction cost.  The orange vehicle therefore remains in the left lane. Because the green vehicle has a lower target speed than the orange vehicle, the interaction cost of remaining in the left lane increases as the orange vehicle approaches from behind. In equilibrium, this interaction cost exceeds the cost of changing lanes, so the green vehicle moves to the right lane, allowing the two vehicles to maintain sufficient separation.

\subsection{Unsignalized intersection crossing}
\label{sec:intersection}
This section examines vehicle behavior at an unsignalized intersection. We consider an eight-vehicle setting in which each vehicle follows a prescribed travel direction and passes through the intersection.  Each vehicle is encouraged to maintain its prescribed target speed while avoiding collisions with other vehicles. The control variable is the acceleration,  and vehicle $i$ minimizes the objective function 
\begin{equation}
\label{eq:intersection_crossing}
J_i(\phi)
=
\mathbb{E}\left[
\int_0^T
\left(
f_i(x_{i,t},y_{i,t},v_{i,t},\phi_{i,t})+\sum_{j\neq i}C(x_{i,t},y_{i,t},x_{j,t},y_{j,t})\right)dt
\right].
\end{equation}
Here $(x_{i,t},y_{i,t})$ denotes the position of vehicle $i$ at time $t$, while $v_{i,t}$ and $\phi_{i,t}$ denote its velocity and acceleration along its assigned travel direction, respectively.   The running cost is given by 
\[f_{i}(\cdot,\cdot,v,\phi)=q_{4} \phi^2 + q_5^{i} (v - v_i^{\mathrm{des}})^2 + q_6
\left(\log\left(1+\exp\left(-q_7 v\right)\right)\right)^2.\]
The second term penalizes deviations from the target speed $v_i^{\mathrm{des}}$, with a vehicle-dependent coefficient $q_{5}^i$. The third term smoothly penalizes negative velocities and thereby discourages reverse motion.
The interaction takes the form 
\[
C(x_i,y_i,x_j,y_j)
=q_8 e^{-r_{ij}^2/2}+q_9
\left(\log\left(1+\exp\left(q_{10}(r_s-r_{ij})\right)\right)\right)^{2},
\]
where $ r_{ij}=\sqrt{(x_i-x_j)^2+(y_i-y_j)^2}.$ The first term penalizes close vehicle configurations, with the interaction cost increasing smoothly as the distance between two vehicles decreases. The second term is a smooth safety-distance penalty that becomes active when the inter-vehicle distance falls below the prescribed separation threshold $r_s$. This term is particularly important in realistic traffic settings, where vehicles occupy physical space rather than being treated as point masses.

Table~\ref{tab:hyperparameters_intersection} reports the model parameters, initial conditions, target speeds, and road-dependent speed-tracking coefficients used in the experiments.

\begin{table}[H]
    \centering
    \caption{Parameter values and initial conditions for the unsignalized intersection-crossing experiments. 
    }
    \label{tab:hyperparameters_intersection}
\scalebox{0.75}{\begin{tabular}{c|cccccc}
    \hline
    Parameter 
    & $q_4$ & $q_6$ & $q_7$ & $q_8$ & $q_9$ & $q_{10}$ \\
    \hline
    Value
    & 0.1 & 1.0 & 10.0 & 10.0 & 5.0 & 10.0  \\
    \hline
\end{tabular}}

\vspace{0.5em}

\scalebox{0.75}{\begin{tabular}{c c c c c c }
\hline
Road type & Position & Initial position & Initial speed & Target speed & Speed-tracking
coefficient $q_{5}^i$\\
\hline
Major & Lead ($\uparrow$) & $(1.0,  -1.5)$ & 1.5 & 2.0 & 300.0 \\
Major & Lead ($\downarrow$)  & $(-1.0, 1.5)$ & 1.5 & 2.0 & 300.0 \\
Major & Follow ($\uparrow$) & $(1.0, -3.0)$  & 1.5 & 2.0 & 300.0 \\
Major & Follow ($\downarrow$)  & $(-1.0, 3.0)$  & 1.5 & 2.0 & 300.0 \\
Minor & Lead ($\rightarrow$) & $(-1.5,  -0.5)$ & 1.5 & 2.0 & 10.0 \\
Minor & Lead ($\leftarrow$)  & $(1.5, 0.5)$ & 1.5 & 2.0 & 10.0 \\
Minor & Follow ($\rightarrow$) & $(-3.0,  -0.5)$ & 1.5 & 2.0 & 10.0 \\
Minor & Follow ($\leftarrow$)  & $(3.0,  0.5)$ & 1.5 & 2.0 & 10.0 \\
\hline
\end{tabular}}

\end{table}

\begin{figure}[h]
\centering

\textbf{With major--minor-road differentiation}\par\vspace{1mm}

\begin{subfigure}[t]{0.24\textwidth}
\centering
\includegraphics[width=\linewidth]{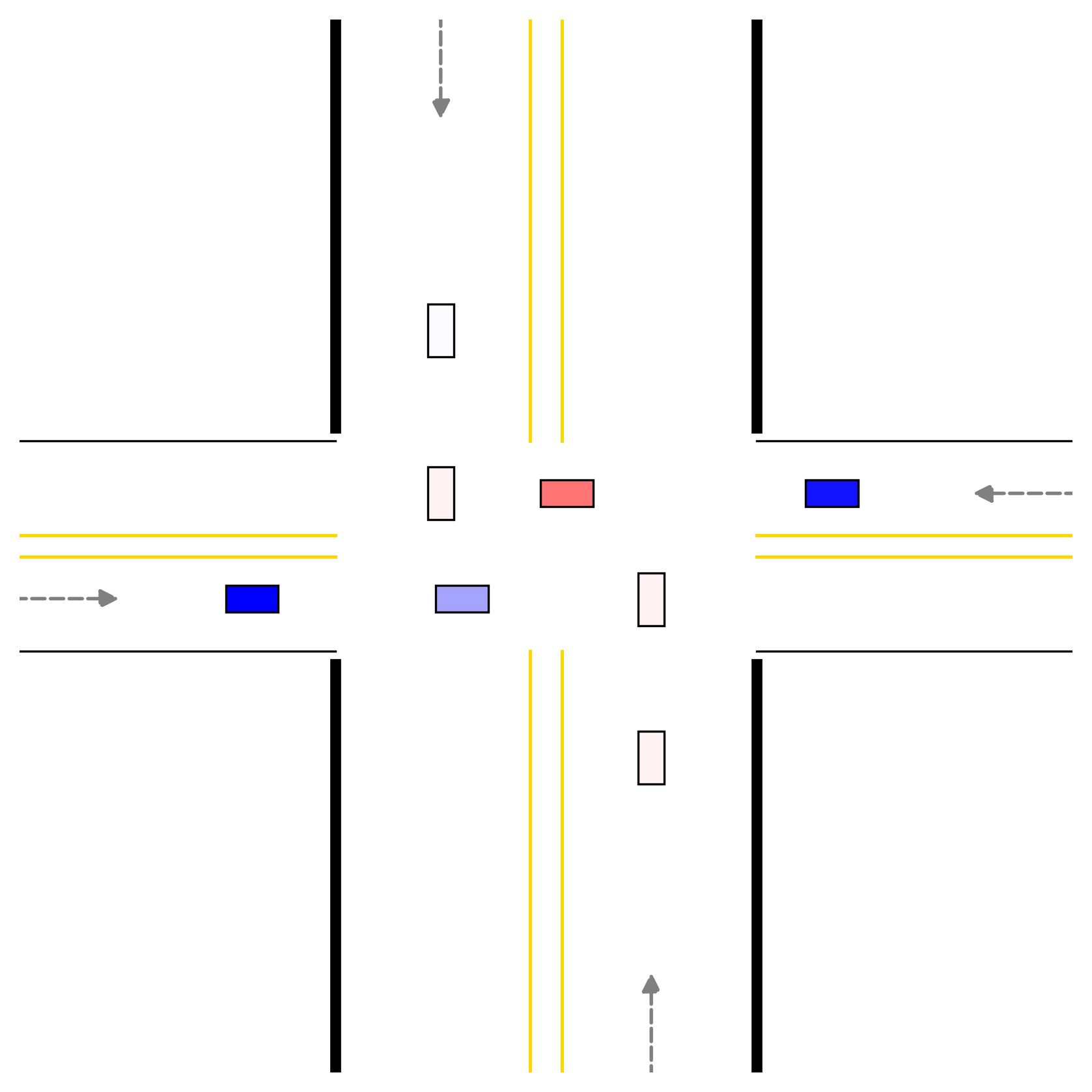}
\caption{$t=0.5$}
\end{subfigure}\hfill
\begin{subfigure}[t]{0.24\textwidth}
\centering
\includegraphics[width=\linewidth]{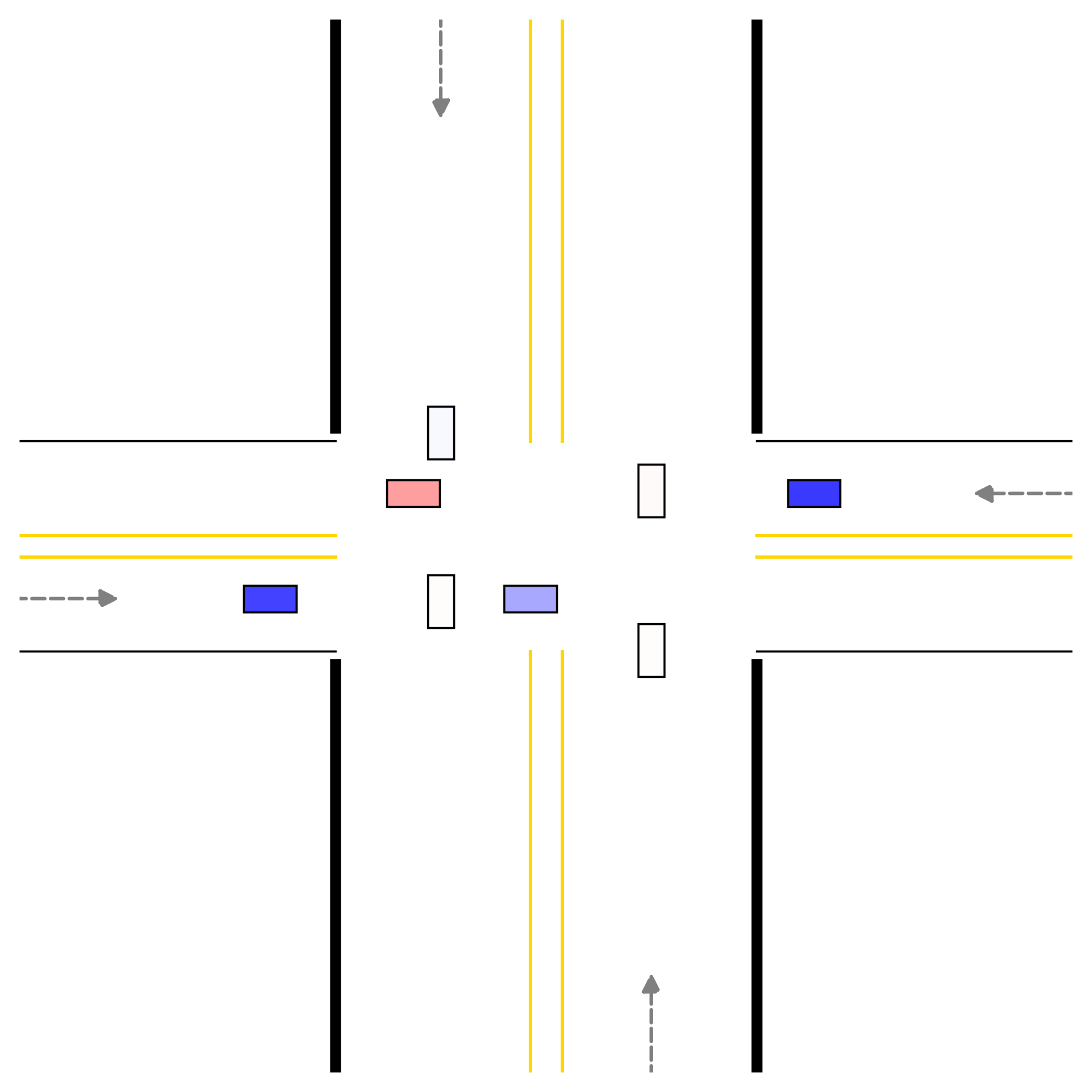}
\caption{$t= 1.0$}
\end{subfigure}\hfill
\begin{subfigure}[t]{0.24\textwidth}
\centering
\includegraphics[width=\linewidth]{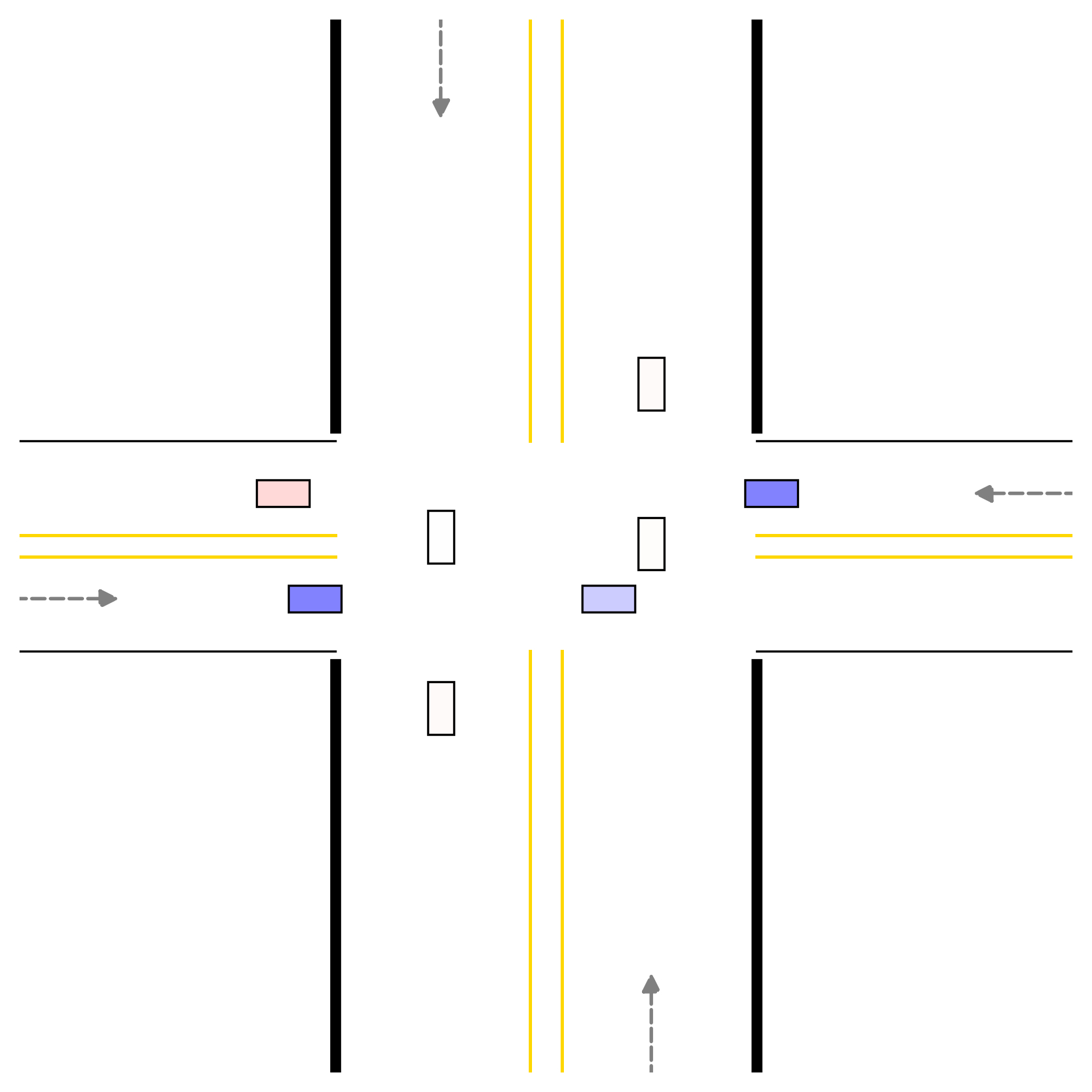}
\caption{$t= 1.5$}
\end{subfigure}\hfill
\begin{subfigure}[t]{0.24\textwidth}
\centering
\includegraphics[width=\linewidth]{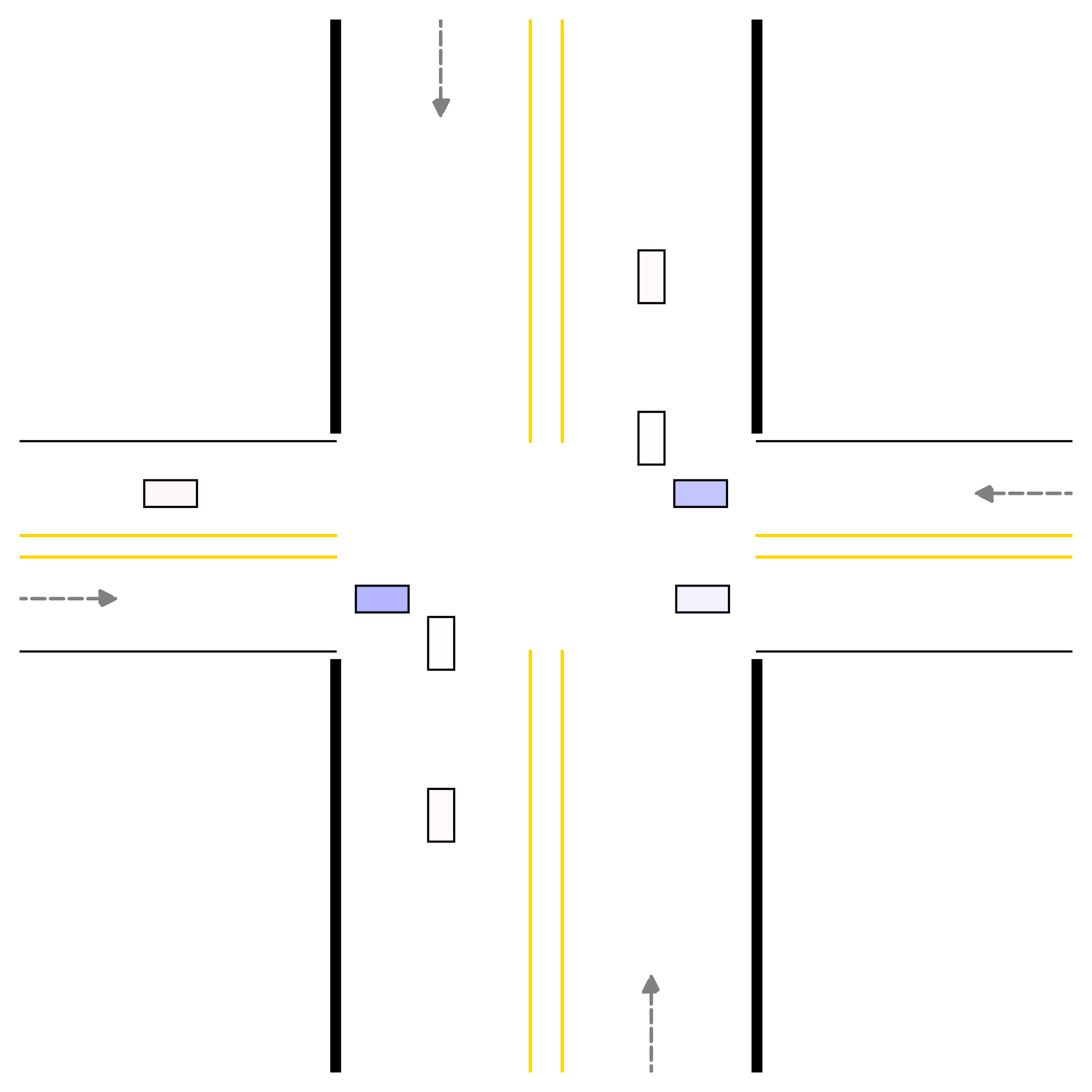}
\caption{$t=2.0$}
\end{subfigure}

\par\vspace{1mm}
\includegraphics[width=0.8\textwidth]
{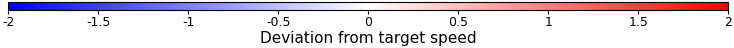}

\vspace{2mm}

\textbf{Without major--minor-road differentiation}\par\vspace{1mm}

\begin{subfigure}[t]{0.24\textwidth}
\centering
\includegraphics[width=\linewidth]{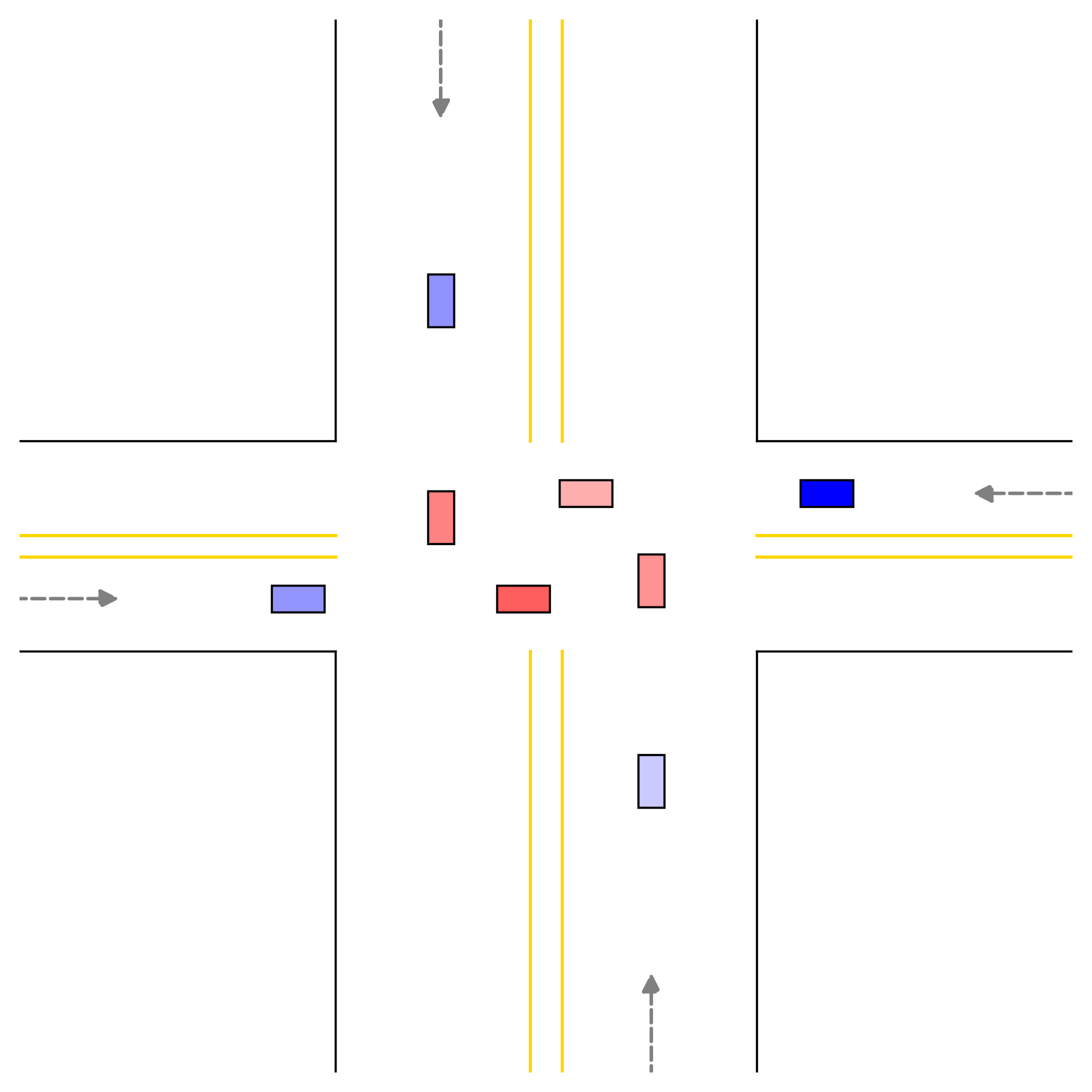}
\caption{$t=0.5$}
\end{subfigure}\hfill
\begin{subfigure}[t]{0.24\textwidth}
\centering
\includegraphics[width=\linewidth]{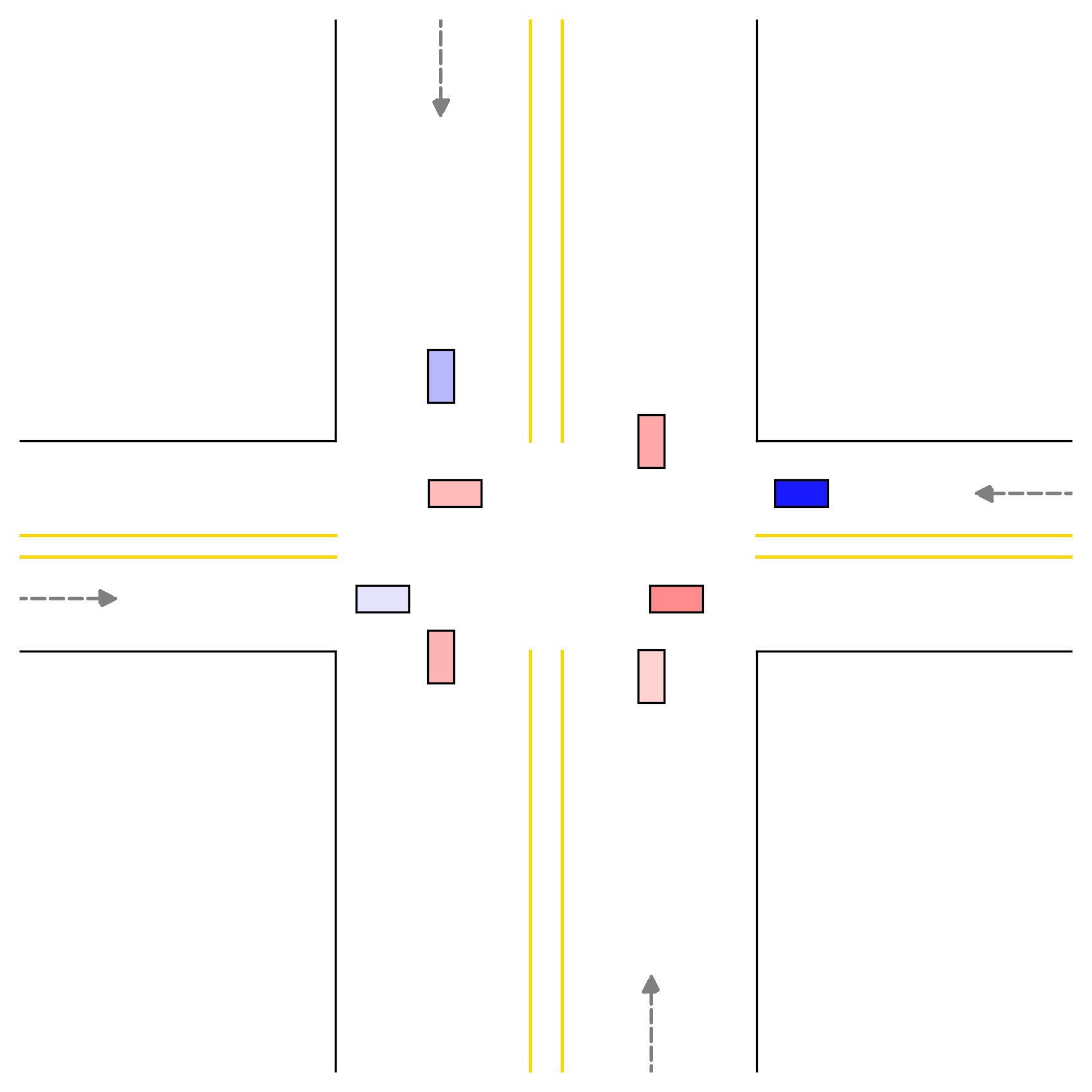}
\caption{$t=1.0$}
\end{subfigure}\hfill
\begin{subfigure}[t]{0.24\textwidth}
\centering
\includegraphics[width=\linewidth]{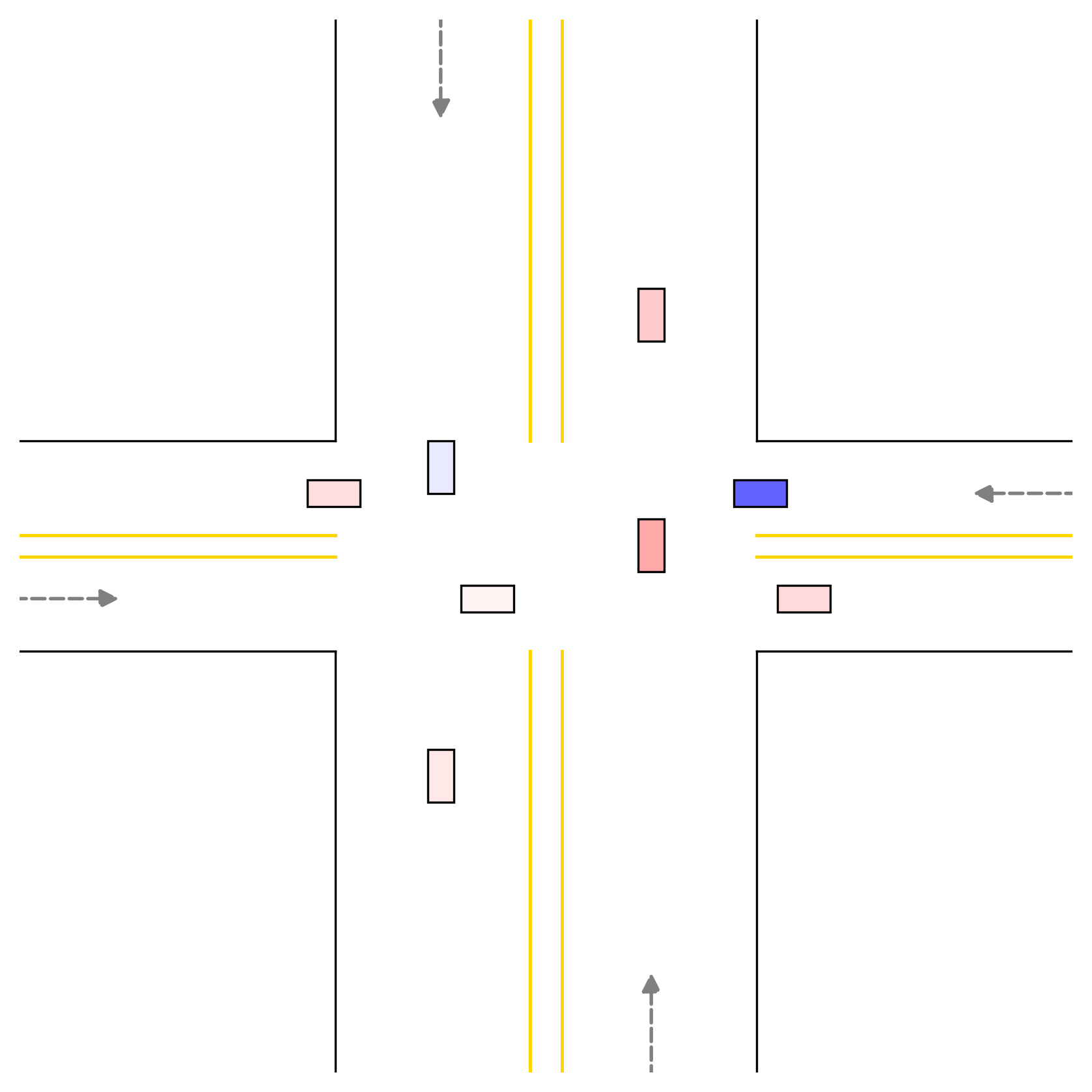}
\caption{$t=1.5$}
\end{subfigure}\hfill
\begin{subfigure}[t]{0.24\textwidth}
\centering
\includegraphics[width=\linewidth]{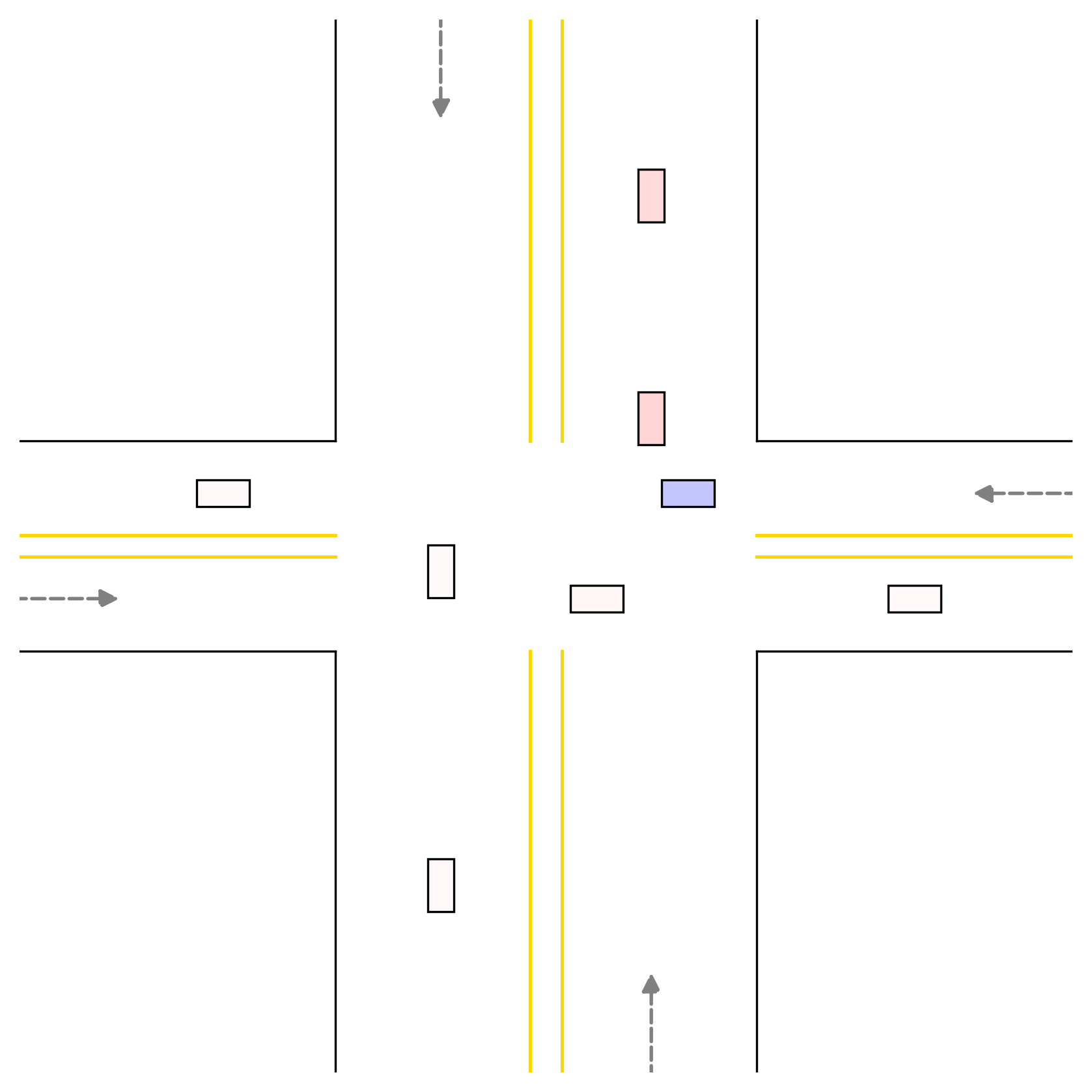}
\caption{$t=2.0$}
\end{subfigure}

\par\vspace{1mm}
\includegraphics[width=0.8\textwidth]
{keyframe_major_minor/speed_deviation_colorbar.png}

\caption{Comparison of the intersection-crossing dynamics at four
representative time instants. The top row shows the
major--minor-road scenario, while the bottom row shows the scenario without
major--minor-road differentiation. Vehicle colors indicate deviations from
their target speeds.}
\label{fig:crossing_comparison}
\end{figure}

We compare the equilibrium intersection-crossing behaviors with and without major–minor-road differentiation under identical initial conditions. In the major--minor-road setting, the vertical road is designated as the major road and the horizontal road as the minor road, with vehicles on the minor road expected to yield to those on the major road. This priority structure is encoded through the vehicle-dependent speed-tracking coefficient $q_5^i$: vehicles on the major road are assigned a larger $q_5^i$, making deviations from their target speeds more costly, whereas vehicles on the minor road are assigned a smaller $q_5^i$, allowing them to decelerate more readily when resolving potential conflicts. In the setting without major--minor-road differentiation, the same speed-tracking penalty is used across roads.
In both scenarios, the vehicles start from the same prescribed positions, with initial speeds slightly below their respective target speed. 
Figure~\ref{fig:crossing_comparison} shows four representative snapshots from each scenario. Vehicle colors indicate deviations from the corresponding target speeds: blue represents speeds below the target, red represents speeds above the target. Colors closer to white indicate smaller deviations, while darker shades correspond to larger deviations in
magnitude.

In the major--minor-road scenario, a clear priority-based crossing pattern emerges. Vehicles traveling on the major road exhibit similar behaviors and maintain speeds close to their target values while passing through the intersection. Among the minor-road vehicles, those that have already entered or are approaching the conflict region accelerate to clear the intersection  and appear in red, indicating speeds above their target values. In contrast, the following vehicles decelerate and wait for the major-road vehicles to pass, appearing in dark blue as their speeds fall substantially below their targets. Once the major-road vehicles have cleared the intersection, the waiting vehicles accelerate and proceed through the intersection.

Without major--minor-road differentiation, no predetermined road hierarchy governs the crossing order. Instead, vehicles resolve potential conflicts through their acceleration and deceleration decisions, with leading vehicles proceeding first and following vehicles yielding as needed. The resulting crossing pattern is therefore determined by the vehicles' relative positions and strategic interactions rather than by road priority.

\section{Proofs}
We provide in this section the proofs of the main results established throughout the paper.

\subsection*{Proof of Proposition~\ref{prop:alpha_PG}}
\noindent
For all  $i\in[N]$, define a new objective  of player $i$:
\begin{equation}
\label{eq:cost_i_symmetric}
\begin{aligned}
J^s_i(\phi&) = \sE \Biggl[\int_0^T 
  \Biggl(
  f_i\big(
  X^{\phi_i}_{i,t},
  \phi_{i}(t,X^{\phi_i}_{i,t})\big) 
   +\sum_{j\not =i}\frac{\lambda_{ij}+\lambda_{ji}}{2}K(X^{\phi_i}_{i,t}-X^{\phi_j}_{j,t})
  \Biggr)\d t  + g_i(X^{\phi_i}_{i,T})\Biggl].
\end{aligned}
\end{equation}
By \cite[Theorem 3.1]{guo2025towards}, the game with the objectives $(J^s_i)_{i\in [N]}$ is a potential game with the potential function $\Phi$ defined in \eqref{eq:potential_fun_symmetric}, which implies
\[
\Phi(\phi_i,\phi_{-i})-\Phi(\phi_i',\phi_{-i}) = J_i^s(\phi_i,\phi_{-i})-J^s_i(\phi_i',\phi_{-i}).
\]
Note that
\begin{align*}
    |J_i^s(\phi)-J_i(\phi)| 
    \leq& \sE\left|\int_0^T\sum_{j\neq i}\left(\frac{1}{2}(\lambda_{ij}+\lambda_{ji})-\lambda_{ij}\right)K(X^{\phi_i}_{i,t}-X^{\phi_j}_{j,t})\d t \right|\le   \frac{T\|K\|_{L^\infty}}{2}\sum_{i\neq j}|\lambda_{ij}-\lambda_{ji}|.
\end{align*}
The desired conclusion then follows from
\begin{align*}
    &\left|\left(J_i(\phi_i,\phi_{-i})-J_i(\phi_i',\phi_{-i})\right)-\left(\Phi(\phi_i,\phi_{-i})-\Phi(\phi_i',\phi_{-i})\right)\right|\\
    =& \left|\left(J_i(\phi_i,\phi_{-i})-J_i(\phi_i',\phi_{-i})\right)-\left(J_i^s(\phi_i,\phi_{-i})-J^s_i(\phi_i',\phi_{-i})\right)\right|\\
    \leq& \left|J_i(\phi_i,\phi_{-i})-J_i^s(\phi_i,\phi_{-i})\right|+\left|J_i(\phi_i',\phi_{-i})-J^s_i(\phi_i',\phi_{-i})\right|  
    \leq T\|K\|_{L^\infty}\sum_{i\neq j}|\lambda_{ij}-\lambda_{ji}|. \hfill \tag*{\Halmos}
\end{align*}

\subsection*{Proof of Theorem~\ref{thm:rescaled-alpha-potential}}
\noindent
By Proposition \ref{prop:alpha_PG},    $\Phi^p$ satisfies 
for all $i\in [N]$,
$\phi_{-i}\in \Pi_{-i}$, and
$\phi_i,\phi'_i\in \Pi_i$, 
\begin{equation*}
|(p_iJ_i(\phi'_i,\phi_{-i})
-p_iJ_i(\phi_i,\phi_{-i}))
-(\Phi^p(\phi'_i,\phi_{-i})
-\Phi^p(\phi_i,\phi_{-i}))|\le \tilde  \alpha_i(p),
\end{equation*}
with 
$
\tilde \alpha_i(p)=T\|K\|_{L^\infty} \sum_{j\ne i} |p_i\lambda_{ij}-p_j\lambda_{ji}|. 
$
Hence a minimizer $\phi$ of $\Phi^p$ satisfies for all $i\in [N]$ and  $\phi'_i\in \Pi$,
$$
J_i(\phi'_i,\phi_{-i})
-J_i(\phi_i,\phi_{-i})\ge \frac{1}{p_i}(\Phi^p(\phi'_i,\phi_{-i})
-\Phi^p(\phi_i,\phi_{-i}))-\frac{1}{p_i}\tilde  \alpha_i(p)\ge -\frac{1}{p_i}\tilde  \alpha_i(p)\ge \alpha(p),
$$
which implies that $\phi$ is an   $\alpha(p)$-NE of   the game  \eqref{eq:state_i}-\eqref{eq:cost_i}.

We proceed to show that $\Phi^p$ admits a minimizer. 
 Since $\Pi_i$ is compact, $\Pi$ equipped with the product topology is also compact in the compact-open topology. It suffices to show that $\Phi^p$ is continuous.
    The extreme value theorem then guarantees that $\Phi^p$ admits a minimizer on $\Pi$.

To show the continuity of $\Phi^p$, consider a sequence $(\phi^n)_{n\in \sN}\subset \Pi$ with $ (\phi^n_i)_{i\in [N]}$ converging to some $\phi^\infty\in \Pi$ in the compact-open topology. 
In the sequel, we write $X^{n}_i= X^{\phi^n_i}_i$ and 
$\Delta X^{n}_i= X^{\phi^n_i}_i-X^{\phi^\infty_i}_i$ for all $n\in \sN\cup\{\infty\}$ and $i\in [N]$, and 
denote by $C\ge 0$ a generic constant 
independent of $n\in \sN\cup\{\infty\}$ and $i\in [N]$, which may take different values in each line.  

  Since $A=\prod_{i\in [N]}A_i$ is compact,  $(\phi^n)_{n\in \sN\cup\{\infty\}}$ is uniformly bounded. 
    The standard moment estimate in \cite[Theorem   3.4.3]{zhang2017backward} implies that  for all $n\in \sN$ and $i\in [N]$, 
  $\sE[\sup_{t\in [0,T]}|X^{n}_{i,t}|^4]\le C$.
  For each $n\in \sN$ and $R>0$, 
  define the stopping time 
$\tau_i^{n,R}\coloneqq \inf\{t\in[0,T]\mid |X_t^{n}|>R \textrm{ or } |X_t^{\infty}|>R\}\wedge T.
$ 
  By the stability estimate \cite[Theorem 3.2.4]{zhang2017backward} and the locally Lipschitz continuity of $\phi^\infty$, 
\begin{align*}
    \sE[\sup_{t\in [0,T]}|\Delta X^{n}_{i,t\wedge \tau_i^{n,R}}|^2]
    &\le C\bigg(\sE\left[\int_0^{\tau_i^{n,R}} |b_i(X^n_{i,t},\phi^n_i(t,X^n_{i,t}))-b_i(X^n_{i,t},\phi^\infty_i(t,X^n_{i,t}))|^2\d t \right]
    \\
    &\quad + \sE\left[\int_0^{\tau_i^{n,R}} |\sigma_i(X^n_{i,t},\phi^n_i(t,X^n_{i,t}))-\sigma_i(X^n_{i,t},\phi^\infty_i(t,X^n_{i,t}))|^2\d t \right]\bigg)
    \\
    &\le C\sE\left[\int_0^{\tau_i^{n,R}} |  \phi^n_i(t,X^n_{i,t}))-\phi^\infty_i(t,X^n_{i,t})|^2\d t \right], 
\end{align*}
which converges to zero as $n\to \infty$, due to the   the convergence of 
$(\phi^n)_{n\in \sN}$ on compact sets.
Then by the Cauchy-Schwarz inequality,
\begin{align*}
      \sE[\sup_{t\in [0,T]}|\Delta X^{n}_{i,t}|^2]
      &=\sE[\sup_{t\in [0,T]}|\Delta X^{n}_{i,t}|^21_{\tau_i^{n,R}=T}]
      +\sE[\sup_{t\in [0,T]}|\Delta X^{n}_{i,t}|^21_{\tau_i^{n,R}<T}]
      \\
      &\le \sE[\sup_{t\in [0,T]}|\Delta X^{n}_{i,\tau_i^{n,R}}|^2]
      +\sE[\sup_{t\in [0,T]}|\Delta X^{n}_{i,t}|^4]^{1/2}P({\tau_i^{n,R}<T})^{1/2}.
\end{align*}
Since $\{\tau_i^{n,R}<T\} \subset \left\{\sup_{t\in[0,T]} |X_t^{n}|\ge R\right\}\cup \left\{\sup_{t\in[0,T]} |X_t^{\infty}|\ge R\right\}$, 
by   Markov's inequality and  the uniform moment bounds of $(X^n_i)_{n\in \sN}$, we have  
$\lim_{n\to \infty}\sE[ \sup_{t\in [0,T]} |X^{n}_{i,t}-X^{\infty}_{i,t}|^2]=0$.
Similar arguments   show that 
$ 
\lim_{n\to \infty}\sup_{t\in [0,T]}\sE[|\phi^n_i(t, X^{n}_{i,t})-\phi^\infty_i(t, X^{\infty}_{i,t})|^2]=0    
$. 
Hence, possibly after passing to a subsequence,    
 $(X^{n}_{i,t},\phi^n_i(t, X^{n}_{i,t}))_{n\in \sN}$ also converges to $(X^{\infty}_{i,t},\phi^{\infty}_i(t, X^{\infty}_{i,t})) $  $\d \sP\otimes \d t $-a.e., and 
 $(X^{n}_{i,T})_{n\in \sN}$   converges to $X^{\infty}_{i,T} $  a.s. 

Let  $F:\sR^{dN}\times \prod_{i\in [N]}A_i\to \sR$ and $G: \sR^{dN}\to \sR$  be continuous functions with at most quadratic growth. Then
$ 
(F\bigl(X_t^{n},\phi^n(t,X_t^{n})\bigr))_{n\in \sN}$ converges to $
F\bigl(X_t^\infty,\phi^\infty(t,X_t^\infty)\bigr)$ $\d \sP\otimes \d t $-a.e.,
and
$ 
(G\bigl(X_T^{n}\bigr))_{n\in \sN}$ converges to $ G\bigl(X_T^\infty\bigr)$ a.s. 
The quadratic growth of $F$ and $G$, 
the $L^2$ convergence of 
$(X^{n}_{i,t},\phi^n_i(t, X^{n}_{i,t}))_{n\in \sN}$, 
and  Vitali's theorem imply
\[
\lim_{n\to \infty}
\sE\left[\int_0^T F\bigl(X_t^{n},\phi^n(t,X_t^{n})\bigr)dt + G\bigl(X_T^{n}\bigr)\right]
=
\sE\left[\int_0^T F\bigl(X_t^\infty,\phi^\infty(t,X_t^\infty)\bigr)dt + G\bigl(X_T^\infty\bigr)\right].
\]
Applying this   to     $F=F^p$ and $G=G^p$ in \eqref{eq:potential_fun_symmetric_p} yields the desired continuity of $\Phi^p$ and finishes the proof.
\hfill \Halmos

\subsection*{Proof of Theorem~\ref{thm:optimal-rescaling}}
\noindent
Since there are no interactions between distinct components, the objective in \eqref{eq:optimal-rescaling-value} is invariant under multiplying all $p_i$, $i\in S_k$, by the same positive constant, independently for each $k$. Hence, without loss of generality, we may impose
\[
    \sum_{i\in S_k}p_i=|S_k|,
    \qquad k=1,\ldots,m.
\]
Introducing the epigraph variable $r$ then gives \eqref{eq:optimal-rescaling-compact}, initially with
$p\in(0,\infty)^N$ and $r\ge0$.

We next show that the strict positivity constraint can be relaxed to $p\in[0,\infty)^N$. Suppose that a feasible solution satisfies $p_i=0$ for some $i\in S_k$. The $i$-th constraint gives $\sum_{j\ne i}p_j\lambda_{ji}=0$.
Hence $p_j=0$ for every $j$ such that $\lambda_{ji}>0$, i.e., for every out-neighbor $j$ of $i$ in $\mathcal G_\lambda$. Repeating this argument along directed paths starting from $i$ and using the strong connectivity of the subgraph induced by $S_k$, we obtain $p_j=0$ for any  $j\in S_k$. 
This contradicts
\[
    \sum_{j\in S_k}p_j=|S_k|>0.
\]
Thus every feasible solution of \eqref{eq:optimal-rescaling-compact} satisfies
$p_i>0$ for all $i$.

Finally, the componentwise normalizations imply
$0\le p_i\le |S_k|$, for all $i\in S_k$.
Also, choosing $p_i=1$ for all $i$ is feasible with $r=r_0$, and hence $r^\star\le r_0$.
Therefore one may restrict $r$ to $[0,r_0]$ without changing the optimal value. The feasible set of \eqref{eq:optimal-rescaling-compact} is then closed and bounded, and hence compact. The minimum is therefore attained, and, by the preceding argument, any optimizer belongs to $(0,\infty)^N$.
\hfill \Halmos

\subsection*{Proof of Theorem~\ref{thm:exact-rescaling}}
\noindent 
The equivalence between the Kolmogorov cyclic condition  and the existence of $p\in(0,\infty)^N$ satisfying \eqref{eq:detailed-balance-lambda} 
follows from the standard characterization of diagonally symmetrizable matrices; see, e.g., \cite[Proposition~4.15 and Lemma~4.12]{mckee2020symmetrizable}.

Suppose first that \eqref{eq:detailed-balance-lambda} holds. Then, for
every $i\in[N]$,
$ 
\sum_{j\ne i}|p_i\lambda_{ij}-p_j\lambda_{ji}|=0,
$ 
and hence $r^\star=0$.
Conversely, suppose that $r^\star=0$. Then there exists a sequence
$(p^n)_{n\ge1}\subset(0,\infty)^N$ such that
\[
r(p^n)
\coloneqq
\max_{i\in[N]}
\frac{1}{p_i^n}
\sum_{j\ne i}
|p_i^n\lambda_{ij}-p_j^n\lambda_{ji}|
\longrightarrow 0.
\]
For every $i\ne j$,
\[
\left|
\lambda_{ij}
-
\frac{p_j^n}{p_i^n}\lambda_{ji}
\right|
\le r(p^n).
\]
Thus, if $\lambda_{ji}=0$, then $\lambda_{ij}=0$; interchanging $i$
and $j$ shows that the zero pattern is symmetric.
Now consider any interacting pair $i\ne j$, namely    $\lambda_{ij}+\lambda_{ji}>0$. By the symmetric zero
pattern, $\lambda_{ij},\lambda_{ji}>0$, and therefore
$\lim_{n\to \infty}
\frac{p_j^n}{p_i^n}
=
\frac{\lambda_{ij}}{\lambda_{ji}}.
$ 
Hence, along any cycle
$i_1,\ldots,i_m,i_{m+1}=i_1$ with positive interaction weights,
\[
1
=
\prod_{k=1}^m
\frac{p_{i_{k+1}}^n}{p_{i_k}^n}
\longrightarrow
\prod_{k=1}^m
\frac{\lambda_{i_k i_{k+1}}}
     {\lambda_{i_{k+1}i_k}}, \quad n\to \infty,
\]
which yields \eqref{eq:kolmogorov-cycle}. If a cycle contains a zero
interaction coefficient, then the symmetric zero pattern implies that
the corresponding reverse coefficient is also zero, so both sides of
\eqref{eq:kolmogorov-cycle} vanish. Thus the Kolmogorov cyclic
condition holds.

Finally, let $p^\star\in(0,\infty)^N$ satisfy
\eqref{eq:detailed-balance-lambda}. Then
$\max_{i\in[N]}
\frac{1}{p_i^\star}
\sum_{j\ne i}
|p_i^\star\lambda_{ij}-p_j^\star\lambda_{ji}|=0.
$ 
Therefore Theorem~\ref{thm:rescaled-alpha-potential} gives
$\alpha(p^\star)=0$, and every minimizer of $\Phi^{p^\star}$ is an NE of \eqref{eq:state_i}--\eqref{eq:cost_i}.
\hfill \Halmos

\subsection*{Proof of Proposition~\ref{prop:bisection-convergence}}
\noindent
After $n$ bisection iterations, the interval
$[r_{\rm low},r_{\rm high}]$ contains $r^\star$ and has length
$r_0/2^n$. Since the algorithm returns $r^n=r_{\rm high}$,
$0\le r^n-r^\star\le \frac{r_0}{2^n}$.
Moreover, $p^n$ is feasible for $r^n$, and hence
\[
    r^\star
    \le
    \max_{i\in[N]}
    \frac{1}{p_i^n}
    \sum_{j\ne i}
    |p_i^n\lambda_{ij}-p_j^n\lambda_{ji}|
    \le r^n.
\]
Multiplying by $T\|K\|_{L^\infty}$ gives the bound for
$\alpha(p^n)-\alpha^\star$.
 \hfill \Halmos

\subsection*{Proof of Proposition~\ref{prop:efficiency-rescaled-potential}}
\noindent
Under Assumption~\ref{assum:compact_policy}, the existence of $\phi^p\in\arg\min_{\phi\in\Pi}\Phi^p(\phi)$ follows from Theorem~\ref{thm:rescaled-alpha-potential}.

By the nonnegativity of the cost components $C_i(\phi)$ and $\mathcal K_{ij}(\phi)$, the definitions of $\underline{\gamma}(p)$ and $\overline{\gamma}(p)$ yield
\begin{equation} 
\label{eq:potential-social-cost-sandwich} 
\underline{\gamma}(p)\operatorname{SC}(\phi) \le \Phi^p(\phi) \le \overline{\gamma}(p)\operatorname{SC}(\phi), \qquad \forall \phi\in\Pi. 
\end{equation}
Under Assumption~\ref{assum:compact_policy}, an argument similar to the proof of Theorem~\ref{thm:rescaled-alpha-potential} shows that the social cost admits a minimizer.  Let $\phi^{\rm opt} \in \arg\min_{\phi\in\Pi}\operatorname{SC}(\phi)$. 
Using \eqref{eq:potential-social-cost-sandwich} and the optimality of $\phi^{p}$ for $\Phi^{p}$ gives
\[ 
\operatorname{SC}(\phi^{p}) \le \frac{1}{\underline{\gamma}(p)} \Phi^{p}(\phi^{p}) \le \frac{1}{\underline{\gamma}(p)} \Phi^{p}(\phi^{\rm opt}) \le \frac{\overline{\gamma}(p)} {\underline{\gamma}(p)} \operatorname{SC}(\phi^{\rm opt}).
\]
Dividing by $\operatorname{SC}^{\rm opt}>0$ and using Definition~\ref{def:social-inefficiency}, proves \eqref{eq:efficiency-guarantee}.
 \hfill \Halmos

\subsection*{Proof of Corollary~\ref{cor:pos-bound}}
\noindent
    By Theorem~\ref{thm:exact-rescaling}, the rescaled potential with rescale parameter $p$ is an exact potential of the rescaled game, and any of its global minimizers is a Nash equilibrium of the original game. 
    Hence $\operatorname{PoS}\le\operatorname{InEff}(\phi^p)$. 
   For every pair satisfying $\lambda_{ij}+\lambda_{ji}>0$, \eqref{eq:detailed-balance-lambda} gives
    \[ 
    \frac{p_i\lambda_{ij}+p_j\lambda_{ji}} {2(\lambda_{ij}+\lambda_{ji})} = \frac{p_ip_j}{p_i+p_j} \le \min\{p_i,p_j\}. 
    \]
    It follows that
    $ 
    \overline{\gamma}(p)=\max_{i\in[N]}p_i, $ and $\underline{\gamma}(p) \ge \min_{1\le i<j\le N} \frac{p_ip_j}{p_i+p_j}. 
    $ 
    Therefore, Proposition~\ref{prop:efficiency-rescaled-potential} yields the first inequality in \eqref{eq:pos-kolmogorov-simple}. Furthermore, 
    \[ 
    \frac{p_ip_j}{p_i+p_j} \ge \frac12\min\{p_i,p_j\} \ge \frac12\min_{k\in[N]}p_k. 
    \]
which gives the second inequality in \eqref{eq:pos-kolmogorov-simple}.

If the interaction coefficients are symmetric, one can take $p_i=1$ for all $i\in[N]$ and obtain  $\operatorname{PoS}\le 2$. 
\hfill \Halmos

\section{Conclusions and Future Work} 
\label{sec:conclude}

We develop an \(\alpha\)-potential game framework for decentralized decision-making of heterogeneous autonomous vehicles, reducing approximate NE computation to potential minimization. We introduce player-specific rescaling to sharpen the equilibrium approximation and   under suitable conditions, recover exact NE  despite asymmetric interactions. We further derive social-efficiency guarantees relative to centralized full cooperation. We implement the framework using a neural-network policy-gradient method and demonstrate its performance across collision and obstacle avoidance, lane changing, and intersection crossing. Numerical results highlight the ability of the framework to capture heterogeneous behavior and strong local interactions, and show that rescaling can substantially reduce exploitability in asymmetric games.

Future work will address model uncertainty and data-driven learning, particularly how errors in estimated dynamics and interaction parameters affect equilibrium accuracy and social efficiency. Further priorities include convergence analysis of the policy-gradient algorithms and systematic comparisons with HJB-based approaches. Additional directions include equilibrium selection when multiple equilibria exist, extensions to policies using local neighborhood information, and the incorporation of explicit safety constraints. Validation in larger, more realistic traffic environments will further assess scalability and practical applicability.

\section*{Acknowledgments}

The corresponding author acknowledges support from the National Science Foundation under the award number CMMI-1943998. Anran Hu acknowledges financial support by InnoHK initiative, The Government of the HKSAR and the AIFT Lab.

\bibliographystyle{elsarticle-harv}
\bibliography{TSref}

@article{candogan2013near,
  title={Near-potential games: Geometry and dynamics},
  author={Candogan, Ozan and Ozdaglar, Asuman and Parrilo, Pablo A},
  journal={ACM Transactions on Economics and Computation (TEAC)},
  volume={1},
  number={2},
  pages={1--32},
  year={2013},
  publisher={ACM New York, NY, USA}
}

@article{monderer1996potential,
  title={Potential games},
  author={Monderer, Dov and Shapley, Lloyd S},
  journal={Games and economic behavior},
  volume={14},
  number={1},
  pages={124--143},
  year={1996},
  publisher={Elsevier}
}

@article{mckee2020symmetrizable,
  title={Symmetrizable integer matrices having all their eigenvalues in the interval $[-2, 2] $},
  author={McKee, James and Smyth, Chris},
  journal={Algebraic Combinatorics},
  volume={3},
  number={3},
  pages={775--789},
  year={2020}
}

@book{albrecht2024multi,
  title={Multi-agent reinforcement learning: Foundations and modern approaches},
  author={Albrecht, Stefano V and Christianos, Filippos and Sch{\"a}fer, Lukas},
  year={2024},
  publisher={MIT Press}
}

@article{zhang2024stackelberg,
  author  = {Zhang, Q. and Langari, R. and Tseng, H. E.
             and Mohan, S. and Szwabowski, S. and Filev, D.},
  title   = {{Stackelberg} Differential Lane Change Game
             Based on {MPC} and Inverse {MPC}},
  journal = {IEEE Transactions on Intelligent Transportation Systems},
  year    = {2024},
  volume  = {25},
  pages   = {8473--8485},
  url     = {https://ieeexplore.ieee.org/document/10507755}
}

@article{zhang2025integrated,
  author  = {Zhang, Lei and Zheng, Jiacheng and Zhang, Zhiqiang
             and Wang, Zhenpo and Wang, Mingqiang},
  title   = {Integrated Motion Planning for On-Ramp Merging
             Based on {Stackelberg} Game Modeling Considering
             Interactive Characteristics},
  journal = {IEEE Transactions on Vehicular Technology},
  year    = {2025},
  volume  = {74},
  number  = {8},
  pages   = {11762--11776},
  doi     = {10.1109/TVT.2025.3554978}
}

@article{yao2025personalized,
  author  = {Yao, Tianluo and Jin, Hui},
  title   = {A Personalized Lane-Changing Decision System
             Based on Improved {Stackelberg} Game
             and Traffic Flow Information},
  journal = {IEEE Transactions on Intelligent Transportation Systems},
  year    = {2025},
  volume  = {26},
  number  = {5},
  pages   = {6789--6801},
  doi     = {10.1109/TITS.2025.3531921}
}

@article{fu2023cooperative,
  author  = {Fu, M. and Li, S. and Guo, M. and Yang, Z.
             and Sun, Y. and Qiu, C. and Wang, X. and Li, X.},
  title   = {Cooperative Decision-Making of Multiple Autonomous
             Vehicles in a Connected Mixed Traffic Environment:
             A Coalition Game-Based Model},
  journal = {Transportation Research Part C: Emerging Technologies},
  year    = {2023},
  volume  = {157},
  pages   = {104415},
  doi     = {10.1016/j.trc.2023.104415}
}

@article{fu2025regional,
  author  = {Fu, M. and Li, S. and Guo, M.
             and Wang, X. and Li, X. and Wang, W.},
  title   = {Regional Cooperative Decision-Making Based
             on Coalition Game for Multilane Merging
             in Mixed Traffic},
  journal = {IEEE Transactions on Intelligent Transportation Systems},
  volume  = {26},
  number  = {12},
  pages   = {22665 - 22679},
  year    = {2025}
}

@article{huang2024robust,
  author  = {Huang, J. and Wu, Z. and Xue, W.
             and Lin, D. and Chen, Y.},
  title   = {Non-Cooperative and Cooperative Driving Strategies
             at Unsignalized Intersections:
             A Robust Differential Game Approach},
  journal = {IEEE Transactions on Intelligent Transportation Systems},
  year    = {2024},
  volume  = {25},
  pages   = {9535--9549},
  url     = {https://ieeexplore.ieee.org/document/10440183}
}

@article{shu2025decision,
  author  = {Shu, Keqi and Ning, Minghao and Alghooneh, Ahmad
             and Li, Shen and Pirani, Mohammad
             and Khajepour, Amir},
  title   = {Decision Making in Urban Traffic:
             A Game Theoretic Approach for Autonomous Vehicles
             Adhering to Traffic Rules},
  journal = {IEEE Transactions on Intelligent Transportation Systems},
  year    = {2025},
  doi     = {10.1109/TITS.2025.3553077},
  url     = {https://arxiv.org/abs/2505.10690}
}

@article{jing2025decentralized,
  author  = {Jing, D. and Yao, E. and Chen, R.
             and Men{\'e}ndez, M.},
  title   = {Decentralized Human-Like Ramp Merging
             Decision-Making and Control Based
             on a Stochastic Potential Game},
  journal = {IEEE Transactions on Intelligent Transportation Systems},
  year    = {2025},
  volume  = {26},
  number  = {10},
  pages   = {16724--16734},
  url     = {https://ieeexplore.ieee.org/document/11033200}
}

@article{wang2023equilibrium,
  author  = {Wang, Hua and Wang, Jing and Chen, Shukai
             and Meng, Qiang},
  title   = {Equilibrium Traffic Dynamics with Mixed
             Autonomous and Human-Driven Vehicles
             and Novel Traffic Management Policies:
             The Effects of Value-of-Time Compensation
             and Random Road Capacity},
  journal = {Transportation Science},
  year    = {2023},
  volume  = {57},
  number  = {5},
  pages   = {1177--1208},
  doi     = {10.1287/trsc.2021.0469}
}

@article{shen2022distributed,
  author  = {Shen, Jinglai and Kammara, Eswar Kumar H.
             and Du, Lili},
  title   = {Fully Distributed Optimization-Based {CAV}
             Platooning Control Under Linear Vehicle Dynamics},
  journal = {Transportation Science},
  year    = {2022},
  volume  = {56},
  number  = {2},
  pages   = {381--403},
  doi     = {10.1287/trsc.2021.1100}
}

@article{shou2022markov,
  author  = {Shou, Zhenyu and Chen, Xu and Fu, Yongjie
             and Di, Xuan},
  title   = {Multi-agent reinforcement learning for
             {Markov} routing games: A new modeling
             paradigm for dynamic traffic assignment},
  journal = {Transportation Research Part C: Emerging Technologies},
  year    = {2022},
  volume  = {137},
  pages   = {103560},
  doi     = {10.1016/j.trc.2022.103560},
  url     = {https://doi.org/10.1016/j.trc.2022.103560}
}

@article{legal2023platooning,
  title   = {Legal Framework for Rear-End Crashes
             in Mixed-Traffic Platooning:
             A Matrix Game Approach},
  journal = {Future Transportation},
  author = {Chen, Xu and Di, Xuan},
  year    = {2023},
  volume  = {3},
  number  = {2},
  pages   = {417-428}
}

@article{wang2015game,
  title={Game theoretic approach for predictive lane-changing and car-following control},
  author={Wang, Meng and Hoogendoorn, Serge P and Daamen, Winnie and van Arem, Bart and Happee, Riender},
  journal={Transportation Research Part C: Emerging Technologies},
  volume={58},
  pages={73--92},
  year={2015},
  publisher={Elsevier}
}

@article{talebpour2015modeling,
  title={Modeling lane-changing behavior in a connected environment: A game theory approach},
  author={Talebpour, Alireza and Mahmassani, Hani S and Hamdar, Samer H},
  journal={Transportation Research Procedia},
  volume={7},
  pages={420--440},
  year={2015},
  publisher={Elsevier}
}

@article{liu2023potential,
  title={Potential game-based decision-making for autonomous driving},
  author={Liu, Mushuang and Kolmanovsky, Ilya and Tseng, H Eric and Huang, Suzhou and Filev, Dimitar and Girard, Anouck},
  journal={IEEE Transactions on Intelligent Transportation Systems},
  volume={24},
  number={8},
  pages={8014--8027},
  year={2023},
  publisher={IEEE}
}

@article{varga2024upper,
  title={On the Upper Bound of Near Potential Differential Games},
  author={Varga, Balint},
  journal={Results in Applied Mathematics},
  volume={22},
  pages={100453},
  year={2024},
  publisher={Elsevier}
}

@article{baros2025mean,
  title={Mean-Field Generalisation Bounds for Learning Controls in Stochastic Environments},
  author={Baros, Boris and Cohen, Samuel N and Reisinger, Christoph},
  journal={arXiv preprint arXiv:2508.16001},
  year={2025}
}

@article{daini2024traffic,
  title={Traffic control via fleets of connected and automated vehicles},
  author={Daini, Chiara and Delle Monache, Maria Laura and Goatin, Paola and Ferrara, Antonella},
  journal={IEEE Transactions on Intelligent Transportation Systems},
  volume={26},
  number={2},
  pages={1573-1582},
  year={2024},
  publisher={IEEE}
}

@article{bailo2018optimal,
  title={Optimal consensus control of the Cucker-Smale model},
  author={Bailo, Rafael and Bongini, Mattia and Carrillo, Jos{\'e} A and Kalise, Dante},
  journal={IFAC-PapersOnLine},
  volume={51},
  number={13},
  pages={1--6},
  year={2018},
  publisher={Elsevier}
}

@article{yan2025markov,
  title={Markov potential game construction and multi-agent reinforcement learning with applications to autonomous driving},
  author={Yan, Huiwen and Liu, Mushuang},
  journal={arXiv preprint arXiv:2503.22867},
  year={2025}
}

@article{di2025mfgreview,
  title={Mean field games for urban mobility: a review},
  author={Di, Xuan and Xu, Zhenhui and Shen, Tielong},
  journal={SCIENCE CHINA: Information Sciences},
  volume={68},
  number={11},
  pages={210201:1–210201:26},
  year={2025},
  publisher={Springer}
}

@article{festa2018mean,
  title={A mean field game approach for multi-lane traffic management},
  author={Festa, Adriano and G{\"o}ttlich, Simone},
  journal={IFAC-PapersOnLine},
  volume={51},
  number={32},
  pages={793--798},
  year={2018},
  publisher={Elsevier}
}

@article{guo2025towards,
  title={Towards an analytical framework for dynamic potential games},
  author={Guo, Xin and Zhang, Yufei},
  journal={SIAM Journal on Control and Optimization},
  volume={63},
  number={2},
  pages={1213--1242},
  year={2025},
  publisher={SIAM}
}

@article{guo2025distributed,
  title={Distributed games with jumps: An $\alpha $-potential game approach},
  author={Guo, Xin and Li, Xinyu and Zhang, Yufei},
  journal={arXiv preprint arXiv:2508.01929},
  year={2025}
}

@article{guo2025alpha,
author = {Guo, Xin and Li, Xinyu and Zhang, Yufei},
title = {An \({\alpha }\)-Potential Game Framework for \( {N}\)-Player Dynamic Games},
journal = {SIAM Journal on Control and Optimization},
volume = {63},
number = {4},
pages = {2964-3005},
year = {2025},
}

@article{zhou2017rolling,
	title={Rolling horizon stochastic optimal control strategy for ACC and CACC under uncertainty},
	author={Zhou, Yang and Ahn, Soyoung and Chitturi, Madhav and Noyce, David A},
	journal={Transportation Research Part C},
	volume={83},
	pages={61--76},
	year={2017},
	publisher={Elsevier}
}

@article{huang2006Largepopulation,
  title = {Large Population Stochastic Dynamic Games: Closed-Loop {{McKean-Vlasov}} Systems and the {{Nash}} Certainty Equivalence Principle},
  author = {Huang, Minyi and Malham{\'e}, Roland P and Caines, Peter E},
  year = {2006},
  journal = {Communications in Information and Systems},
  volume = {6},
  number = {3},
  pages = {221--252},
  publisher = {International Press of Boston}
}

@article{lasry2007mean,
  title={Mean field games},
  author={Lasry, Jean-Michel and Lions, Pierre-Louis},
  journal={Japanese Journal of Mathematics},
  volume={2},
  number={1},
  pages={229--260},
  year={2007},
  publisher={Springer}
}

@article{di2021survey,
  title={A survey on autonomous vehicle control in the era of mixed-autonomy: From physics-based to {AI}-guided driving policy learning},
  author={Di, Xuan and Shi, Rongye},
  journal={Transportation Research Part C},
  volume={125},
  pages={103008},
  year={2021},
  publisher={Elsevier}
}

@article{huang2020game,
	title={A Game-Theoretic Framework for Autonomous Vehicles Velocity Control: Bridging Microscopic Differential Games and Macroscopic Mean Field Games},
	author={Huang, Kuang and Di, Xuan and Du, Qiang and Chen, Xi},
	journal={Discrete and Continuous Dynamical Systems - Series B},
	volume={25},
	number={12},
	pages={4869-4903},
	year={2020}
}

@article{li2018nonlinear,
  title={Nonlinear consensus-based connected vehicle platoon control incorporating car-following interactions and heterogeneous time delays},
  author={Li, Yongfu and Tang, Chuancong and Peeta, Srinivas and Wang, Yibing},
  journal={IEEE Transactions on Intelligent Transportation Systems},
  volume={20},
  number={6},
  pages={2209--2219},
  year={2018},
  publisher={IEEE}
}

@article{wei2017dynamic,
	title={Dynamic programming-based multi-vehicle longitudinal trajectory optimization with simplified car following models},
	author={Wei, Yuguang and Avci, Cafer and Liu, Jiangtao and Belezamo, Baloka and Aydin, Nizamettin and Li, Pengfei Taylor and Zhou, Xuesong},
	journal={Transportation research part B: methodological},
	volume={106},
	pages={102--129},
	year={2017},
	publisher={Elsevier}
}

@article{gong2016constrained,
	title={Constrained optimization and distributed computation based car following control of a connected and autonomous vehicle platoon},
	author={Gong, Siyuan and Shen, Jinglai and Du, Lili},
	journal={Transportation Research Part B: Methodological},
	volume={94},
	pages={314--334},
	year={2016},
	publisher={Elsevier}
}

@article{yu2018human,
	title={A human-like game theory-based controller for automatic lane changing},
	author={Yu, Hongtao and Tseng, H Eric and Langari, Reza},
	journal={Transportation Research Part C: Emerging Technologies},
	volume={88},
	pages={140--158},
	year={2018},
	publisher={Elsevier}
}

@article{wang2024multi,
  title={Multi-agent DRL-controlled connected and automated vehicles in mixed traffic with time delays},
  author={Wang, Zhuwei and Xue, Yi and Liu, Lihan and Zhang, Haijun and Qu, Chunhui and Fang, Chao},
  journal={IEEE Transactions on Intelligent Transportation Systems},
  volume={25},
  number={11},
  pages={17676-17688},
  year={2024},
  publisher={IEEE}
}

@article{wang2024iterative,
  title={Iterative learning-based cooperative motion planning and decision-making for connected and autonomous vehicles coordination at on-ramps},
  author={Wang, Bowen and Gong, Xinle and Lyu, Peiyuan and Liang, Sheng},
  journal={IEEE Transactions on Intelligent Transportation Systems},
  volume={25},
  number={7},
  pages={8105--8120},
  year={2024},
  publisher={IEEE}
}

@article{jing2024decentralized,
  title={Decentralized human-like control strategy of mixed-flow multi-vehicle interactions at uncontrolled intersections: A game-theoretic approach},
  author={Jing, Dian and Yao, Enjian and Chen, Rongsheng},
  journal={Transportation Research Part C},
  volume={167},
  pages={104835},
  year={2024},
  publisher={Elsevier}
}

@article{liu2024decentralized,
  title={Decentralized platoon formation for a fleet of connected and autonomous trucks},
  author={Liu, Dahui and Eksioglu, Burak and Schmid, Matthias and Huynh, Nathan and Comert, Gurcan},
  journal={Expert Systems with Applications},
  volume={249},
  pages={123650},
  year={2024},
  publisher={Elsevier}
}

@article{hua2025multi,
  title={Multi-agent reinforcement learning for connected and automated vehicles control: Recent advancements and future prospects},
  author={Hua, Min and Qi, Xinda and Chen, Dong and Jiang, Kun and Liu, Zemin Eitan and Sun, Hongyu and Zhou, Quan and Xu, Hongming},
  volume={22},
  pages={16266-16286},
  journal={IEEE Transactions on Automation Science and Engineering},
  year={2025},
  publisher={IEEE}
}

@article{oelschlager1985law,
  title={A law of large numbers for moderately interacting diffusion processes},
  author={Oelschl{\"a}ger, Karl},
  journal={Zeitschrift f{\"u}r Wahrscheinlichkeitstheorie und verwandte Gebiete},
  volume={69},
  number={2},
  pages={279--322},
  year={1985},
  publisher={Springer}
}

@article{chen2026deep,
  title={Deep Fictitious Play-Based Potential Differential Games for Learning Human-Like Interaction at Unsignalized Intersections},
  author={Chen, Kehua and Lin, Ryan Feng and Zhang, Shucheng and Wang, Yinhai},
  journal={IEEE Transactions on Intelligent Transportation Systems},
  year={2026},
  publisher={IEEE}
}

@inproceedings{karmarkar1984new,
  title={A new polynomial-time algorithm for linear programming},
  author={Karmarkar, Narendra},
  booktitle={Proceedings of the sixteenth annual ACM symposium on Theory of computing},
  pages={302--311},
  year={1984}
}

@article{cohen2021solving,
  title={Solving linear programs in the current matrix multiplication time},
  author={Cohen, Michael B and Lee, Yin Tat and Song, Zhao},
  journal={Journal of the ACM (JACM)},
  volume={68},
  number={1},
  pages={1--39},
  year={2021},
  publisher={ACM New York, NY, USA}
}

@inproceedings{vaidya1989speeding,
  title={Speeding-up linear programming using fast matrix multiplication},
  author={Vaidya, Pravin M},
  booktitle={30th annual symposium on foundations of computer science},
  pages={332--337},
  year={1989},
  organization={IEEE}
}

@incollection{zhang2017backward,
  title={Backward stochastic differential equations},
  author={Zhang, Jianfeng},
  booktitle={Backward Stochastic Differential Equations: From Linear to Fully Nonlinear Theory},
  pages={79--99},
  year={2017},
  publisher={Springer}
}

@article{guo2025markov,
  title={Markov $\alpha$-Potential Games},
  author={Guo, Xin and Li, Xinyu and Maheshwari, Chinmay and Sastry, Shankar and Wu, Manxi},
  journal={IEEE Transactions on Automatic Control},
  year={2025},
  publisher={IEEE}
}

@article{wei2022merging,
  title={Game theoretic merging behavior control for autonomous vehicle at highway on-ramp},
  author={Wei, Chao and He, Yuanhao and Tian, Hanqing and Lv, Yanzhi},
  journal={IEEE Transactions on Intelligent Transportation Systems},
  volume={23},
  number={11},
  pages={21127--21136},
  year={2022},
  publisher={IEEE}
}

@article{han2022strategic,
  title={Strategic and tactical decision-making for cooperative vehicle platooning with organized behavior on multi-lane highways},
  author={Han, Xu and Xu, Runsheng and Xia, Xin and Sathyan, Anoop and Guo, Yi and Bujanovi{\'c}, Pavle and Leslie, Ed and Goli, Mohammad and Ma, Jiaqi},
  journal={Transportation Research Part C: Emerging Technologies},
  volume={145},
  pages={103952},
  year={2022},
  publisher={Elsevier}
}

@article{liu2025cooperative,
  title={Cooperative control method for connected and automated vehicle platoon based on arbitrary time headway switched system},
  author={Liu, Gongzhe and Zheng, Nan and Wang, Hao},
  journal={Transportation Research Part C: Emerging Technologies},
  volume={180},
  pages={105353},
  year={2025},
  publisher={Elsevier}
}

@article{jond2022differential,
  title={Differential game-based optimal control of autonomous vehicle convoy},
  author={Jond, Hossein B and Plato{\v{s}}, Jan},
  journal={IEEE Transactions on Intelligent Transportation Systems},
  volume={24},
  number={3},
  pages={2903--2919},
  year={2022},
  publisher={IEEE}
}

@article{liu2024cooperative,
  title={Cooperative decision-making for cavs at unsignalized intersections: A marl approach with attention and hierarchical game priors},
  author={Liu, Jiaqi and Hang, Peng and Na, Xiaoxiang and Huang, Chao and Sun, Jian},
  journal={IEEE Transactions on Intelligent Transportation Systems},
  volume={26},
  number={1},
  pages={443--456},
  year={2024},
  publisher={IEEE}
}

@article{guo2024heuristic,
  title={Heuristic-based multi-agent deep reinforcement learning approach for coordinating connected and automated vehicles at non-signalized intersection},
  author={Guo, Zihan and Wu, Yan and Wang, Lifang and Zhang, Junzhi},
  journal={IEEE Transactions on Intelligent Transportation Systems},
  volume={25},
  number={11},
  pages={16235--16248},
  year={2024},
  publisher={IEEE}
}

@article{chen2025game,
  title={A Game-Theoretical Framework for Safe Decision Making and Control of Mixed Autonomy Vehicles},
  author={Chen, Mingyang and Li, Bingbing and Zhang, Sunan and Zhang, Hao and Zhuang, Weichao and Yin, Guodong and Chen, Boli},
  journal={IEEE Transactions on Intelligent Transportation Systems},
  volume={27},
  number={1},
  pages={1338--1351},
  year={2025},
  publisher={IEEE}
}

@article{heshami2024towards,
  title={Towards Self-Organizing connected and autonomous Vehicles: A coalitional game theory approach for cooperative Lane-Changing decisions},
  author={Heshami, Seiran and Kattan, Lina},
  journal={Transportation Research Part C: Emerging Technologies},
  volume={166},
  pages={104789},
  year={2024},
  publisher={Elsevier}
}

 \end{document}